\documentclass[10pt,a4paper, twoside]{article}
\usepackage[notref,notcite]{showkeys}

\usepackage{listings}
\usepackage{amsmath,amssymb,amsthm,amsfonts,mathrsfs,amscd,environ}
\usepackage{dsfont}
\usepackage{latexsym,enumerate,color,geometry,extarrows}
\usepackage{xcolor}
\usepackage{verbatim,fancyhdr,enumitem}
\usepackage{hyperref}
\hypersetup{
hidelinks,
hypertex=true,
colorlinks=true,
linkcolor=red,
citecolor=blue,
urlcolor=blue
}
 \NewEnviron{ews}{%
\begin{equation}\begin{split}
  \BODY
\end{split}\end{equation}
}

\NewEnviron{ews*}{%
\begin{equation*}\begin{split}
  \BODY
\end{split}\end{equation*}
}

\def\beg{\begin}
\def\bequ{\begin{equation}}
\def\enqu{\end{equation}}
\def\bes{\begin{split}}
\def\ens{\end{split}}
\def\bews{\begin{ews}}
\def\beqn{\begin{eqnarray}}
\def\enqn{\end{eqnarray}}
\def\beq*{\begin{equation*}}
\def\enq*{\end{equation*}}
\def\bqn*{\begin{eqnarray*}}
\def\eqn*{\end{eqnarray*}}
\def\bary{\begin{array}}
\def\eary{\end{array}}
\def\bpma{\begin{pmatrix}}
\def\epma{\end{pmatrix}}
\def\bvma{\begin{Vmatrix}}
\def\evma{\end{Vmatrix}}

 \newtheorem{thm}{Theorem}[section]
 
 \newtheorem{cor}[thm]{Corollary}

 \numberwithin{equation}{section}

\def\al{\alpha}
\def\be{\beta}
\def\ga{\gamma}
\def\de{\delta}

\def\ze{\zeta}
\def\et{\eta}
\def\th{\theta}

\def\ka{\kappa}
\def\la{\lambda}
\def\rh{\rho}

\def\si{\sigma}

\def\ps{\psi}

\def\Ga{\Gamma}

\def\De{\Delta}

\def\Ph{\Phi}
\def\Ps{\Psi}
\def\Om{\Omega}

\def\R{\mathbb R}
\def\P{\mathbb P}
\def\E{\mathbb E}

\def\N{\mathbb N}

\def\sF{\mathscr F}

\def\sB{\mathscr B}

\def\sL{\mathscr L}

\def\sP{\mathscr P}

\def\cL{\mathcal L}
\def\cH{\mathcal H}
\def\cW{\mathcal W}

\def\cT{\mathcal T}
\def\cM{\mathcal M}

\def\d{\mathrm{d}}

\def\ff{\frac}
\def\ra{\rightarrow}
\def\nn{\nabla}
\def\pp{\partial}
\def\<{\langle}
\def\>{\rangle}
\def\sq{\sqrt}
\def\tld{\tilde}
\def\we{\wedge}
\def\1{\mathds{1}}

\def\trac{\mathrm{tr}}

\def\div{\mathrm{div}}

\def\supp{\displaystyle{\mathrm{supp}}}

\allowdisplaybreaks

\title{{\bf A Local Bifurcation Theorem  for McKean-Vlasov Diffusions}
}

\author{
{\bf Shao-Qin Zhang }\\
\footnotesize{School of Statistics and Mathematics, Central University of Finance and Economics, Beijing 100081, China}\\
\footnotesize{Email: zhangsq@cufe.edu.cn}\\
}

\begin{document}

\maketitle

\begin{abstract}
We establish an existence result of a solution to  a class of probability measure-valued equations, whose solutions can be associated with stationary distributions of many McKean-Vlasov diffusions with gradient-type drifts.  Coefficients of the probability measure-valued equation may be discontinuous in the weak topology and the total variation norm. Owing to that the bifurcation point of the probability measure-valued equation is relevant to the phase transition point of the associated  McKean-Vlasov diffusion, we establish a local Krasnosel’skii bifurcation theorem.  Regularized determinant for the Hilbert-Schmidt operator is used to derive our criteria for the bifurcation point.  Concrete examples, including the granular media equation and the Vlasov-Fokker-Planck equation with quadratic interaction, are given to illustrate our results.

\end{abstract}\noindent


\noindent


\vskip 2cm

\section{Introduction}

By passing to the mean field limit for a system of interacting diffusions, a stochastic differential equation (SDE) whose coefficients depend on the own law of the solution was introduced by McKean in \cite{McK}. This SDE is also called distribution dependent SDE or mean-field SDE, see e.g. \cite{BLPR,RZ,Wan18}. The associated empirical measure of the interacting diffusions converges in the weak sense  to a probability measure with density, which is called the propagation of chaos property, and the density  satisfies a nonlinear parabolic partial differential equation called McKean-Vlasov equation in the literature, see e.g. \cite{ChD,Szn}. The existence of several stationary distributions to McKean-Vlasov SDEs is referred to phase transition.  \cite{Daw} established for the first time the phase transition for the equation with a particular double-well confinement and Curie-Weiss interaction on the line. Precisely,  stationary distributions of the following SDE was investigated in \cite{Daw}: 
\beg{align}\label{Daw}
\d X_t=-(X_t^3-X_t)\d t-\be(X_t-\E X_t)\d t+\si\d B_t,
\end{align}
where $B_t$ is a one dimensional Brownian motion in the probability space $(\Om,\sF,\P)$, $\E$ is the expectation with respect to $\P$,  and $\si,\be$ are positive constants. The stationary distributions of \eqref{Daw} can be obtained by solving the following equation
\beg{equation}\label{Daw-fixp}
\nu(\d x)=\ff {\exp\left\{-\ff 2 {\si^2}\left(\ff {x^4} 4-\ff {x^2} 2\right)-\ff {\be} {\si^2}\int_{\R}(x-z)^2 \nu(\d z)\right\}} {\int_{\R} \exp\left\{-\ff 2 {\si^2}\left(\ff {x^4} 4-\ff {x^2} 2\right)-\ff {\be} {\si^2}\int_{\R}(x-z)^2 \nu(\d z)\right\}\d x}\d x.
\end{equation}
The confinement potential $\ff {x^4} 4-\ff {x^2} 2$ has two minima, which is known as double-well potential. It is  proved in \cite{Daw} that for fixed $\be>0$, there exists $\si_c>0$ so that \eqref{Daw} has a unique stationary distribution if $\si>\si_c$  and has three stationary distributions if $0<\si<\si_c$. Beside \cite{Daw},  phase transitions for McKean-Vlasov SDEs are studied by many works, e.g. \cite{Tam} provided a criteria of the phase transition for equations on the whole space; equations with multi-wells confinement were investigated extensively  by Tugaut et  al. in \cite{DuTu,HT10a,Tug10,Tug13,Tug14a,Tug14b}; quantitative results on  phase transitions  of McKean-Vlasov diffusions on the torus were provided in \cite{CGPS,ChPa}; the relation between phase transition and functional inequality was investigated in  \cite{DGPS}; non-uniqueness of stationary distributions for general distribution dependent SDEs was discussed in \cite{ZSQ}. Phase transitions of nonlinear Markov jump processes were studied in \cite{Chen,FenZ}.

Bifurcation theory has been used to analyse the phase transition. For instance, \cite{CGPS} showed that as the intensity of the diffusion term  or the intensity of the interaction potential crosses a critical point,  new stationary distributions branch out from the uniform distribution, which is a homogeneous steady state of McKean-Vlasov diffusion on torus without confinement potential.  Bifurcation analysis was also given by \cite{Tam} for McKean-Vlasov SDEs on the whole space with odd interaction potentials.  However, the assumption that the interaction potential is odd is unphysical, and excludes the model in \cite{Daw}. The Crandall-Rabinowitz theorem used in \cite{CGPS,Tam} requires that the Fredholm operator induced by the interaction potential should have one dimensional null space.

Denote by $\sP(\R^d)$ the space of all probability measures on $\R^d$. In this paper, we analyse solutions of equations on $\sP(\R^d)$ of the following form 
\beg{align}\label{fix-p}
\mu(\d x)=\ff {\exp\{-V_0(x)-  V(x,\mu) \}} {\int_{\R^d} \exp\{- V_0(x)-  V(x,\mu)\}\d x}\d x,
\end{align}
where $\mu$ is a probability measure, $V_0:\R^d\mapsto\R$ is measurable and $V$ is a function on $\R^d\otimes\sP(\R^d)$. We first establish an existence result of a solution to \eqref{fix-p}. When $V(x,\mu)=\int_{\R^d}V(x,y)\mu(\d y)$ for some measurable function $V:\R^d\otimes\R^d\mapsto \R$, we introduce a parameter  $\al$ to \eqref{fix-p} to characterize the intensity of the diffusion term  or the intensity of the interaction potential (as $\si$ or $\be$ in \eqref{Daw-fixp}). Precisely, we consider 
\beg{align}\label{eq-bif}
\mu(\d x)=\ff {\exp\{-\th(\al)V_0(x)-  \al \int_{\R^d}V(x,y)\mu(\d y) \}} {\int_{\R^d} \exp\{-\th(\al) V_0(x)- \al \int_{\R^d}V(x,y)\mu(\d y)\}\d x}\d x,
\end{align}
where $\th$ is a continuously differentiable function, and establish a local Krasnosel’skii bifurcation theorem (see e.g. \cite{Kie,Kra}) to \eqref{eq-bif}. This local bifurcation theorem allows the interaction potential to induce a Fredholm operator with multidimensional null space. 

Equation  \eqref{eq-bif} (or \eqref{fix-p}) generalises \eqref{Daw-fixp}. Besides \eqref{Daw}, stationary distributions and  phase transitions of McKean-Vlasov SDEs with gradient-type drifts can be obtained  by analyzing \eqref{eq-bif}. To be specific, when $V_0$ and $V$ are continuous differentiable functions,  solutions of \eqref{eq-bif} can be associated with stationary distributions of  McKean-Vlasov SDEs on $\R^d$ of the following form
\beg{equation}\label{SDE-1.4}
\d X_t=b(X_t,\sL_{X_t})\d t+\sq{\ff 2 {\al}}  \si(X_t,\sL_{X_t})\d B_t,
\end{equation}
where $\sL_{X_t}$ denote the distribution of $X_t$,  $\si:\R^d\times\sP(\R^d)\mapsto\R^d\otimes\R^d$ so that for $\mu\in\sP(\R^d)$, $\si(\cdot,\mu)$ is continuous differentiable, 
\beg{align}\label{bbb}
b(x,\mu)&:=-\sum_{i,j=1}^d\left(\sum_{k=1}^d\si_{ik}\si_{jk}\right)(x,\mu)\left(\ff {\th(\al)} {\al}\pp_j V_0(x)+\int_{\R^d}\pp_j V(\cdot,y)(x)\mu(\d y)\right)e_i\nonumber\\
&\quad\, +\al^{-1}\sum_{i,j,k=1}^d\pp_j(\si_{ik}(\cdot,\mu)\si_{jk}(\cdot,\mu))(x)e_i,~x\in\R^d,\mu\in\sP(\R^d),
\end{align}
and $\{e_i\}_{i=1}^d$ is the standard orthonormal basis of $\R^d$. In this setting, the parameter $\al$ can be used to characerize the absolute temperature when \eqref{SDE-1.4} is used to model systems in statistical physics, see e.g.  \cite[Page 38]{Daw} or \cite[Page 43]{Mcc}. \eqref{Daw-fixp} and \eqref{Daw} can be derived from \eqref{eq-bif} and \eqref{SDE-1.4} by setting  $d=1$, $V(x,y)=(x-y)^2$, $V_0(x)=\ff {x^4} 4-\ff {x^2} 2$, $\si(x,\mu)\equiv\sq{\ff {\be} 2}$, and $\al=\ff {\be} {\si^2}$. Moreover,  \eqref{eq-bif} can also be used to analyze stationary distributions of  degenerate SDEs which can not be covered by \eqref{SDE-1.4} with the drift term satisfying  \eqref{bbb}. For instance,   
\beg{exa}\label{exa11}
For the degenerate system on $\R^2$
\begin{equation}\label{Hamil-deg}
\begin{cases}
\d X_t =  Y_t\d t\\
\d Y_t  = -(X_t^3-X_t)\d t-\be\int_{\R}(X_t-z)\sL_{X_t}(\d z)\d t-Y_t\d t+\si \d B_t,
\end{cases}
\end{equation}
where $B_t$ is the one dimensional Brownian motion and $\be,\si$ are positive constants, stationary distributions of this system can be obtain by solving the following equation
\beg{equation}\label{Ham-exp} 
\nu(\d x,\d y)=\ff {\exp\left\{-\ff 2 {\si^2}\left(\ff {y^2} 2+\ff {x^4} 4-\ff {x^2} 2\right)-\ff {\be} {\si^2}\int_{\R}(x-z)^2 \nu_1(\d z)\right\}\d x\d y} {\int_{\R} \exp\left\{-\ff 2 {\si^2}\left(\ff {y^2} 2+\ff {x^4} 4-\ff {x^2} 2\right)-\ff {\be} {\si^2}\int_{\R}(x-z)^2 \nu_1(\d z)\right\}\d x\d y},
\end{equation}
where $\nu_1(\d z)=\int_{\R}\nu(\d z,\d w)$ is a marginal of $\nu$. We set accordingly that $V_0(x,y)=\ff {y^2} 2+\ff {x^4} 4-\ff {x^2} 2$, $V(x,y,z,w)=(x-z)^2$, $\al=\ff {\be} {\si^2}$ and $\th(\al)=\ff 2 {\be}\al$. Fixing $\be$,   critical points of phase transitions for \eqref{Ham-exp} or \eqref{Hamil-deg} as $\si$ varying are related to bifurcation  points of the parameter $\al$.
\end{exa}
\noindent This example will be revisited in Section 2 and Section 3. For more examples, one can see \cite{CGPS,Daw,DuTu,Tam,Tug14a} or examples in Section 2 and Section 3.
 
This paper is structured as follows. In Section 2, we prove the existence of a solution for \eqref{fix-p}, see Theorem \ref{exis-thm1}. This theorem is established by using the Lyapunov condition and the Schauder fixed point theorem. Our assumptions allow that $V(\cdot,\mu)$ is in some first order Sobolev space and $V(x,\cdot)$ is continuous w.r.t.  some weighted variation distance, see {\bf Assumption (H)} below.   In Section 3, a local bifurcation theorem is established, see Theorem \ref{thm-bif}.  By using the regularized determinant for the Hilbert-Schmidt operator (see e.g. \cite{Simon}), we give a criteria of the  bifurcation point, which is based on the algebraic multiplicity of an eigenvalue for the integral operator induced by the kernel $V(x,y)$.   

{\bf Notation:} The following notations are used in the sequel.\\
$\bullet$ We denote by $L^p$ (resp. $L^p(\mu)$) the space of functions for which the $p$-th power of the absolute value is Lebesgue integrable (resp. integrable w.r.t. the measure $\mu$), and $W^{k,p}$ (resp. $W^{1,p}_{loc}$) the $k$ order  (resp. local) Sobolev space on $\R^d$; $C_0$ (resp.  $C_0^\infty$) the space of all the continuous (resp. smooth) functions with compact support on $\R^d$; $C^k$ the $k$ continuously differentiable functions on $\R^d$. For a probability measure $\mu$, we denote 
\beg{align*}
&\cW^{k,p}_{q, \mu}=\{f\in W^{k,p}_{loc}~|~\nn f,\cdots,\nn^k f\in L^q(\mu)\}.
\end{align*}
We use $\cW^{k,p}_{ \mu}$  to denote $\cW^{k,p}_{p, \mu}$.  We denote by $\sL(L^2( \mu))$ and $\sL_{HS}(L^2( \mu))$ the space of all bounded operators and the space of all Hilbert-Schmidt operators on $L^2( \mu)$ respectively. \\
$\bullet$ For measurable function $f$ on $\R^d$, we define for $p,q\in [1,+\infty]$  
\beg{align*}
\|f\|_{L^p_xL^{q}_y}&=\left(\int_{\R^d}\|f(\cdot,y)\|_{L^{p}(\bar\mu)}^{q}\bar\mu(\d y)\right)^{\ff 1 {q}},\\
\|f\|_{L^p_yL^{q}_x}&=\left(\int_{\R^d}\|f(x,\cdot)\|_{L^{p}(\bar\mu)}^{q}\bar\mu(\d x)\right)^{\ff 1 {q}},
\end{align*}
where $\bar\mu$ is a probability measure defined in \eqref{barmu}.\\
Let $\chi $ be a decreasing and  continuously differentiable function on $[0,+\infty)$ such that $\1_{[0\leq r\leq 1]}\leq \chi(r)\leq \1_{[0\leq r\leq 2]}$ and $|\chi'(r)|\leq 2$. Denote by $\ze_n(x)=\chi(|x|/n)$.\\
We denote by $B_N$  the open ball with  radius $N$ and centre at $0$.\\
$\bullet$ We denote by $\sP(\R^d)$ the space of all probability measures on $\R^d$. For any measurable function $V\geq 1$,   
$$\sP_{V}:=\{\mu\in\sP(\R^d)~|~\|\mu\|_{V}:=\mu(V)<\infty\}, $$
endowed with the weight total variance distance:
$$\|\mu-\nu\|_{V}=\sup_{|f|\leq V} \left|\mu(f)-\nu(f)\right|,~\mu,\nu\in\sP_V.$$
For a probability measure $\mu$ and a measurable function $f$, we denote by $f\mu$ the sign measure $(f\mu)(\d x)=f(x)\mu(\d x)$. \\
$\bullet$ For a function on $\R^d\times\sP(\R^d)$, saying  $f$, we denote by $f(\mu)$ the function $f(\cdot,\mu)$ with $\mu$ given, and $\nn f(x,\mu)$ the gradient w.r.t. the variable on $\R^d$ at $(x,\mu)$. For a function $f:L^2(\mu_1)\times J\mapsto L^2(\mu_2)$, where $\mu_1,\mu_2$ are probability measures and $J$ is an interval in $\R$, letting $(\rh,\al)\in L^2(\mu_1)\times J$, we denote by $f(x;\rh,\al)$ the value of the function $f(\rh,\al)$ at $x\in\R^d$, and $\nn f(\rh,\al)$ the Fr\'echet derivative and $\pp_\al f(\rh,\al)$ the derivative of the second variable.

\section{Existence}


In this section, we investigate the existence of a solution to \eqref{fix-p}. To this aim, we choose a reference probability measure 
\beg{align}\label{barmu}
\bar\mu(\d x):=\ff {e^{-\bar V(x)}}{\int_{\R^d} e^{-\bar V(x)}\d x}\d x, 
\end{align}
and reformulate \eqref{fix-p} into another form:
\[\mu(\d x)=\ff {\ps(x,\mu)} {\bar\mu(\ps(\mu))}\bar\mu(\d x),\]
where 
\[\ps(x,\mu)=\exp\left\{-V_0(x)-V(x,\mu)+\bar V(x) \right\},\]
and potentials $V_0,V,\bar V$ satisfy following assumptions
\beg{description}
\item {\bf Assumption (H)} 
\item (H1)   The potentials $V_0$ and $\bar V$ are measurable functions such that  $e^{-V_0},e^{-\bar V}\in L^1$, and there exist $p\in (d,+\infty]$ and $q\geq 1$ such that $V_0,\bar V\in \cW^{1,p}_{q,\bar\mu}$. 

\item (H2) There is a measurable function $W_0\geq 1$ such that $W_0\in L^1(\bar\mu)$, $V:\R^d\times \sP_{W_0}\mapsto \R$ is measurable and for all $\mu\in \sP_{W_0}$, $V(\cdot,\mu)\in W^{1,p}_{loc}$. 
There exist nonnegative functions $F_0,F_1,F_2,F_3$ such that $F_0\in L^{\infty}_{loc}$,  $F_2\in L^q(\bar\mu)\bigcap L^p_{loc}$, $F_1,F_3$ are increasing on $[0,+\infty)$ with $\lim_{r\ra 0^+}F_1(r)=0$, and
\beg{align}
|V(x,\mu)-V(x,\nu)|&\leq F_0(x)F_1(\|\mu-\nu\|_{W_0}),\label{V2FF}\\
|V(x,\bar\mu)|&\leq C(F_0(x)+1),\label{bmu-F}\\
|\nn V(x,\mu)|&\leq F_2(x)F_3(\|\mu\|_{W_0}),~\mu,\nu\in \sP_{W_0}.\label{nnV2}
\end{align}
\item (H3) There is a nonnegative and increasing function $F_4$ on $[0,+\infty)$ such that 
\beg{align}\label{V-F3}
-V_0(x)+\be F_0(x)\leq -\bar V(x)+F_4(\be),~\be\geq 0.
\end{align}
\end{description}
Under the assumption {\bf (H)}, we can prove that $\ps(\mu)\in L^\infty$, see Lemma \ref{ps(mu)} below. Then {(H1)} implies that $\ps(\mu)\in L^1(\bar\mu)$. Let
\beg{align}\label{hatT}
\hat\cT(x,\mu)=\ff {\ps(x,\mu)} {\bar\mu(\ps(\mu))}.
\end{align}
We also denote by $\hat \cT$ the mapping $\hat\cT:~\mu\mapsto \hat\cT(\cdot,\mu)$ when there is no confusion. For every $0\leq f\in L^1(\bar\mu)$ with $\bar\mu(f)=1$, we define
$$\mathscr{I}: f \mapsto \mathscr{I}(f)\equiv f\bar\mu\in\sP(\R^d).$$
For a fixed point of $\hat\cT\circ \mathscr{I}$, saying $\rh$, the  probability measure $\rh\bar\mu$ satisfies \eqref{fix-p}. Hence, we investigate the fixed point of $\hat\cT\circ \mathscr{I}$ instead of \eqref{fix-p}.

Giving $\mu\in\sP_{W_0}$, we introduce the following differential operator:
\beg{align}\label{Lmu0}
L_{\mu}g &:= \De g-\< \nn (V_0+V(\mu)),\nn g\>\nonumber\\
&\,=\De g+\<\nn \log(\ps(\mu)e^{-\bar V}),\nn g\>, g\in C_0^\infty.
\end{align}
Due to (H1) and (H2), $V_0,V(\mu)\in \cW^{1,p}_{q,\bar\mu}$. Thus $L_\mu$ is well-defined. We assume that $L_\mu$ satisfies the following Lyapunov condition. 
\beg{description}
\item {\bf Assumption (W)} 

\item  (W1) There is a measurable function $W\geq 1$ such that $\displaystyle\lim_{|x|\ra +\infty}W(x)=+\infty$ and
\beg{align}\label{W0W}
\sup_{x\in\R^d}\ff {W_0(x)} {W(x)}<\infty,\quad &\varlimsup_{|x|\ra+\infty}\ff {W_0(x)} {W(x)} =0.
\end{align}
\item  (W2) There exist a positive measurable function $W_1\in W^{2,1}_{loc}$ and strictly increasing functions $G_1,G_2$ on $[0,+\infty)$ such that  $G_2$ is convex and  
\beg{align}\label{G12}
\varlimsup_{r\ra+\infty}\ff {G_1(r)} {G_2(r)}<1&,\\
L_{\mu}W_1\leq G_1(\|\mu\|_W)-G_2(W)&,~\mu\in\sP_W.\label{LYP}
\end{align}
\end{description}
The condition (W1) implies that $\sP_W\subset\sP_{W_0}$. Thus, $L_\mu$ is well-defined for $\mu\in\sP_W$. We have the following theorem on the fixed point of $\hat\cT\circ \mathscr{I}$.

\beg{thm}\label{exis-thm1}
Assume  that {\bf(H)} holds with $F_0\in L^1(W_0\bar\mu)$, and {\bf(W)} holds with $W_1\in \cW^{2,p_1}_{\bar\mu}$ for some $p_1\geq \ff { q} { q-1}$. Then  $\hat\cT\circ \mathscr{I}$ has a fixed point in $\cW^{1,p}_{q,\bar\mu}\cap L^\infty\cap L^1(W\bar\mu)$. 
\end{thm}

To illustrate this theorem, we give the following corollaries and examples.  The first corollary can be used to investigate the existences of stationary distributions for the granular media equation, see e.g. \cite{CMV,Wan23}. 
\beg{cor}\label{Exa-Gr}
Consider the following equation:
\beg{equation}\label{exa-granular}
\mu(\d x)=\ff {\exp\{-V_0(x)+\int_{\R^d}H(x-y)\mu(\d y)\}\d x} {\int_{\R^d}\exp\{-V_0(x)+\int_{\R^d}H(x-y)\mu(\d y)\}\d x},
\end{equation}
where  $V_0,H\in C^1(\R^d)$, $V_0,\nn V_0$ have polynomial growth: there is $\ga_0>0$ such that
\beg{equation}\label{V-poly}
\varlimsup_{|x|\ra +\infty}\ff {|V_0(x)|+|\nn V_0(x)|} {(1+|x|)^{\ga_0}}=0,
\end{equation}
and there exist positive constants $C_i,i=0,\cdots,5$, $\ga_i,i=1,2,3,4$ with $\ga_3\geq \ga_4$
\beg{equation}\label{ga1234}
\ga_1>\ga_2\vee\ga_3\vee (\ga_4+1)
\end{equation}
such  that  for all $x\in\R^d$
\beg{align}
V_0(x)&\geq C_0(1+|x|)^{\ga_1}-C_1,\label{exa-V0}\\
\<\nn V_0(x),x\>&\geq C_2(1+|x|)^{\ga_1}-C_3,\label{exa-nnV0}\\
|H(x-y_1)-H(x-y_2)|&\leq C_4(1+|x|)^{\ga_2}((1+|y_1|)^{\ga_3}+(1+|y_2|)^{\ga_3}),\label{HH}\\
|\nn H(x)|&\leq C_5(1+|x|)^{\ga_4}.\label{exa-nnH}
\end{align}
Let 
\[\bar\mu(\d x)= \ff {e^{-\ff {C_0} 2(1+|x|)^{\ga_1}}\d x} {\int_{\R^d}e^{-\ff {C_0} 2(1+|x|)^{\ga_1}}\d x},~W(x)=(1+|x|)^{\ga_1}.\]
Then for any $q\in [1,+\infty)$, \eqref{exa-granular} has a solution $\mu$ with  $\ff {\d \mu} {\d \bar\mu}\in L^\infty\cap L^1(W\bar\mu)\cap   \cW^{1,\infty}_{q,\bar\mu}$.
\end{cor}
\beg{proof}
We first check {\bf (H)}.  Let $\bar V(x)=\ff {C_0} 2(1+|x|)^{\ga_1}$, $W_0(x)=(1+|x|)^{\ga_3}$. Then $V_0,\bar V\in \cW^{1,\infty}_{q,\bar\mu}$ for any $q\geq 1$, and $W_0\in L^1(\bar\mu)$. Thus $(H1)$ holds. For any $\mu_1,\mu_2\in \sP_{W_0}$, let $\pi$ be the Wasserstein coupling of $\mu_1,\mu_2$, i.e.
\[\pi(\d y_1,\d y_2)=(\mu_1\we\mu_2)(\d y_1)\de_{y_1}(\d y_2)+\ff {(\mu_1-\mu_2)^+(\d y_1)(\mu_1-\mu_2)^-(\d y_2)} {(\mu_1-\mu_2)^-(\R^d)}.\] 
Then it follows from \eqref{HH} that 
\beg{align*}
&\left|\mu_1(H(x-\cdot))-\mu_2(H(x-\cdot))\right|\\
&=\left|\int_{\R^d\times\R^d}(H(x-y_1)-H(x-y_2))\pi(\d y_1,\d y_2)\right|\\
&=\left|\int_{\R^d\times\R^d}(H(x-y_1)-H(x-y_2))\ff {(\mu_1-\mu_2)^+(\d y_1)(\mu_1-\mu_2)^-(\d y_2)} {(\mu_1-\mu_2)^-(\R^d)}\right|\\
&\leq C_4(1+|x|)^{\ga_2} \int_{\R^d\times\R^d}(W_0(y_1)+W_0(y_2))\ff {(\mu_1-\mu_2)^+(\d y_1)(\mu_1-\mu_2)^-(\d y_2)} {(\mu_1-\mu_2)^-(\R^d)}\\
&=C_4(1+|x|)^{\ga_2} \left((\mu_1-\mu_2)^+(W_0)+(\mu_1-\mu_2)^-(W_0)\right)\\
&=C_4(1+|x|)^{\ga_2}\|\mu_1-\mu_2\|_{W_0}.
\end{align*}
Due to  \eqref{exa-nnH},    there is a contant $C>0$ such that
\[|H(x)|\leq |H(0)|+C_5 (1+|x|)^{\ga_4}|x|\leq C(1+|x|)^{\ga_4+1}.\]
Combining this with \eqref{HH}, we find that
\begin{align*}
|\bar\mu(H(x-\cdot))|&\leq |H(x)|+C_4(1+|x|)^{\ga_2}(1+\bar\mu((1+|\cdot|)^{\ga_3}))\\
&\leq  \left(C+C_4(1+\|\bar\mu\|_{W_0})\right)(1+|x|)^{(\ga_4+1)\vee \ga_2}.
\end{align*}
It follows from the dominated convergence theorem, $\ga_3\geq \ga_4$ and \eqref{exa-nnH} that  
\beg{align*}
\left|\left(\nn\int_{\R^d}H(\cdot-y)\mu(\d y)\right)(x)\right|&=\left|\mu((\nn H)(x-\cdot))\right|\leq C_5\mu((1+|x-\cdot|)^{\ga_4})\\
&\leq C_5(1+|x|)^{\ga_4}\mu((1+|\cdot|)^{\ga_4})\leq C_5(1+|x|)^{\ga_4}\|\mu\|_{W_0}^{\ff {\ga_4} {\ga_3}}.
\end{align*}
We set $F_0(x)=(1+|x|)^{(\ga_4+1)\vee \ga_2}$, $F_1(r)=C_4r$, $F_2(x)=(1+|x|)^{\ga_4}$, $F_3(r)=C_5r^{\ff {\ga_4} {\ga_3}}$. Then by the H\"older inequality and $\ga_1>(\ga_4+1)\vee \ga_2$, there exists a constant $C>0$ such that
\begin{align*}
-V_0(x)+\be F_0(x)&\leq -\ff {C_0} 2(1+|x|)^{\ga_1}+C_1+C\be^{\ff {\ga_1} {\ga_1-(\ga_4+1)\vee \ga_2} },~\be>0.
\end{align*}
Hence, {\bf (H)} holds.

Set $W_1(x)=|x|^2$, $W(x)=(1+|x|)^{\ga_1}$, and  
\beg{align}\label{Lmu00}
L_\mu g(x)=\De g(x) -\left(\nn V_0(x)+\mu((\nn H)(x-\cdot))\right)\cdot(\nn g)(x)~,~\mu\in\sP_W,~g\in C^2.
\end{align}
Then $W_1\in W^{2,\infty}\cap \cW^{2,p_1}_{\bar\mu}$ for any $p_1\geq 1$. By using the H\"older inequality,  there exist positive constants $\tld C_1,\tld C_2,\tld C_3$ such that
\beg{align*}
\left(L_\mu W_1\right)(x)&=2d-2\<\nn V_0(x),x\>-2\<\mu((\nn H)(x,\cdot)),x\>\\
&\leq 2d-2C_2(1+|x|)^{\ga_1}+2C_3+2C_5|x| \mu((1+|x-\cdot|)^{\ga_4})\\
&\leq 2d-2C_2(1+|x|)^{\ga_1}+2C_3+2^{(\ga_4-1)^++1}C_5|x|\left((1+|x|)^{\ga_4}+\mu((1+|\cdot|)^{\ga_4})\right)\\
&\leq -\tld C_1(1+|x|)^{\ga_1}+\tld C_2+\tld C_3 \mu((1+|\cdot|)^{\ga_4})^{\ff {\ga_1} {\ga_1-1} }\\
&\leq -\tld C_1(1+|x|)^{\ga_1}+\tld C_2+\tld C_3 \mu((1+|\cdot|)^{\ga_1})^{\ff {\ga_4} {\ga_1-1} }.
\end{align*}
Thus $G_1(r)=\tld C_3 r^{\ff {\ga_4} {\ga_1-1}}$, $G_2(r)=\tld C_1r-\tld C_2$. Then \eqref{G12} holds due to $\ga_1>\ga_4+1$, and \eqref{W0W} holds since $\ga_1>\ga_3$. Hence, {\bf (W)} holds.

Therefore, for any $q\in [1,+\infty)$, Theorem \ref{exis-thm1} implies that \eqref{exa-granular} has a fixed point $\mu$ and $\ff {\d \mu} {\d \bar\mu}\in L^\infty\cap L^1(W\bar\mu)\cap   \cW^{1,\infty}_{q,\bar\mu}$.

\end{proof}

It is clear that  stationary distributions of \eqref{Daw} satisfy \eqref{Daw-fixp}, and  \eqref{Daw-fixp} satisfies all conditions of Corollary \ref{Exa-Gr} with $d=1$, $V_0(x)=\ff 2 {\si^2}(\ff {x^4} 4-\ff {x^2} 2)$, $H(x)=x^2$ and $\ga_1=4$, $\ga_2=\ga_4=1$, $\ga_3=2$. Solutions of \eqref{exa-granular} can be associated with stationary distributions of  the following SDE
\beg{equation}\label{ad-V0HB}
\d X_t=\nn V_0(X_t)\d t+\int_{\R^d}\nn H(x-z)\sL_{X_t}(\d z)\d t+\sq 2\d B_t.
\end{equation}
Indeed, if $\mu$ is a solution of \eqref{exa-granular}, then the integration by part formula implies that 
\[\mu(L_{\mu}f)=0,~f\in C^2_0,\]
where $L_{\mu}$ is defined by \eqref{Lmu00}. If the strong well-posedness holds for \eqref{ad-V0HB} and the following  equation
\beg{equation}\label{ad-SXX}
\d X_t=\nn V_0(X_t)\d t+\int_{\R^d}\nn H(x-z)\mu(\d z)\d t+\sq 2\d B_t,
\end{equation}
then stationary distributions of \eqref{ad-SXX} are of the form of the right-hand side of \eqref{exa-granular}, and it follows from \cite[Proposition 1.5.5]{BKR} that $\sL_{X_t}\equiv \mu$ if $\sL_{X_0}=\mu$.  
Usually, if $\nn V_0$ has a strongly  dissipative leading term and $\nn H$ is ``weak" comparing with the dissipative term, then Corollary \ref{Exa-Gr} is applicable.  To illustrate this corollary, we  give another concrete example.
\beg{exa}
Let $A_0\in\R^d\otimes\R^d$ so that $A_0$ is a positive definite matrix with $A_0\geq a_0>0$ for some constant $a_0$, and let $v_0\in\R^d$, 
\[V_0(x)=|A_0x|^4-|x|^2\sin(|x|),~~H(x)=(1+|x|^2)\sin(\<v_0,x\>).\] 
Then \eqref{V-poly}-\eqref{exa-nnH} hold, and in this case, solutions of \eqref{exa-granular} are associated with stationary distributions of the following SDE
\beg{align*}
&\d X_t=\sq 2\d B_t-\left(|A_0X_t|^2A_0^2X_t-2\sin(|X_t|)X_t-|X_t|\cos(|X_t|)X_t\right)\d t\\
&+\int_{\R^d} \left(  2 (X_t-z)\sin(\<v_0,X_t-z\>)+(1+|X_t-z|^2) \cos(\<v_0,X_t-z\>)v_0\right)\sL_{X_t}(\d z)\d t.
\end{align*}
\end{exa}
\beg{proof}
It is easy to see that \eqref{V-poly} and \eqref{exa-V0}  hold for $\ga_0>4$ and $\ga_1=4$, and \eqref{exa-nnH} holds for $\ga_4=2$.  For  \eqref{exa-nnV0}, we have 
\beg{align*}
\<V_0(x),x\>&=4|A_0x|^2\<A_0^2x,x\>-2|x|^2\sin (|x|)-|x|^2\cos(|x|)|x|\\
&\geq 4a_0^4|x|^4-2|x|^2-|x|^3,
\end{align*}
which yields that \eqref{exa-nnV0} holds for $\ga_1=4$. For \eqref{HH}, we have 
\beg{align*}
&|H(x-y_1)-H(x-y_2)|\\
&\leq \left(2(|x|+|y_1|+|y_2|)+|v_0|\left(1+(|x|+|y_1|+|y_2|)^2\right) \right)|y_1-y_2|\\
&\leq C\left(1+ |x|+|y_1|+|y_2| \right)^{2}(|y_1|+|y_2|)\\
&\leq C\left(1+ |x|\right)^{2}\left(1+|y_1|+|y_2|\right)^{3}.
\end{align*}
Thus, \eqref{HH} holds for $\ga_2=2$ and $\ga_3=3$.  Moreover, \eqref{ga1234} holds.

\end{proof}
The condition $\ga_1>\ga_4+1$ in \eqref{ga1234} may exclude some models such as Example \ref{exa11}. We discuss the following example, which covers Example \ref{exa11} and is a complement to Corollary \ref{Exa-Gr}.
\beg{exa}
Consider the following equation:
\beg{equation}\label{Hamil-New0}
\mu(\d x,\d y)=\ff {\exp\{-\ff 2{\si^2}\left(\ff {y^2} 2-(1-\de) xy+\ff {x^4} 4-\ff {\de} 2x^2\right)-\ff {\be} {\si^2}\int_{\R}(x-z)^2\mu_1(\d z)\}\d x\d y} {\int_{\R}\exp\{-\ff 2{\si^2}\left(\ff {y^2} 2-(1-\de) xy+\ff {x^4} 4-\ff {\de} 2x^2\right)-\ff {\be} {\si^2}\int_{\R}(x-z)^2\mu_1(\d z)\}\d x\d y},
\end{equation}
where $\mu_1(\d x)=\int_{\R}\mu(\d x,\d y)$ is a marginal of $\mu$, $\be,\si$ are positive constants and $\de\in\R$. Solutions of this equation can be associated with stationary distributions of the following SDE on $\R^2$
\begin{equation}\label{Hamil-New1}
\begin{cases}
\d X_t = (Y_t-(1-\de) X_t)\d t\\
\d Y_t  = \left(-X_t^3+(1+\de(1-\de))X_t-\de Y_t\right)\d t-\be\int_{\R}(X_t-z)\sL_{X_t}(\d z)\d t+\si \d B_t.
\end{cases}
\end{equation}
Let 
\[\bar \mu(\d x)=\ff {e^{-\ff 2 {\si^2}\left(\ff {y^2} 4+\ff {x^4} 8\right)}} {\int_{\R^2} e^{-\ff 2 {\si^2}\left(\ff {y^2} 4+\ff {x^4} 8\right)}\d x\d y}\d x\d y,~W(x,y)=x^4+y^2+1.\]
Then for any $q\in [1,+\infty)$, \eqref{Hamil-New0} has a solutions $\mu$ with  $\ff {\d \mu} {\d \bar\mu}\in L^\infty\cap L^1(W\bar\mu)\cap   \cW^{1,\infty}_{q,\bar\mu}$.
\end{exa}
\beg{proof}
We set  $\bar V(x,y)=\ff 2 {\si^2}\left(\ff {y^2} 4+\ff {x^4} 8\right)$, and 
\[V_0(x,y)=\ff 2{\si^2}\left(\ff {y^2} 2-(1-\de) xy+\ff {x^4} 4-\ff {\de} 2x^2\right),~V(x,\mu)=\ff {\be} {\si^2} \int_{\R}(x-z)^2\mu_1(\d z).\]
Then it can be checked directly that {\bf (H)} holds for 
\[W_0(x,y)=(1+|x|)^2,~F_0(x,y)=F_2(x,y)=\ff {2\be} {\si^2}(1+|x|)^2,~F_1(r)=r,~F_3(r)=r+1.\]
Let $W_1(x,y)=x^2+y^2$ and 
\beg{align}\label{Lmu2.44}
L_{\mu}g(x,y)&=\De g(x,y)-\ff 2 {\si^2}(y-(1-\de)x)\pp_yg(x,y)\nonumber\\
&\quad\,+\left(\ff 2 {\si^2}((1-\de)y-x^3+\de x)-\ff {2\be} {\si^2}\int_{\R}(x-z)\mu_1(\d z)\right)\pp_xg(x,y),~g\in C^2.
\end{align}
Then the H\"older inequality implies that there are positive constants $C(\de,\be,\si^{-1}),C(\be,\si^{-1})$ depending on $\de,\be,\si^{-1}$ such that
\beg{align*}
L_{\mu}W_1(x,y)&=4-\ff 4 {\si^2}(y-(1-\de)x)y\\
&\quad\,+\left(\ff 4 {\si^2}((1-\de)y-x^3+\de x)-\ff {2\be} {\si^2}\int_{\R}(x-z)\mu_1(\d z)\right)x\\
&\leq -\ff 2 {\si^2}\left(x^4+y^2\right)+C(\de,\be,\si^{-1})+C(\be,\si^{-1})\left(\int_{\R}|z|\mu_1(\d z)\right)^{\ff 4 3}\\
&\leq -\ff 2 {\si^2}\left(x^4+y^2+1\right)+C(\de,\be,\si^{-1})+\ff 2 {\si^2}\\
&\quad\,+C(\be,\si^{-1})\left(\int_{\R}(|z|^4+|y|^2+1)\mu(\d z,\d y)\right)^{\ff 1 3}.
\end{align*}
Thus {\bf (W)} holds for $G_1(r)=C(\be,\si^{-1})r^{\ff 1 3}$, $G_2(r)=\ff 2 {\si^2} r-C(\de,\be,\si^{-1})+\ff 2 {\si^2}$.

Therefore, this example is proved by Theorem \ref{exis-thm1}.

\end{proof}

\beg{rem}\label{Rem2.1}
Setting $\de=1$, we obtain \eqref{Ham-exp} from \eqref{Hamil-New0}. Since \eqref{Hamil-New1} is a degenerate SDE, we can not use \cite{ZSQ} directly to prove that the existence of solutions to \eqref{Hamil-New0}. Solutions of \eqref{Hamil-New0} can be stationary distributions not only of \eqref{Hamil-New1}, but also  of the following SDE 
\begin{equation}\label{Hami-New2}
\begin{cases}
\d X_t=\sq 2\d B_t^{(1)}+\left(\ff 2 {\si^2}((1-\de)Y_t-X_t^3+\de X_t)-\ff {2\be} {\si^2}\int_{\R}(X_t-z)\sL_{X_t}(\d z)\right)\d t\\
\d Y_t=\sq 2\d B_t^{(2)}-\ff 2 {\si^2}(Y_t-(1-\de)X_t)\d t
\end{cases}
\end{equation}
where $(B_t^{(1)},B_t^{(2)})$ is a two dimensional Brownian motion.  \eqref{Hami-New2} can be derived  from \eqref{SDE-1.4} and \eqref{bbb} by setting $d=2$, $\si\equiv \sq{\al}I$ (identity matrix on $\R^2$), $\al=\ff {\be} {\si^2}$, $\th(\al)=\ff 2 {\be}\al$, and 
\[V_0(x,y)=\left(\ff {y^2} 2-(1-\de) xy+\ff {x^4} 4-\ff {\de} 2x^2\right),~V(x,\mu)= \int_{\R}(x-z)^2\mu_1(\d z).\]
Equation \ref{Hami-New2} is the McKean-Vlasov diffusion associated with $L_{\mu}$ defined by \eqref{Lmu2.44}.  We can apply \cite[Theorem 2.2]{ZSQ} to \eqref{Hami-New2} and obtain the existence of solutions to \eqref{Hamil-New0}. From this, one can see the difference between this paper and \cite{ZSQ}. In \cite{ZSQ}, we consider the existence of stationary distributions for a given SDE, and in this paper, we focus on the existence of solutions to an equation on stationary distributions, which may be associated with different SDEs. 
\end{rem}

The second corollary shows that our criteria can be applied to McKean-Vlasov diffusions with singular drifts. 
\beg{cor}\label{exa-singular0}
Consider the following equation
\beg{equation}\label{eq-exa1}
\mu(\d x)=\ff {\exp\{-V_0(x)+\sum_{i,j=1}^{m}h_i(x)H_{ij}\int_{\R^d}\th_j(y)\mu(\d y)\}\d x} {\int_{\R^d}\exp\{-V_0(x)+\sum_{i,j=1}^{m}h_i(x)H_{ij}\int_{\R^d}\th_j(y)\mu(\d y)\}\d x},
\end{equation}
where $V_0\in C^1(\R^d)$ satisfies \eqref{V-poly}-\eqref{exa-nnV0}$, m\in \N$, $H_{ij}\in\R$, $\{h_i\}_{i=1}^m$ and $\{\th_j\}_{j=1}^m$ are measurable functions. Suppose  there are nonnegative constants $C$, $\ga_2,\ga_3,\ga_4$  so that $\ga_5\in [0,1)$, $\ga_4\leq \ga_3<\ga_1$ and $\ga_1>\ga_2+\ga_4+1$,   
\begin{align*}
|\th_i(x)|&\leq C(1+|x|)^{\ga_2},\qquad\quad  |h_i(x)|\leq C(1+|x|)^{\ga_3},\\
|\nn h_i(x)|&\leq C(1+|x|^{\ga_4}+|x|^{-\ga_5}),~x\in\R^d-\{0\},~i=1,\cdots,m.
\end{align*}
Let 
\[\bar\mu(\d x)=\ff {e^{-\ff {C_0} {2}(1+|x|)^{\ga_1}}\d x} {\int_{\R^d} e^{-\ff {C_0} {2}(1+|x|)^{\ga_1}}\d x},~~W(x)=|x|^{\ga_1}+1.\] 
Then for any $p\in (d,\ff {d} {\ga_5})$ and $q\in [1,\ff {d} {\ga_5})$, \eqref{eq-exa1} has a solutions $\mu$ with  $\ff {\d \mu} {\d \bar\mu}\in L^\infty\cap L^1(W\bar\mu)\cap   \cW^{1,p}_{q,\bar\mu}$.
\end{cor}
\beg{proof}
We first check {\bf (H)}. Set 
\[V(x,\mu)=\sum_{i,j=1}^m h_i(x)H_{ij}\mu(\th_j),~~~\bar V(x)=\ff {C_0} {2}(1+|x|)^{\ga_1},~~~W_0(x)=(1+|x|)^{\ga_2}.\] 
Then $V_0,\bar V\in \cW^{1,\infty}_{q,\bar\mu}$ for any $q\geq 1$, $W_0\in L^1(\bar\mu)$. Thus $(H1)$ holds.  For all $\mu,\nu\in \sP_{W_0}$
\beg{align*}
|V(x,\mu)-V(x,\nu)|&\leq \sum_{i,j=1}^m |h_i(x)|\cdot|H_{ij}|\cdot|\mu(\th_j)-\nu(\th_j)|\\
&\leq \left(C^2\sum_{i,j=1}^m|H_{ij}|\right)(1+|x|)^{\ga_3}\|\mu-\nu\|_{W_0},\\
|V(x,\bar\mu)|&\leq \left(C^2\sum_{i,j=1}^m|H_{ij}|\right)(1+|x|)^{\ga_3}\bar\mu((1+|\cdot|)^{\ga_2}),\\
|\nn V(\cdot,\mu)(x)|&\leq \left(C^2\sum_{i,j=1}^m|H_{ij}|\right)(1+|x|^{\ga_4}+|x|^{-\ga_5})\|\mu\|_{W_0}.
\end{align*}
Set $F_0(x)=(1+|x|)^{\ga_3}$, $F_1(r)=F_3(r)=\left(C^2\sum_{i,j=1}^m|H_{ij}|\right)r$, $F_2(x)=(1+|x|^{\ga_4}+|x|^{-\ga_5})$. Then $F_2\in L^q(\bar\mu)\bigcap L^p_{loc}$ for any $p\in (d,\ff {d} {\ga_5})$ and $q\in [1,\ff {d} {\ga_5})$, and $(H2)$ holds. Due to the H\"older inequality and $\ga_3<\ga_1$, there is $\tld C>0$ such that
\beg{align*}
-V_0(x)+\be F_0(x)&\leq -C_0(1+|x|)^{\ga_1}+C_2+\be(1+|x|)^{\ga_3}\\
&\leq -\ff {C_0} 2(1+|x|)^{\ga_1}+\tld C\be^{\ff {\ga_1} {\ga_1-\ga_3}},~\be\geq 0.
\end{align*}
It is clear that $F_0\in L^1(W_0\bar\mu)$. Hence, {\bf (H)} holds. 

Set $W_1(x)=|x|^2$, $W(x)=(|x|+1)^{\ga_1}$ and 
\beg{align*}
L_\mu&=\De -\nn V_0\cdot\nn+\sum_{i,j=1}^m\nn h_iH_{ij}\mu(\th_j)\cdot\nn~,~\mu\in \sP_{W_0}.
\end{align*}
Then $W_1\in \cW^{2,p_1}_{\bar\mu}$ for any $p_1\in[ 1,+\infty]$. By using the H\"older inequality, $0<1-\ga_5<\ga_4+1$, and $\ga_1>\ga_4+\ga_2+1$,  there exist positive constants $\tld C_1,\tld C_2,\tld C_3$ such that 
\beg{align*}
L_{\mu}|x|^2&=2d-2\<\nn V_0(x),x\>+2\sum_{i,j=1}^m\<\nn h_i(x),x\>H_{ij}\mu(\th_j)\\
&\leq -2C_2(1+|x|)^{\ga_1}+2(C_3+d)\\
&\quad\, +\left(2C^2\sum_{i,j=1}^m|H_{ij}|\right) (1+|x|^{\ga_4+1}+|x|^{1-\ga_5})\mu((1+|\cdot|)^{\ga_2})\\
&\leq -\tld C_1(|x|+1)^{\ga_1}+\tld C_2+\tld C_3(\mu((1+|\cdot|)^{\ga_2})^{\ff {\ga_1} {\ga_1-\ga_4-1}}\\
&\leq -\tld C_1(|x|+1)^{\ga_1}+\tld C_2+\tld C_3(\mu((1+|\cdot|)^{\ga_1})^{\ff {\ga_2} {\ga_1-\ga_4-1}}\\
&=-\tld C_1 W(x)+\tld C_2+\tld C_3\|\mu\|_{W}^{\ff {\ga_2} {\ga_1-\ga_4-1}}.
\end{align*}
This, together with $\ga_2<\ga_1-\ga_4-1$, yields that {\bf (W)} holds with $G_1(r)=\tld C_3r^{\ff {\ga_2} {\ga_1-\ga_4-1}}$ and $G_2(r)=\tld C_1r-\tld C_2$. 

Therefore, for any $p\in (d,\ff {d} {\ga_5})$ and $q\in [1,\ff {d} {\ga_5})$, Theorem \ref{exis-thm1} implies that \eqref{eq-exa1} has a fixed point $\mu$ and $\ff {\d \mu} {\d \bar\mu}\in L^\infty\cap L^1(W\bar\mu)\cap   \cW^{1,p}_{q,\bar\mu}$.

\end{proof}

We give a concrete example to finish this subsection, which is a direct consequence of Corollary \ref{exa-singular0}, and we omit the proof. 
\beg{exa}\label{exa-singular}
Consider the following equation
\beg{equation}\label{eq-exa2}
\mu(\d x)=\ff {\exp\{-C_1\ff {|x|^4} 4+C_2 \ff {|x|^2} 2+|x|^{\ga_3-1}\int_{\R^d}\<x,\th(y)\>\mu(\d y)\}\d x} {\int_{\R^d}\exp\{-C_1\ff {|x|^4} 4+C_2 \ff {|x|^2} 2+|x|^{\ga_3-1}\int_{\R^d}\<x,\th(y)\>\mu(\d y)\}\d x},
\end{equation}
where $C_1,C_2$ are positive constants, $\ga_3\in (0,4)$,  $\th$ is a $\R^d$-valued measurable function and there exist  $C_3\geq 0$, $\tld \ga_2\in (0,4)$ such that  
\beg{equation}\label{ad-th0}
|\th(y)|\leq C_3\left(1+|y|^{\ff {\tld\ga_2(4-\ga_3)} 4}\right),~y\in\R^d.
\end{equation}
Let 
\[\bar\mu(\d x)=\ff {e^{-\ff {C_1} {8}|x|^4}\d x} {\int_{\R^d} e^{-\ff {C_1} {8}|x|^4}\d x},~~W(x)=|x|^4+1.\] 
Then for any $p\in (d,\ff {d} {(1-\ga_3)^+})$ and $q\in [1,\ff {d} {(1-\ga_3)^+})$, \eqref{eq-exa2} has a solutions $\mu$ with  $\ff {\d \mu} {\d \bar\mu}\in L^\infty\cap L^1(W\bar\mu)\cap   \cW^{1,p}_{q,\bar\mu}$.
\end{exa}
\beg{rem}
We can obtain Dawson's model from this example by setting  $d=1$,  $C_1=\ff 2 {\si^2}$, $C_2=\ff {2(1-\be)} {\si^2}$, $\ga_3=1$, $\th(y)=-\ff {2\be} {\si^2} y$, and $\th$ satisfies \eqref{ad-th0} with $\tld\ga_2=\ff 3 2$. 

According to \eqref{V2FF}, $V(x,\cdot)$ is continuous with respect to the metric $\|\cdot\|_{W_0}$. We emphasize that  $V(x,\cdot)$ may be discontinuous in the weak topology and the total variation norm. For instance, in Example \ref{exa-singular},  let
\[\th(y)=\sum_{i=1}^d\left(1+|y_i|^{\ff {\tld\ga_2(4-\ga_3)} 4}\right){\rm sgn}\left(\sin\left(\ff 1 {y_i}\right)\right)e_i,~y=(y_1,y_2,\cdots,y_d)\in\R^d,\]
where $\ga_1,\ga_2$ are positive constants in Example \ref{exa-singular}, ${\rm sgn}(\cdot)$ is  the sign function and $\{e_i\}_{i=1}^d$ is the standard orthonormal basis of $\R^d$. Then $\th$ is a discontinuous function satisfying \eqref{ad-th0}, and the mapping $\mu\mapsto\int_{\R^d}\th(y)\mu(\d y)$ is not continuous in the weak topology and the total variation norm. It is required in \cite[(H3)]{ZSQ} that  coefficients are continuous in distribution variable w.r.t. the weak topology. Thus, Example \ref{exa-singular} with $\th$ defined as above can not be covered by \cite{ZSQ}. 
Moreover, solutions of \eqref{eq-exa2} can be associated with  stationary solutions to the  following McKean-Vlasov equation:
\begin{align*}
\d X_t&=\d B_t-(C_1|X_t|^2X_t-C_2X_t)\d t\\
&\quad\, +\left(|X_t|^{\ga_3-1}\int_{\R^d}\th(y)\sL_{X_t}(\d y)+(\ga_3-1)|X_t|^{\ga_3-3}\int_{\R^d} \<X_t,\th(y)\> \sL_{X_t}(\d y)X_t\right)\d t.
\end{align*}
When $\ga_3<1$, the drift term is singular in spatial variable and discontinuous in distribution variable. 
\end{rem} 

\subsection{Proof of Theorem \ref{exis-thm1}}

We use the Schauder fixed point theorem to prove Theorem \ref{exis-thm1}.  So, we first investigate the continuity of $\hat\cT$ (see Lemma \ref{ps(mu)}), and find a nonempty closed   convex  subset $\cM_M$ in $L^1(W_0\bar\mu)$ (see \eqref{cMM}), which is also  an invariant subset  of $\hat\cT\circ \mathscr{I}$ (see Lemma \ref{T-inva}). Then we prove that $\hat\cT\circ \mathscr{I}$ is compact on $\cM_M$ (see Lemma \ref{T-comp}), and the Schauder fixed point theorem can be applied to $\hat\cT\circ \mathscr{I}$ on $\cM_M$.

\beg{lem}\label{ps(mu)}
Assume {\bf (H)}. Then  for each $\mu\in\sP_{W_0}$, there is $\ps(\mu)\in \cW^{1,p}_{q,\bar\mu}\cap L^\infty$. Furthermore, if $F_0\in L^1(W_0\bar\mu)$, then $\hat \cT$ is continuous from $\sP_{W_0}$ to $L^1(W_0\bar\mu)$. Consequently, $\hat\cT\circ \mathscr{I}$ is continuous on $\{f\in L^1(W_0\bar\mu)~|~f\bar\mu\in \sP_{W_0}\}$, which inherits the metric induced by the norm $\|\cdot\|_{L^1(W_0\bar\mu)}$, and the following mapping is also continuous on $\sP_{W_0}$:
\beg{align}\label{map-Tmu}
\tld\cT:~\mu\in\sP_{W_0}\mapsto \tld\cT_\mu:=\hat\cT(\mu)\bar\mu.
\end{align}
\end{lem}

\beg{proof}
For each $\mu\in \sP_{W_0}$, it follows from \eqref{V2FF}, \eqref{bmu-F} and \eqref{V-F3} that 
\beg{equation}\label{VVF3}
\beg{split}
-V_0(x)+|V(x,\mu)| & =-V_0(x)+|V(x,\mu)-V(x,\bar\mu)|+  |V(x,\bar\mu)|\\
&\leq -V_0(x)+ F_0(x)F_1(\|\mu-\bar\mu\|_{W_0})+|V(x,\bar\mu)|\\
&\leq -\bar V(x)+F_4( F_1(\|\mu-\bar\mu\|_{W_0})+ C)+  C.
\end{split}
\end{equation}
This implies that $\ps(\mu)\in L^\infty$. Lemma \ref{ef} implies that $\ps(\mu)\in W^{1,p}_{loc}$ with
$$\nn \ps(\mu)=\ps(\mu)(-\nn V_0-\nn V(\mu)+\nn\bar V).$$
Taking into account that $\ps(\mu)\in L^\infty$, $\nn V_0,\nn \bar V\in L^q(\bar\mu)$ and  \eqref{nnV2} which yields $V(\mu)\in \cW^{1,p}_{q,\bar\mu}$, we find that $\ps(\mu)\in \cW^{1,p}_{q,\bar\mu}\cap L^\infty$.

Next, we prove the continuity of $\hat\cT$.  It follows from the H\"older inequality and \eqref{VVF3} that for each $\mu\in\sP_{W_0}$
\beg{equation}\label{b-mups}
\beg{split}
\left(\bar\mu(\ps(\mu))\right)^{-1}&\leq \left(\bar\mu(e^{-V_0+\bar V})\right)^{-2}\bar\mu\left(\exp\left(-V_0+\bar V+V(\mu)\right)\right)\\
&\leq \left(\int_{\R^d} e^{-V_0(x)}\d x\right)^{-2}\exp\left[F_4( F_1(\| \mu-\bar\mu\|_{W_0})+C)+  C\right].
\end{split}
\end{equation}
Thus $\left(\bar\mu(\ps(\cdot))\right)^{-1}$ is locally bounded on $ \sP_{W_0}$.
For  $\mu_1,\mu_2\in\sP_{W_0}$,
\beg{align}\label{T-T}
\hat\cT(\mu_2)-\hat\cT(\mu_1) = \ff {\ps(\mu_2)-\ps(\mu_1)} {\bar\mu(\ps(\mu_2))}+\ff {\ps(\mu_1)(\bar\mu(\ps(\mu_1)-\ps(\mu_2)))} {\bar\mu(\ps(\mu_1))\bar\mu(\ps(\mu_2))}.
\end{align}
Fix $\mu_1$. Then we derive from \eqref{b-mups} and \eqref{T-T} that,  to prove the continuity of $\hat\cT$, it is sufficient to prove $\ps(\mu)$ is continuous in $\mu$.  It follows from \eqref{funi-1}, \eqref{V2FF} and \eqref{VVF3} that
\beg{align*}
 \left|\ps(\mu_2)-\ps(\mu_1) \right| &\leq  \left|V(\mu_2)-V(\mu_1)\right| e^{-V_0-V(\mu_2)\we V(\mu_1)+\bar V}\\
&\leq F_0 F_1(\|\mu_2-\mu_1\|_{W_0})\exp[F_4( F_1(\|\mu_1-\bar\mu\|_{W_0})\vee\|\mu_2-\bar\mu\|_{W_0}) + C)+  C].
\end{align*}
Due to  $F_0\in L^1(W_0\bar\mu)$ and $\lim_{r\ra 0^+}F_1(r)=0$, we find that   
\beg{align*}
\lim_{\|\mu_2-\mu_1\|_{W_0}\ra 0}\left\|\ps(\mu_2)-\ps(\mu_1) \right\|_{L^1(W_0\bar\mu)}=0.
\end{align*}
This also implies that 
$$\lim_{\|\mu_2-\mu_1\|_{W_0}\ra 0}\bar\mu\left(\left|\ps(\mu_1)-\ps(\mu_2)\right|\right)=0,$$
since $W_0\geq 1$.  

Finally, for any nonnegative functions $f_1,f_2\in L^1(W_0\bar\mu)$ with $\bar\mu(f_1)=\bar\mu(f_2)=1$, we have that
\beg{equation}\label{PW0-L1}
\beg{split}
\left\|f_1 -f_2\right\|_{L^1(W_0\bar\mu)}&=\sup_{\|g\|_\infty\leq 1}\left|\bar\mu(W_0\left(f_1-f_2 \right)g)\right|=\sup_{|\tld g|\leq W_0}\left|\bar\mu( \left(f_1-f_2 \right)\tld g)\right|\\
&= \|f_1\bar\mu-f_2\bar\mu\|_{W_0}.
\end{split}
\end{equation}
We derive from  this equality that $\hat\cT\circ \mathscr{I}$ is continuous on $\{f\in L^1(W_0\bar\mu)~|~f\bar\mu\in \sP_{W_0}\}$, and the mapping $\tld\cT$ is continuous on $\sP_{W_0}$

\end{proof}

For $W_0,W$ satisfying \eqref{W0W} and each $M\in (\bar\mu(W),+\infty)$, we introduce the following set 
\beg{align}\label{cMM}
\cM_{M}=\{f\in L^1(W_0\bar\mu)~|~f\geq 0,~\bar\mu(f)=1,\bar\mu(Wf)\leq M\}.
\end{align}
Then $\cM_{M}$ is a nonempty closed and convex subset of $ L^1(W_0\d\bar\mu)$ and $\mathscr{I}(\cM_{M})\subset\sP_W$. Due to \eqref{W0W}, for each $f\in\cM_{M}$
$$\bar\mu(|f|W_0)\leq \left(\sup_{x\in\R^d}\ff {W_0(x)} {W(x)}\right)\bar\mu(fW)\leq M\left(\sup_{x\in\R^d}\ff {W_0(x)} {W(x)}\right).$$
Thus, $\cM_{M}$ is also bounded in $L^1(W_0\bar\mu)$. 

\beg{lem}\label{T-inva}
Assume that {\bf(H)} holds,  and  {\bf(W)} holds for $W_1\in \cW^{2,p_1}_{\bar\mu}$ and $p_1\geq \ff { q} { q-1}$. Then  there is $M_0>0$ such that $\hat\cT(f\bar\mu)\in \cM_{M}$ for every $f\in\cM_M$ and $M>M_0$.
\end{lem}
\beg{proof}
For $\mu\in\sP_W$, we have that $\mu\in\sP_{W_0}$ since \eqref{W0W}. Thus $\nn\log(\ps(\mu)e^{-\bar V})\in L^q(\bar\mu)$ due to {\bf(H)}. Then 
\beg{align}\label{aW1log-L1}
\nn W_1,~\nn^2 W_1,
\< \nn \log(\ps(\mu)e^{-\bar V}),\nn W_1\>\in L^1(\bar\mu).
\end{align}
According to Lemma \ref{ps(mu)}, $\ps( \mu)\in L^\infty$.  This, together with  \eqref{aW1log-L1}, yields that  $L_{\mu}W_1\in L^1(\bar\mu)\subset L^1(\tld\cT_\mu)$  and 
\[\int_{\R^d} |\nn W_1|(x)\ps(x,\mu)\bar\mu(\d x)\leq \|\ps(\mu)\|_\infty\bar\mu(|\nn W_1|)<\infty.\]
Then, as \eqref{add-in00}, 
\beg{align*}
\left|\tld\cT_\mu(L_{ \mu}W_1)\right|&=\lim_{n\ra +\infty}\left|\tld\cT_\mu(\ze_n(L_{\mu}W_1))\right|\\
&\leq \lim_{n\ra+\infty} \left(\ff 2 n\|\ps(\mu)\|_\infty \bar\mu(|\nn W_1|) \right)=0.
\end{align*}
This, together with \eqref{LYP} and the Jensen inequality, yields that
\beg{align*}
0&=\tld\cT_\mu(L_{ \mu}W_1)\leq G_1(\| \mu\|_W)-\tld\cT_\mu(G_2(W))\\
&\leq G_1(\| \mu\|_W)-G_2(\tld\cT_\mu(W))= G_1(\| \mu\|_{W})-G_2(\|\tld\cT_\mu\|_W).
\end{align*}
According to \eqref{G12}, there is $M_0>0$ such that
$$G_1(r)<G_2(r),\quad r>M_0.$$
Then for every $M>M_0$ and  $\mu\in\sP_W$ such that $\|\mu\|_W\leq M$, there is
\beg{align*}
G_2(\|\tld\cT_\mu\|_W)\leq G_1(\|\mu\|_W)\leq G_1(M)<G_2(M),
\end{align*}
where we have used in the second inequality that $G_1$ is increasing. This implies that $\|\tld\cT_\mu\|_W\leq M$, since $G_2$ is increasing.  For each $f\in\cM_{M}$, 
\beg{align}\label{iso}
\|f\bar\mu\|_{W}=\sup_{|g|\leq W_0}\bar\mu(gf)=\bar\mu(Wf)\leq M,
\end{align}
which implies that $f\bar\mu\in\sP_{W}$. Hence, $\hat\cT(f\mu)\in \cM_{M}$ for every $f\in\cM_M$ and $M>M_0$.

\end{proof}

It follows from   \eqref{PW0-L1} that $\mathscr{I}$ is an isometric mapping from $\cM_M$ onto $\tld\cM_M$:
$$\tld\cM_M\equiv\{f\bar\mu ~|~f\in\cM_{M}\}\subset \sP_{W_0}$$
which is equipped with the weighted  variation metric $\|\cdot\|_{W_0}$.  For $M>M_0$, since $\hat\cT(f\bar\mu)\in \cM_{M}$ for every $f\in\cM_M$, $\hat\cT\circ \mathscr{I}$ is a mapping from $\cM_M$ to itself.  If $F_0\in L^1(W_0\bar\mu)$, then $\hat\cT\circ \mathscr{I}$ is continuous on $\cM_M$, according to Lemma \ref{ps(mu)}. 

Next, we prove that $\hat\cT\circ \mathscr{I} $ is compact on $\cM_{ M}$.  

\beg{lem}\label{T-comp}
Suppose that the   assumption of Theorem \ref{exis-thm1} holds. Then, for every $M>M_0$, $\hat\cT\circ \mathscr{I}$ is compact on $\cM_{ M}$.
\end{lem}
\beg{proof}
Let $\{f_n\}_{n\geq 1}$ be a sequence in $\cM_{ M}$. We have prove that $\hat\cT\circ \mathscr{I}$ is continuous on $\cM_M$, due to Lemma \ref{ps(mu)}.  To prove that $\hat\cT\circ \mathscr{I}$ is compact on $\cM_{ M}$, it is sufficient to  prove that there is a subsequence $\{\hat\cT (f_{n_k}\bar\mu)\}_{k\geq 1}$ converging in $L^1(W_0\bar\mu)$. 

We first prove that for every $N>0$, $\{\hat\cT(f_n\bar\mu)\ze_N^2\}_{n\geq 1}$ is bounded in $W^{1,p}(B_{2N})$.  Since \eqref{b-mups}, \eqref{PW0-L1} and  $\cM_M$ is bounded in $L^1(W_0\bar\mu)$, it is sufficient to prove that $\ps(f_n\bar\mu)\ze_N^2$ is bounded in $W^{1,p}(B_{2N})$.  Due to $\log(\ps(f_n\bar\mu))\in W^{1,p}_{loc}$ and $\ps(f_n\bar\mu)\in W^{1,p}_{loc}$, we have that $\ps(f_n\bar\mu)\ze_N^2\in W^{1,p}(B_{2N})$. It follows from \eqref{nnV2} that
\beg{equation}\label{nnpN}
\beg{split}
&\left|\nn(\ps(f_n\bar\mu)\ze_N^2)\right|\\
&\leq \left|\nn \ps(f_n\bar\mu)\right|\ze_N^2+2\ps(f_n\bar\mu)\ze_N|\nn\ze_N|\\
&\leq \ps(f_n\bar\mu)\ze_N\left(\left(|\nn V_0|+|\nn V(f_n\bar\mu)|+|\nn \bar V|\right)\ze_N+ 2|\nn\ze_N|\right)\\
&\leq \ps(f_n\bar\mu)\ze_N\left(\left(|\nn V_0|+F_2 F_3(\|f_n\bar\mu\|_{W_0})+|\nn \bar V|\right)\ze_N+2|\nn\ze_N|\right)\\
&= \ps(f_n\bar\mu)\ze_N\left(\left(|\nn V_0|+F_2 F_3(\bar\mu(W_0f_n))+|\nn \bar V|\right)\ze_N+2|\nn\ze_N|\right)\\
&\leq \ps(f_n\bar\mu)\ze_N\left(\left(|\nn V_0|+F_2 F_3(\tld C\bar\mu(f_nW))+|\nn \bar V|\right)\ze_N+2|\nn\ze_N|\right)\\
&\leq \ps(f_n\bar\mu)\ze_N\left(\left(|\nn V_0|+F_2 F_3(\tld CM)+|\nn \bar V|\right)\ze_N+2|\nn\ze_N|\right),
\end{split}
\end{equation} 
where $\tld C=\sup_{x\in\R^d}\ff {W_0} {W}(x)$.  Since $\nn V_0,F_2,\nn\bar V\in L^p_{loc}$, there is a positive constant $C_{V_0,F_2,F_3,\bar V,N,\tld C,M}$ which is independent of $n$ such that
$$ \left\|\left(|\nn V_0|+F_2 F_3(\tld CM)+|\nn \bar V|\right)\ze_N\right\|_{L^p}\leq C_{V_0,F_2,F_3,\bar V,N,\tld C,M}.$$
Putting this into \eqref{nnpN},  there is a constant $\hat C$ which depends on $V_0,F_2,F_3,\bar V,N,\tld C,M$ and  is independent of $n$ such that   
\beg{align}\label{sup-nnps}
\sup_{n\geq 1}\left\|\nn(\ps(f_n\bar\mu)\ze_N^2)\right\|_{L^p}\leq \hat C\sup_{n\geq 1}\|\ps(f_n\bar\mu)\ze_N\|_{\infty}.
\end{align}
It follows from $p>d$ and the Morrey embedding theorem (\cite[Theorem 9.12]{Bre}) that 
$$\|V_0\ze_{2N}\|_\infty+\|\bar V\ze_{2N}\|_{\infty}\leq C\left(\|V_0\ze_{2N}\|_{W^{1,p}}+\|\bar V\ze_{2N}\|_{W^{1,p}}\right) <\infty.$$
By \eqref{V2FF} and \eqref{bmu-F}, we have that
\beg{align*}
\|V(f_n\bar\mu)\ze_{2N}\|_\infty &\leq \|(V(f_n\bar\mu)-V(\bar\mu))\ze_{2N}\|_\infty+  \|V(\bar\mu)\ze_{2N}\|_\infty\\
&\leq \|F_0\ze_{2N}\|_\infty F_2(\|f_n\bar\mu-\bar\mu\|_{W_0})+C(\|F_0\ze_{2N}\|_\infty+1)\\
&\leq \|F_0\ze_{2N}\|_\infty F_2(\tld C\bar\mu(Wf_n)+\bar\mu(W_0))+C(\|F_0\ze_{2N}\|_\infty+1)\\
&\leq \|F_0\ze_{2N}\|_\infty F_2(\tld CM+\bar\mu(W_0) )+C(\|F_0\ze_{2N}\|_\infty+1).
\end{align*}
Combining this with $F_0\in L^{\infty}_{loc}$,  we have that $\sup_{n}\|V(f_n\bar\mu)\ze_{2N}\|_\infty<\infty$. Consequently, 
\beg{align*}
\sup_{n\geq 1}\| \ps(f_n\bar\mu)\ze_N^2 \|_{L^p}&\leq C_N\sup_{n\geq 1}\| \ps(f_n\bar\mu)\ze_N \|_{\infty}\\
&\leq C_N\sup_{n\geq 1}\| e^{(|V_0|+|V(f_n\bar\mu)|+|\bar V|)\ze_{2N}}\ze_N \|_{\infty}<\infty.
\end{align*}
This and \eqref{sup-nnps} imply that
\beg{align*}
\sup_{n\geq 1}\|\ps(f_n\bar\mu)\ze_N^2\|_{W^{1,p}(B_{2N})}=\sup_{n\geq 1}\|\ps(f_n\bar\mu)\ze_N^2\|_{W^{1,p}}<\infty.
\end{align*}

Finally, we find a Cauchy subsequence from $\{\hat\cT(f_{n}\bar\mu )\}_{n\geq 1}$ by using Cantor’s diagonal argument. We have proven above  that, for all $N\in\N$,  $\{\hat\cT(f_{n}\bar\mu)\ze_N^2\}_{n\geq 1}$ is bounded in $W^{1,p}(B_{2N})$. For $N=1$, it follows from the Rellich-Kondrachov theorem (\cite[Theorem 9.16]{Bre}) that there is a subsequence $\{f_{n_{1,k}}\}_{k\geq 1}$ such that $\{\hat\cT(f_{n_{1,k}}\bar\mu)\ze_1^2\}_{k\geq 1}$ is a Cauchy sequence in $C(\overline{B_{2}})$. For $N\in\N$,  if a subsequence $\{f_{n_{N,k}}\}_{k\geq 1}$ of $\{f_n\}_{n\geq 1}$ has been  selected, then  by using the Rellich-Kondrachov theorem, we have a subsequence $\{f_{n_{{N+1},k}}\}_{k\geq 1}$ of $\{f_{n_{N,k}}\}_{k\geq 1}$ such that $\{\hat\cT(f_{n_{N+1,k}}\bar\mu)\ze_{N+1}^2\}_{k\geq 1}$ is a Cauchy sequence in $C(\overline{B_{2(N+1)}})$. By induction, we obtain subsequences $\{f_{n_{N,k}}\}_{N\geq 1, k\geq 1}$ such that for any $N\geq 1$, $\{f_{n_{N+1,k}}\}_{k\geq 1}$ is a subsequence of $\{f_{n_{N,k}}\}_{k\geq 1}$, and  $\{\hat\cT(f_{n_{N,k}}\bar\mu)\ze_N^2\}_{k\geq 1}$ is a Cauchy sequence in $C(\overline{B_{2N}})$. We choose a subsequence $\{f_{n_{k,k}}\}_{k\geq 1}$, which will be denoted by $\{f_{n_k}\}_{k\geq 1}$ for simplicity. Then for any $N\in\N$, $\{\hat\cT(f_{n_{k}}\bar\mu)\ze_{N}^2\}_{k\geq 1}$ is a Cauchy sequence in $C(\overline{B_{2N}})$. For each $N\in\N$, since $\{\hat\cT(f_{n_{k}} \bar\mu)\}_{k\geq 1}\subset \cM_M$, we have for any $k,k'\in\N$  that
\beg{align*}
&\bar\mu\left(W_0\left|\hat\cT(f_{n_k}\bar\mu )-\hat\cT(f_{n_{k'}}\bar\mu )\right|(1-\ze_N^2)\right)\\
&\leq \bar\mu\left(W_0\left|\hat\cT(f_{n_k}\bar\mu )-\hat\cT(f_{n_{k'}}\bar\mu )\right|\1_{[|x|\geq N]}\right)\\
&\leq  \left(\sup_{|x|\geq N}\ff {W_0} {W}(x)\right)\bar\mu\left(W\left|\hat\cT(f_{n_k}\bar\mu)-\hat\cT(f_{n_{k'}}\bar\mu )\right|\right)\\
&\leq 2M\left(\sup_{|x|\geq N}\ff {W_0} {W}(x)\right).
\end{align*}
Then  
\beg{equation}\label{L1-Cau}
\beg{split}
\bar\mu\left(W_0\left|\hat\cT(f_{n_k}\bar\mu )-\hat\cT(f_{n_{k'}}\bar\mu )\right|\right)
&\leq \bar\mu\left(W_0\left|\hat\cT(f_{n_k}\bar\mu )-\hat\cT(f_{n_{k'}}\bar\mu )\right|\ze_N^2\right)\\
&\quad +\bar\mu\left(W_0\left|\hat\cT(f_{n_k}\bar\mu )-\hat\cT(f_{n_{k'}}\bar\mu )\right|(1-\ze_N^2)\right)\\
&\leq \bar\mu(W_0)\left\|\left(\hat\cT(f_{n_k}\bar\mu )-\hat\cT(f_{n_{k'}}\bar\mu )\right)\ze_N^2\right\|_{B_{2N},\infty}\\
&\quad +2M\left(\sup_{|x|\geq N}\ff {W_0} {W}(x)\right).
\end{split}
\end{equation}
Hence, letting $k,k'\ra+\infty$ first and then $N\ra+\infty$, we derive from \eqref{L1-Cau} and \eqref{W0W} that $\{\hat\cT(f_{n_k}\bar\mu )\}_{k\geq 1}$ is a Cauchy sequence in $L^1(W_0\bar\mu)$.

\end{proof}

\noindent {\bf{Proof of Theorem \ref{exis-thm1}.~~}}
According to Lemma \ref{T-inva} and Lemma \ref{T-comp}, in $L^1(W_0\bar\mu)$, $\cM_{M}$ is a nonempty closed bounded and convex subset, and $\hat\cT\circ \mathscr{I}$ is compact from $\cM_{M}$ to $\cM_{M}$ for $M>M_0$.  Therefore, the Schauder fixed point theorem yields that $\hat\cT\circ \mathscr{I}$ has a fixed point in $\cM_M$ for $M>M_0$.  For any fixed point $f\in\cM_M$, we have that $f\in L^1(W\bar\mu)$, and according to Lemma \ref{ps(mu)}, $f\in W^{1,p}_{q,\bar\mu}\cap L^\infty$.

\qed

\section{Bifurcation}

Let $0<\hat\al<\check{\al}<+\infty$.  For $\al\in (\hat{\al},\check{\al})$ and $\th\in C^1((\hat{\al},\check{\al});(0,+\infty))$, we investigate the parameter $\al_0\in (\hat{\al},\check{\al})$ which more than one branch of solutions to \eqref{eq-bif} emanate from. For better reference, let us introduce the following definition of the bifurcation point for a parameter-dependent problem on $\sP(\R^d)\times J$, where $J\subset \R$ is an open interval.
\beg{defn}\label{bif-proba}
Let $\al_0\in J$, $\mathcal{P}\subset \sP(\R^d)$, $\mathcal{F}$ be a mapping from  $\mathcal{P}\times J$ to signed measures and $\{(\nu_\al,\al)~|~\al\in J,\nu_\al\in\mathcal{P}\}$ such that $\nu_\al$ is weakly continuous in $\al$ and 
\[\mathcal{F}(\nu_\al,\al)=0,~\al\in J.\]  
Then $(\nu_{\al_0},\al_0)$ is said to be a bifurcation point for the equation $\mathcal{F}(\mu,\al)=0$ if and only if there is a sequence $\{\mu_{\al_n}\}_{n\geq 1}\subset\mathcal{P}$  with $\al_n\in J$ such that $\mathcal{F}(\mu_{\al_n},\al_n)=0$, $\mu_{\al_n}\neq \nu_{\al_n}$ and 
\[\lim_{n\ra+\infty}\al_n=\al_0,~\lim_{n\ra+\infty}\mu_{\al_n}=\nu_{\al_0},\]
where $\mu_{\al_n}$  converges w.r.t. the weak topology.   For simplicity, we call that $(\mu_{\al_0},\al_0)$ is a bifurcation point.
\end{defn}
A bifurcation point for $\mathcal{F}(\mu,\al)=0$ is a cluster point of solutions for $\mathcal{F}(\mu,\al)=0$  except the known solution path $\{(\nu_\al,\al)~|~\al\in J,\nu_\al\in\mathcal{P}\}$. Usually, the bifurcation is investigated for a  parameter-dependent problem on the Banach space and the known solution path for a parameter-dependent problem is set to be the trivial solution,  see \cite{Deim,Kie}.  We first reformulate \eqref{eq-bif}.  Let $\bar V$ be a measurable function with $e^{-\bar V}\in L^1$, and let $\bar\mu$ be defined by \eqref{barmu}. Then we reformulate $\hat \cT\circ \mathscr{I}$ into the following form
\[ \hat\cT\circ\mathscr{I}(\rh,\al)=\ff {\exp\{-\th(\al)V_0-\al \int_{\R^d}V(\cdot,y)\rh(y)\bar\mu(\d y)+\bar V\}} {\int_{\R^d}\exp\{-\th(\al)V_0(x)-\al \int_{\R^d}V(x,y)\rh(y)\bar\mu(\d y)+\bar V(x)\}\bar\mu(\d x)}.\]
In this section, we denote $\cT =\hat\cT\circ\mathscr{I}$ for simplicity. Fix $\al\in (\hat{\al},\check{\al})$. The existence of the fixed point for $\cT(\cdot,\al)$ can be investigated by using results in Section 2.  If there is a family of fixed points for $\cT$, saying $\{\rh_\al\}_{\al\in (\hat{\al},\check{\al})}$, then we can set 
\beg{align}\label{PHPH}
\Ph(\rh,\al)=\rh_\al^{-1}\left((\rh+1)\rh_\al-\cT((\rh+1)\rh_\al,\al)\right), 
\end{align}
and $0$ is a trivial solution of 
\[\Ph(\cdot,\al)=0,~\al\in (\hat{\al},\check{\al}).\]
Moreover, for $\rh \in L^1(\mu_{\al})$ satisfying $\Ph(\rh,\al)=0$,  $(\rh+1)\rh_\al$ is a fixed point of $\cT(\cdot,\al)$ and the probability measure $(\rh+1)\rh_\al\bar\mu$ is a solution of \eqref{eq-bif}. Thus, we give a local bifurcation theorem for $\Ph(\rh,\al)=0$.

Before our detailed discussion, we explain our framework and strategy. We decompose  $V(x,y)$ into four part:
\begin{equation}\label{VVV}
V(x,y)=V_1(x)+V_2(x,y)+K_1(y)+K_2(x,y).
\end{equation}
$K_1$ can be canceled in $\cT$, see Remark \ref{canc-K1}. We assume there exist $\al_0\in (\hat{\al},\check{\al})$ and $\rh_{\al_0}\in L^1(\bar\mu)$  such that $\rh_{\al_0}$ is a fixed point of $\cT(\cdot,\al_0)$ without the terms $K_1$ and $K_2$.  We first prove in Lemma \ref{con-rhal} that $\rh_{\al_0}$ can be extended uniquely in  $L^2(\bar\mu)$ to a smooth path $\{\rh_{\al}\}_{\al\in [\al_0-\de,\al_0+\de]}$ such that $\rh_\al$ is also a fixed point of $\cT(\cdot,\al)$ without $K_1$ and $K_2$.  $K_2$ is assumed to be orthogonal to the path $\{\rh_{\al}\}_{\al\in [\al_0-\de,\al_0+\de]}$, see (A3) for the precise meaning. The condition (A3) ensures that $ \rh_{\al}$ remains a  fixed point of $\cT(\cdot,\al)$ for $\al\in [\al_0-\de,\al_0+\de]$, see Lemma \ref{lem-ps-T}. Then $\Ph$ is well-defined. In Corollary  \ref{cor:mu0mu}, we prove that $\{\rh_{\al}\}_{\al\in [\al_0-\de,\al_0+\de]}$ can be compared with $\rh_{\al_0}$, then the bifurcation analysis can be given in $L^2(\mu_{\al_0})$ for $\Ph=0$.  We prove that $\Ph(\cdot,\al)$ is Fr\'echet differentiable on $L^2(\mu_{\al_0})$ and the derivative $\nn \Ph(0,\cdot)$ is continuous  in $L^2(\mu_{\al_0})$ , see Corollary  \ref{cor:mu0mu} and Lemma \ref{con-nnPh}. Due to the Krasnosel’skii Bifurcation Theorem (\cite[Theorem II.3.2]{Kie}), if $\nn\Ph(0,\al)$ has an odd crossing number at $\al_0$ (Definition \ref{def-odd-cro}), then $\al_0$ is a bifurcation point of $ \Ph=0$ in the following sense:
\beg{defn}\label{bif-def-00}
$\al_0$ is a bifurcation point of $\Ph=0$ if and only if $(0,\al_0)$ a cluster point of nontrivial solutions  $(\rh,\al)\in L^2(\mu_{\al_0})\times [\al_0-\de,\al_0+\de]$, $\rh\neq 0$, of $\Ph(\rh,\al)=0$, i.e. 
\beg{equation}\label{bif-L2}
(0,\al_0)\in \overline{\{(\rh,\al)|\Ph(\rh,\al)=0,~\rh\neq 0,\al\in [\al_0-\de,\al_0+\de]\}},
\end{equation}
where the closure in the right-hand side is taken in $L^2(\mu_{\al_0})\times\R$. 
\end{defn}
\noindent Then we prove that $(\mu_{\al_0},\al_0)$ is a bifurcation point of \eqref{eq-bif} in the sense of Definition \ref{bif-proba}  in  Lemma  \ref{lem-bif-bif}. We use the regularized determinant for the Hilbert-Schmidt operator (see \cite{GLZ08,Simon}) to derive a criteria for $\nn\Ph(0,\al)$ has an odd crossing number at $\al_0$. Then our criteria for a bifurcation point of $\Ph=0$ is established, i.e. Theorem \ref{thm-bif}. 

We first discuss the well-definedness and the regularity of $\Ph$ in the following subsection. Then the bifurcation result is presented in Subsection 3.2. All proofs are presented in Subsection 3.3 and Subsection 3.4. Through out this section, we denote $\mu_\al=\rh_\al\bar\mu$.


\subsection{Well-definedness and regularity of $\Ph$}
In this subsection, assumptions are introduced, and the well-definedness and the regularity of $\Ph$ are discussed. All proofs of lemmas and corollaries in this subsection are presented in Subsection 3.3. 

Denote by $R_\th$ the range of the function $\th$.
Assume that
\beg{description}

\item {(A1)} The potentials $V_0,\bar V,V_1,V_2$ satisfy $\sup_{\th_1\in R_\th}e^{-\th_1 V_0}\in L^1$, $e^{-\bar V}\in L^1$, and
\beg{align}
\int_{\R^d}\left(|V_0|^r+|V_1|^r+e^{  \be \|V_2(x,\cdot)\|_{L^2(\bar\mu)}}\right)\bar\mu(\d x)<+\infty,~ r\geq 1, \be >0,\label{expVV}
\end{align}
and there exists a positive function $C_0$ on $R_\th\times (\hat{\al},\check{\al})\times [0,+\infty)$ so that $C_0$ is increasing in each variable and  for $\th_1\in R_\th,\be_1\in (\hat{\al},\check{\al}),\be_2\geq 0$, there is 
\beg{equation}\label{V-F3-ad}
-\th_1 V_0(x)+\be_1 |V_1(x)|+\be_2\|V_2(x,\cdot)\|_{L^2(\bar\mu)} \leq -\bar V(x)+ C_0(\th_1,\be_1,\be_2),~\text{a.e.}~ x\in\R^d.
\end{equation} 
\end{description}
\beg{rem}\label{rem:1}
Noticing that the fixed point of $ \cT(\cdot,\al)$ is a probability density,  the decomposition \eqref{VVV} can be replaced by the following form without changing fixed points of $\cT$:
\[\int_{\R^d}V(x,y)\rh(y)\mu(\d y)=V_1(x)+\int_{\R^d}\left( V_2(x,y)+K_1(y)+K_2(x,y)\right)\rh(y)\bar\mu(\d y).\]
For simplicity, we denote by  $V(x,\rh \bar\mu)$ the potential without $K_1,K_2$, i.e.
\[V(x,\rh\bar\mu)=V_1(x)+\int_{\R^d}V_2(x,y)\rh(y)\bar\mu(\d y).\] 

The condition \eqref{V-F3-ad} is borrowed from \eqref{V-F3}. Then, as proving in Lemma \ref{ps(mu)} (see \eqref{VVF3} and \eqref{b-mups}), we have for $\cT$ without $K_1,K_2$   that 
\beg{equation}\label{ine-mu-cT-K}
\beg{split}
&\bar\mu\left(e^{-\th(\al)V_0-\al V(\rh\bar\mu)+\bar V}\right) \\
&\qquad\geq \left(\int_{\R^d}e^{-\th(\al)V_0(x)}\d x\right)^{2}\left(\bar\mu\left(e^{-\th(\al)V_0+\al (|V_1|+\|V_2\|_{L_y^2}\|\rh\|_{L^2(\bar\mu)})+\bar V}\right)\right)^{-1}\\
&\qquad\geq \|e^{-\th(\al)V_0}\|_{L^1}^{2} \exp\{-C_0(\th(\al),\al,\al\|\rh\|_{L^2(\bar\mu)})\},
\end{split}
\end{equation}
and 
\beg{equation}\label{ine-cT-K}
\beg{split}
\cT(\rh,\al)=\ff {e^{-\th(\al)V_0-\al V(\rh\bar\mu)+\bar V}} {\bar\mu\left(e^{-\th(\al)V_0-\al V(\rh\bar\mu)+\bar V}\right)}\leq \|e^{-\th(\al)V_0}\|_{L^1}^{-2} \exp\{2C_0(\th(\al),\al,\al\|\rh\|_{L^2(\bar\mu)})\}.
\end{split}
\end{equation}
Hence, $\cT(\cdot,\al)$ without $K_1$ and $K_2$ is a mapping from $L^2(\bar\mu)$ to $L^\infty$. 

\end{rem}

We first introduce the following local uniqueness and regularity result on the fixed point of the mapping  $\cT$ with $K_1,K_2\equiv 0$. For a probability measure $\mu_\al$, let   
\beg{align*}
\pi_\al f= f-\mu_{\al}(f),~f\in L^1(\mu_{\al }),
\end{align*}
and let  ${\bf V}_{2,\al}$ be the integral operator in $L^2(\mu_{\al})$ induced by the kernel $V_2 $: 
\beg{align*}
{\bf V}_{2,\al}f(x)=\int_{\R^d}V_2(x,y)f(y)\mu_{\al}(\d y).
\end{align*}
Denote by $J_{\al_0,\de}=[\al_0-\de,\al_0+\de]$.
\beg{lem}\label{con-rhal}
Assume that (A1) holds except \eqref{expVV}, $K_1,K_2\equiv 0$, $V_0\in L^2(\bar\mu)$,  and $V_2\in L^2(\bar\mu\times\bar\mu)$. Suppose that at some $\al_0\in (\hat{\al},\check{\al})$, $ \cT(\cdot,\al_0)$ has a fixed point $\rh_{\al_0}\in L^2(\bar\mu)$  and $I+\al_0\pi_{\al_0}{\bf  V}_{2,\al_0}\pi_{\al_0}$ is invertible on $L^2(\mu_{\al_0})$. Then there exists $\de>0$ such that for each $\al\in J_{\al_0,\de}$, there is a unique $\rh_{\al}\in L^2(\bar\mu)$ satisfying  that $\rh_\al=\cT(\rh_\al,\al)$, $J_{\al_0,\de}\ni\al\mapsto \rh_\al$ is continuously differentiable in $L^2(\bar\mu)$, and
\beg{align}\label{sup-rh}
&\sup_{\al\in J_{\al_0,\de}}\|\rh_\al\|_\infty<+\infty,\\
\pp_\al\log\rh_{\al_0}=  -& (I+\al_0\pi_{\al_0} {\bf V}_{2,\al_0}\pi_{\al_0})^{-1} \pi_{\al_0}\left(\th'(\al_0)V_0+V(\mu_{\al_0})\right).\label{pp-logrh0}
\end{align}
If \eqref{expVV} holds furthermore, then 
\beg{align}\label{sup-pp-rh}
\sup_{\al\in J_{\al_0,\de}}|\pp_\al\log\rh_\al| \in L^r(\bar\mu),~r\geq 1,
\end{align}
and for any $r\geq 1$, $\rh_\al,\pp_\al\rh_\al,\pp_\al\log\rh_\al$ are continuous in $\al$ from $J_{\al_0,\de}$ to $L^r(\bar\mu)$. 
\end{lem}

We also assume that  $V_2,K_1,K_2$ satisfy the following conditions. 
\beg{description}

\item (A2) $K_1\in L^2(\bar\mu)$. For all $\be>0$, 
\beg{align}\label{V1-z-V}
\int_{\R^d}\exp\left\{\be \|K_2(x,\cdot)\|_{L^2(\bar\mu)}\right\}\bar\mu(\d x)<\infty.
\end{align}
There are $\ga_1,\ga_2>2$ such that $\|V_2\|_{L^2_xL^{\ga_1}_y}$ and $\|K_2\|_{L^2_xL^{\ga_2}_y}$ are finite.

\end{description}
For $\mu_\al$ given by Lemma \ref{con-rhal}, we assume that $K_2$ is orthogonal to $\{\mu_\al\}_{\al\in J_{\al_0,\de}}$, i.e.
\beg{description}
\item (A3) For almost $x\in\R^d$, 
\beg{align}\label{K2mu0}
\int_{\R^d}K_2(x,y)\mu_\al(\d y)=0,~\al\in J_{\al_0,\de}.
\end{align} 
\end{description}

\beg{rem}\label{canc-K1}
Under (A2), for $\rh\in L^2(\bar\mu)$, $\bar\mu(\rh K_1)$ is a constant, then $K_1$ can be canceled in $\cT$, see the proof of Lemma \ref{lem-ps-T} below. 
\end{rem}
\beg{rem}\label{dec-H}
Introducing a parameter to  Equation \eqref{exa-granular}:
\beg{equation}\label{al-granular}
\mu(\d x)=\ff {\exp\{-\th(\al)V_0(x)+\al\int_{\R^d}H(x-y)\mu(\d y)\}\d x} {\int_{\R^d}\exp\{-\th(\al) V_0(x)+\al \int_{\R^d}H(x-y)\mu(\d y)\}\d x}
\end{equation}
where $V_0, H\in C^1(\R^d)$  are  symmetric functions, i.e. for all $x\in\R^d$, $V_0(-x)=V_0(x)$ and $H(-x)=H(x)$.  $H(x-y)$ can be decomposed into even and odd parts. Indeed, we can set 
\[V_2(x,y)=-\ff {H(x-y)+H(x+y)} 2,\qquad K_2(x,y)=-\ff {H(x-y)-H(x+y)} 2.\]
Then  $V_2(x,y)+K_2(x,y)=-H(x-y)$ and 
\beg{align*}
V_2(x,y)=V_2(y,x),&\qquad V_2(-x,y)=V_2(x,y),\\
K_2(x,y)=K_2(y,x),&\qquad K_2(x,-y)=-K_2(x,y).
\end{align*}
In this case, $V_1=K_1\equiv 0$. Then the fixed point of $\cT(\cdot,\al)$ without $K_2$, saying $\rh_\al$, is symmetry. Taking into account that $K_2$ is anti-symmetric, we find that (A3) holds.  
\end{rem}

To illustrate (A1)-(A3), we present the following concrete example.
\beg{exa}\label{exa-xsin00}
Consider 
\beg{equation}\label{exa-xsin0}
\mu(\d x)=\ff {\exp\left\{-\ff 2 {\si^2}\left[\ff {x^4} 4-\ff {x^2} 2+\ff {\be} 2\int_{\R}\left((x-y)^2-(\sin x-\sin y)^2\right)\mu(\d y)\right]\right\}\d x} {\int_{\R}\exp\left\{-\ff 2 {\si^2}\left[\ff {x^4} 4-\ff {x^2} 2+\ff {\be} 2\int_{\R}\left((x-y)^2-(\sin x-\sin y)^2\right)\mu(\d y)\right]\right\}\d x}.
\end{equation}
We set $\al=\ff {\be} {\si^2}$,  $\th(\al)=\ff 2 {\be} \al$, $V_0(x)=\ff {x^4} 4-\ff {x^2} 2$, 
\beg{gather*}
V_1(x)=x^2-\sin^2x,~V_2(x,y)\equiv 0,\\
K_1(y)=y^2-\sin^2 y,~K_2(x,y)=-2xy+2\sin x\sin y,
\end{gather*}
For any finite interval $(\hat{\al},\check{\al})\subset (0,+\infty)$,  we have $R_{\th}=(\ff {2\hat{\al}} {\be}, \ff {2\check{\al}} {\be})$. Let $\bar V(x)=\ff {\hat{\al}} {4\be}x^4$. Then the H\"older inequality implies directly that (A1) and (A2) holds.   In this setting, \eqref{exa-xsin0} can be reformulated as follows
\beg{align}\label{dmu-VVV}
\mu(\d x)&= \ff {\exp\left\{-\ff 2 {\si^2}\left[\ff {x^4} 4-\ff {x^2} 2+\ff {\be} 2\int_{\R}\left(x^2-\sin^2 x -2xy +2\sin x\sin y\right)\mu(\d y)\right]\right\}\d x} {\int_{\R}\exp\left\{-\ff 2 {\si^2}\left[\ff {x^4} 4-\ff {x^2} 2+\ff {\be} 2\int_{\R}\left(x^2-\sin^2 x -2xy +2\sin x\sin y\right)\mu(\d y)\right]\right\}\d x}\nonumber\\
&= \ff {\exp\left\{-\ff {2\al} {\be}V_0(x)-\al V_1(x)-\al\int_{\R}K_2(x,y)\mu(\d y)\right\}\d x} {\int_{\R}\exp\left\{-\ff {2\al} {\be}V_0(x)-\al V_1(x)-\al\int_{\R}K_2(x,y)\mu(\d y)\right\}\d x}.
\end{align}
Let 
\[\mu_{\al}(\d x)=\ff {e^{-\ff {2\al} {\be} V_0(x)-\al V_1(x)}} {\int_{\R}e^{-\ff {2\al} {\be} V_0(x)-\al V_1(x)}\d x}\d x.\]
Since $V_0,V_1$ are even functions, we have that $\int_{\R}K_2(x,y)\mu_{\al}(\d y)=0$. Thus (A3) holds.
\end{exa}

 We denote by 
\beg{align*}
\Ps(x;w,\al)=\exp\left\{ -\al  \int_{\R^d}(V_2+K_2)(x,y)w(y)\mu_{\al}(\d y) \right\} .
\end{align*} 
The following lemma shows that the zero is a trivial solution  of $ \Ph(\cdot,\al)=0$, and $ \Ph(\cdot,\al)$ is continuously Fr\'echet differentiable.    
\beg{lem}\label{lem-ps-T}
Assume that (A1) holds, $K_1\in L^2(\bar\mu)$ and $K_2$ satisfies \eqref{V1-z-V}. Let $\{\rh_\al\}_{\al\in J_{\al_0,\de}}\subset  L^2(\bar\mu)$ be a family of fixed points for $\cT$ with $K_1,K_2\equiv 0$. Suppose $K_2$ is  orthogonal to $\{\rh_\al\}_{\al\in J_{\al_0,\de}}$. Then $\Ph(0,\al)=0$. 
Moreover, $\Ph(\cdot,\al)$ is continuously Fr\'echet differentiable from $L^2(\mu_\al)$ to $L^2(\bar\mu)$ with the Fr\'echet derivative given by 
\beg{equation}\label{nnPh0}
\beg{split}
\nn_{w} \Ph(w_1,\al)&= w- \ff {\Ps(w_1,\al)\log\Ps(w,\al)} {\mu_\al(\Ps(w_1,\al))}\\
&\quad + \ff {\Ps(w_1,\al)\mu_\al(\Ps(w_1,\al)\log\Ps(w,\al))} {\mu_\al(\Ps(w_1,\al))^2},~w,w_1\in L^2( \mu_\al).
\end{split}
\end{equation}
In particular, 
\beg{equation}\label{nnPh(0)} 
\nn_{w} \Ph(0,\al)=w  +\al \pi_{\al}({\bf V}_{2,\al}+{\bf K}_{2,\al})w, 
\end{equation}
where ${\bf K}_{2,\al}$ be the integral operator on $L^2(\mu_{\al})$ induced by the kernel $K_2 $: 
\beg{align*}
{\bf K}_{2,\al}f(x)=\int_{\R^d}K_2(x,y)f(y)\mu_{\al}(\d y),
\end{align*}
and $\nn \Ph(0,\al)$ is a  Fredholm operator on $L^2(\bar\mu)$ and $L^2(\mu_\al)$.
\end{lem}
If  assumptions of Lemma \ref{con-rhal} holds in addition, the reference spaces $\{L^2(\mu_\al)\}_{\al\in J_{\al_0,\de}}$ can be reduced to one, i.e. 
\beg{cor}\label{cor:mu0mu}
Assume that (A1) holds,  $K_1\in L^2(\bar\mu)$ and $K_2$ satisfies \eqref{V1-z-V}.   Suppose that at some  $\al_0\in (\hat{\al},\check{\al})$, $ \cT(\cdot,\al_0)$ without $K_1,K_2$ has a fixed point $\rh_{\al_0}\in L^2(\bar\mu)$, and $I+\al_0\pi_{\al_0}{\bf  V}_{2,\al_0}\pi_{\al_0}$ is invertible on $L^2(\mu_{\al_0})$. Assume  $K_2$ is orthogonal to $\{\rh_\al\}_{\al\in J_{\al_0,\de}}$ which is given by Lemma \ref{con-rhal}. Then,  for each $\al\in J_{\al_0,\de}$ (may be another smaller $\de$), $\Ph(\cdot,\al)$  is  continuously Fr\'echet differentiable from $L^2(\mu_{\al_0})$ to $L^2(\bar\mu)$, \eqref{nnPh0} and \eqref{nnPh(0)} hold for $w,w_1\in L^2(\mu_{\al_0})$, and $\nn \Ph(0,\al)$ is a Fredholm operator on $L^2(\mu_{\al_0})$.
\end{cor}

The following lemma is devoted to the regularity of $\Ph(0,\cdot)$ under conditions (A1)-(A3). We denote by ${\bf V}_2,{\bf K}_2$ integral operators on $L^2(\bar\mu)$ induced by the kernel $V_2(x,y)$ and $K_2(x,y)$, and $\otimes$ the tensor product on $L^2(\bar\mu)$, i.e. 
\[(f\otimes g)w = f\int_{\R^d} g(y)w(y)\bar\mu(\d y),~f,g,w\in L^2(\bar\mu).\]
\beg{lem}\label{con-nnPh}
Assume that (A1) and (A2) hold.  Suppose at some  $\al_0\in (\hat{\al},\check{\al})$, $ \cT(\cdot,\al_0)$ without $K_1,K_2$ has a fixed point $\rh_{\al_0}\in L^2(\bar\mu)$, and $I+\al_0\pi_{\al_0}{\bf  V}_{2,\al_0}\pi_{\al_0}$ is invertible on $L^2(\mu_{\al_0})$. Let $\{\rh_\al\}_{\al\in J_{\al_0,\de}}$ be given by Lemma \ref{con-rhal}, and assume that $K_2$  and $\{\rh_\al\}_{\al\in J_{\al_0,\de}}$ satisfy (A3).  Then for another smaller $\de$, $\pi_{\al}({\bf V}_{2,\al}+{\bf K}_{2,\al})\in C^1(J_{\al_0,\de};\cL_{HS}(L^2(\mu_{\al_0})))$ with
\beg{equation}\label{ppaPh}
\beg{split}
\pp_\al\left(\pi_{\al}({\bf V}_{2,\al}+{\bf K}_{2,\al})\right)& =\pi_{\al}({\bf V}_{2}+{\bf K}_{2})\cM_{\rh_\al+\al\pp_\al\rh_\al}\\
&\quad\,  -\al(\1\otimes\pp_\al\rh_\al) ({\bf V}_{2}+{\bf K}_{2})\cM_{\rh_\al},
\end{split}
\end{equation}
where $\cM_{\rh_\al+\al\pp_\al\rh_\al}$ and $\cM_{\rh_\al}$ are multiplication operators and $\1$ is the constant function with value equals to $1$.
\end{lem}

\subsection{Main result}
Let $L^{2}_{\mathbb{C}}(\mu_{\al_0})$ be the complexification of $L^2(\mu_{\al_0})$. Let $P(\al_0)$ be the eigenprojection of $-\al_0\pi_{\al_0}({\bf V}_{2,\al_0}+{\bf K}_{2,\al_0})\pi_{\al_0}$ associated with the eigenvalue $1$ in $L_{\mathbb{C}}^2(\mu_{\al_0})$ (if $1$ is an eigenvalue):
$$P(\al_0)=-\ff 1 {2\pi {\bf i}}\int_{\Ga}( -\al_0\pi_{\al_0}({\bf V}_{2,\al_0}+{\bf K}_{2,\al_0})\pi_{\al_0}-\et)^{-1}\d\et,$$
where ${\bf i}=\sq{-1}$, and $\Gamma$ is some simple and closed curve  enclosing $1$ but no other eigenvalue.  Denote 
$$\cH_0=P(\al_0)L^2_{\mathbb{C}}(\mu_{\al_0}),\qquad \cH_1=(I-P(\al_0))L^2_{\mathbb{C}}(\mu_{\al_0}).$$
Under the assumption of Corollary \ref{cor:mu0mu}, $\pi_{\al_0}({\bf V}_{2,\al_0}+{\bf K}_{2,\al_0})\pi_{\al_0}$ is a Hilbert-Schmidt operator. Then $\cH_0$ is  finite dimensional. Denote 
$$\tld A_0=-P(\al_0) \al_0\pi_{\al_0}({\bf V}_{2,\al_0}+{\bf K}_{2,\al_0})\pi_{\al_0}\Big|_{\cH_0},\qquad \tld M_0=\al_0 P(\al_0)\cM_{\pp_{\al}\log\rh_{\al_0}}\Big|_{\cH_0}.$$
Then $\tld A_0$ and $\tld M_0$ are matrices on $\cH_0$.

\beg{thm}\label{thm-bif}
Assume (A1) and (A2). Assume that there is $\al_0\in (\hat{\al},\check{\al})$ such that $ \cT$ with $K_2\equiv 0$ has a fixed point $\rh_{\al_0}\in L^2(\bar\mu)$ and $I+\al_0 \pi_{\al_0}{\bf V}_{2,\al_0}\pi_{\al_0}$ is invertible on $L^2(\mu_{\al_0})$. Let $\{\rh_{\al}\}_{\al\in J_{\al_0,\de}}$ be the unique family of fixed points for $ \cT$ without $K_1,K_2$ for some $\de>0$. Suppose that $K_2$ and $\{\rh_{\al}\}_{\al\in J_{\al_0,\de}}$ satisfy (A3).\\
If $0$ is an eigenvalue of $I+\al_0\pi_{\al_0}({\bf V}_{2,\al_0}+{\bf K}_{2,\al_0})\pi_{\al_0}$, $I_{\cH_0}+\tld M_0$ is invertible on $\cH_0$ and the algebraic multiplicity of the  eigenvalue $0$ of $(I_{\cH_0}+\tld M_0)^{-1}\left(\tld A_0^{-1}-I_{\cH_0}\right)$ is odd, then  $\al_0$ is a bifurcation point for $\Ph=0$ in the sense of Definition \ref{bif-def-00}. Consequently, $(\mu_{\al_0},\al_0)$ is a bifurcation point in the sense of Definition \ref{bif-proba} for \eqref{eq-bif}.\\
In particular, if $0$ is a semi-simple eigenvalue of $I+\al_0\pi_{\al_0}({\bf V}_{2,\al_0}+{\bf K}_{2,\al_0})\pi_{\al_0}$ with odd algebraic multiplicity and $I_{\cH_0}+\tld M_0$ is invertible on $\cH_{0}$, then $ \al_0$ is a bifurcation point in both sense.
\end{thm}

The proof of this theorem is presented in Subsection 3.4.  For symmetric $V_2$ and $K_2$, we use the last case in Theorem \ref{thm-bif}. As an application, we investigate 
\beg{equation}\label{eq-HVK}
\mu(\d x)=\ff {\exp\left\{-\ff 2 {\si^2} \left(V_0(x)+ \ff {\be} 2\int_{\R^d}(V_2(x,y)+K_2(x,y))\mu(\d y)\right)\right\}} {\int_\R \exp\left\{-\ff 2 {\si^2} \left(V_0(x)+ \ff {\be} 2\int_{\R^d} (V_2(x,y)+K_2(x,y))\mu(\d y)\right)\right\}\d x}\d x.
\end{equation} 
Assume that $V_0\in C^1$ is symmetric $(V_0(-x)=V_0(x))$, $V_2$ and $K_2$ are of the following form  
\beg{align}\label{ad-VVKK}
V_2(x,y)=\sum_{i,j=1}^{l}J_{ij}v_i(x)v_j(y),\qquad K_2(x,y)=\sum_{i,j=1}^{m}G_{ij}k_i(x)k_j(y)
\end{align}
where the matrices $J=(J_{ij})_{1\leq i,j\leq l}$ and $G=(G_{ij})_{1\leq i,j\leq m}$ are symmetric, $\{\1\}\cup\{v_i\}_{i=1}^l$ is linearly independent and symmetric ($v_i(-x)=v_i(x)$),  and $\{k_i\}_{i=1}^{m}$  is linearly independent and anti-symmetric ($k_i(-x)=-k_i(x)$).  Fix $\be>0$. We set $\al=\ff {\be} {\si^2}$ and $\th(\al)=\ff {2\al} {\be}$.  We also assume that $V_0$, $\{v_i\}_{i=1}^l$ and $\{k_i\}_{i=1}^{m}$  are continuous and satisfy the conditions in Corollary \ref{exa-singular0}. In this case, our criteria for the bifurcation point is presented by using characteristics of matrices.  We remark here that if the kernel $H(x-y)$ in Remark \ref{dec-H} induces a finite rank operator,  then, according to Remark \ref{dec-H} and \cite[Theorem 1.4]{Simon}, the even and odd parts of $H(x-y)$ can be represented in the form of \eqref{ad-VVKK}. For a matrix $A$, we denote by $A^*$ the transpose of $A$. Our criteria for the bifurcation point of \eqref{eq-HVK} with \eqref{ad-VVKK} reads as follows.

\begin{cor}\label{cor-finite}
Fix $\be>0$. Let $V_0,V_2,K_2,J,G$, $\{v_i\}_{i=1}^l$ and $\{k_i\}_{i=1}^{m}$ be stated as above. Assume that $J$ is nonnegative  and  the following equation has a solution at some $\si_0>0$
\beg{equation}\label{nu-si0}
\nu_{\si}(\d x)=\ff {\exp\left\{-\ff 2 {\si^2} \left(V_0(x)+ \ff {\be} 2\sum_{i,j=1}^l J_{ij} v_i(x)\nu_{\si}(v_j)\right)\right\}} {\int_\R \exp\left\{-\ff 2 {\si^2} \left(V_0(x)+ \ff {\be} 2\sum_{i,j=1}^l J_{ij} v_i(x)\nu_{\si}(v_j)\right)\right\}\d x}\d x.
\end{equation}
Denote by $\nu_{\si_0}(\d x)$ the solution of \eqref{nu-si0} at $\si_0$. Then there is  $\de>0$ such that $(\nu_{\si_0},\si_0)$ can be extended uniquely to a path : $\si\in [\si_0-\de,\si_0+\de]\mapsto (\nu_\si,\si)$ which satisfies \eqref{nu-si0}. Moreover, $(\nu_\si,\si)$ also satisfies \eqref{eq-HVK}.\\
Let $G(\si_0)\equiv(G_{ij}(\si_0))_{1\leq i,j\leq m}:=(\nu_{\si_0}(k_ik_j))_{1\leq i,j\leq m}$, $\al_0=\ff {\be} {\si_0^2}$ and $\mu_{\al_0}=\nu_{\si_0}$. Suppose\\
(1) $0$ is an eigenvalue of the matrix $I+\al_0GG(\si_0)$ with odd algebraic multiplicity, \\
(2) the matrices $G(\si_0)$, $J$ and $G$ satisfy 
\beg{equation}\label{eq-rak}
\beg{split}
&{\rm rank}\left(\left[\begin{array}{cc}
0 & I+\al_0GG(\si_0)\\
I+\al_0GG(\si_0) & -(I+\al_0G(\si_0)^{-1}M_K(\al_0))
\end{array}
\right]\right)\\
&\quad\, =m+{\rm rank}(I+\al_0GG(\si_0)),
\end{split}
\end{equation}
where ${\rm rank}(\cdot)$ is the rank of a matrix,  $M_K(\al_0)=\left(\mu_{\al_0}(\pp_\al\log\rh_{\al_0}k_ik_j)\right)_{1\leq i,j\leq m}$,
\begin{align*}
\pp_\al\log\rh_{\al_0}&=-\sum_{i=1}^l\left[(I+\al_0JJ(\al_0))^{-1}\left(\ff 2 {\be} w+\tld w\right)\right]_{i}\pi_{\al_0}(v_i)\\
&\quad - \ff 2 {\be}\left( \pi_{\al_0}V_0-\sum_{i=1}^l w_i\pi_{\al_0}(v_i)\right),\\
J(\al_0)&=(\mu_{\al_0}( v_i v_j )-\pi_{\al_0}( v_i)\pi_{\al_0}( v_j))_{1\leq i,j\leq l},
\end{align*}
and $w,\tld w$ are column vector such that
\beg{align*}
w^*&=(w_i)_{1\leq i\leq l}= \left(  \sum_{j=1}^l (J(\al_0)^{-1})_{ij}\mu_{\al_0}( V_0 \pi_{\al_0}v_j )\right)_{1\leq i\leq l},\\
\tld w^*&=(\tld w_i)_{1\leq i\leq l}=(\sum_{j=1}^lJ_{ij}\mu_{\al_0}(v_j))_{1\leq i\leq l}.
\end{align*}
Then $\al_0$ (or $\si_0$) is a bifurcation point in the sense of Definition \ref{bif-def-00}, which implies that $(\nu_{\si_0},\si_0)$ is a bifurcation point for \eqref{eq-HVK} in the sense of Definition \ref{bif-proba}.
\end{cor}
\beg{rem}
When $V_2=0$ in \eqref{eq-HVK}, i.e. $J=0$ or  $v_i=0$ for $i=1,\cdots,l$, we can just let $M_K(\al_0)=-\ff 2 {\be}\pi_{\al_0}V_0$ (just set $J=0$ and $v_i=0$ in the definition of $M_{K}(\al_0)$) without changing the assertion in this Corollary.  
\end{rem}
The proof of this corollary is presented at the end of Subsection 3.4.   We give the following example to illustrate this corollary and the theorem. 

\beg{exa}
Consider Example \ref{exa-xsin00} with $\be=2$. Due to \eqref{dmu-VVV}, we investigate
\beg{equation}\label{exa-xsin}
\mu(\d x)=\ff {\exp\left\{-\ff 2 {\si^2}\left[\ff {x^4} 4+\ff {x^2} 2-\sin^2 x-2\int_{\R}\left(xy-\sin x\sin y\right)\mu(\d y)\right]\right\}\d x} {\int_{\R}\exp\left\{-\ff 2 {\si^2}\left[\ff {x^4} 4+\ff {x^2} 2-\sin^2x-2\int_{\R}\left(xy-\sin x\sin y\right)\mu(\d y)\right]\right\}\d x}.
\end{equation}
In this case, we set $\al=\ff 2 {\si^2}$,  
\beg{gather*}
V_0(x)=\ff {x^4} 4+\ff {x^2} 2-\sin^2 x,\qquad  k_1(x)=x,\qquad  k_2(x)=\sin x,\\[4pt]
\mu_{\al}(\d x)=\ff {e^{-\al V_0(x)}} {\int_{\R}e^{-\al V_0(x)}\d x}\d x,\qquad G=\left[\begin{array}{cc} -2&0\\ 0& 2\end{array}\right],\\[4pt]
G(\al)=\left[\begin{array}{cc} \mu_{\al}(k_1^2)&\mu_{\al}(k_1k_2)\\[3pt] \mu_{\al}(k_1k_2)& \mu_{\al}(k_2^2)
\end{array}\right], \qquad M_K(\al)=-\left[\begin{array}{cc} \mu_{\al}(k_1^2\pi_{\al} V_{0})&\mu_{\al}(k_1k_2\pi_{\al} V_{0})\\[4pt] \mu_{\al}(k_1k_2\pi_{\al} V_{0})& \mu_{\al}(k_2^2\pi_{\al} V_{0})
\end{array}\right],
\end{gather*}
If there is $\si_0>0$ such that $\al_0:=\ff 2 {\si_0^2}$ satisfies 
\beg{gather}\label{equ-alal}
1-2\al_0\mu_{\al_0}(k_1^2-k_2^2)-4\al_0^2\left(\mu_{\al_0}(k_1^2)\mu_{\al_0}(k_2^2)-(\mu_{\al_0}(k_1k_2))^2\right)=0,\\[5pt]
{\rm rank}\left(\left[\begin{array}{cc}
0 & I+\al_0GG(\al_0)\\
I+\al_0GG(\al_0) & -(I+\al_0G(\al_0)^{-1}M_{K}(\al_0))
\end{array}
\right]\right)=3,\label{mat-GGM}
\end{gather}
then $\al_0$ is  a bifurcation point for \eqref{exa-xsin}.
\end{exa}
\beg{proof}
According to  \eqref{exa-xsin}, we can set $\th(\al)=\al$, $V_1=0,V_2=0$, and for any  $0<\hat{\al}<\al_0<\check{\al}$, 
\[\bar V(x)= \hat{\al}\left(\ff {x^4} 4+\ff {x^2} 2\right).\]
Then $R_{\th}=(\hat{\al}, \check{\al})$ and 
\[-\th_1 V_0(x)\leq -\bar V(x)+\th_1 \sin^2x\leq -\bar V(x)+\th_1,~\th_1\in (\hat{\al}, \check{\al}).\]
Hence, (A1) holds.  We can apply Corollary \ref{cor-finite} directly to this example with $V_2= 0$, $\be=2$, and $V_0,k_1,k_2,G,$ e.t.c  defined accordingly as above. Next, we check (1) and (2) in Corollary \ref{cor-finite}. It follows from \eqref{equ-alal} that  
\beg{align*}
{\rm det}(I+\al_0GG(\al_0))&=1-2\al_0\mu_{\al_0}(k_1^2-k_2^2)-4\al_0^2\left(\mu_{\al_0}(k_1^2)\mu_{\al_0}(k_2^2)-(\mu_{\al_0}(k_1k_2))^2\right)\\
&=0.
\end{align*}
Notice that the trace of the matrix $I+\al_0GG(\al_0)$:
\beg{align*}
\trac(I+\al_0GG(\al_0))&=2-2\al_0\mu_{\al_0}(k_1^2-k_2^2)=1+1-2\al_0\mu_{\al_0}(k_1^2-k_2^2)\\
&=1+4\al_0^2\left(\mu_{\al_0}(k_1^2)\mu_{\al_0}(k_2^2)-(\mu_{\al_0}(k_1k_2))^2\right)\\
&\geq 1.
\end{align*}
We see that $0$ is  an eigenvalue of the matrix $I+\al_0GG(\al_0)$ with algebraic multiplicity one. Consequently, \eqref{eq-rak} holds due to \eqref{mat-GGM}, ${\rm rank}(I+\al_0GG(\al_0))=1$ and $m=2$.

\end{proof}

\beg{rem}
The choice of $V_1,V_2,K_1,K_2$ may be not unique. Indeed, comparing with Example \ref{exa-xsin00}, we choose different $V_0,V_1$ in this example. 

We solve \eqref{equ-alal} numerically and get $\al_0\approx 5.94468752$ by using the function {\rm fsolve$()$} from {\rm scipy.optimize} of Python with the initial value $=4$,  and the solution converges.  The same approximate value of $\al_0$ is also obtained by using {\rm optimize.root()} with  the initial value $=4$ and method $=$ 'lm'. For this approximate value of $\al_0$, we have
\beg{align*}
&G(\al_0)\approx\left[\begin{array}{cc} 0.39618539&0.35382333\\[3pt] 0.35382333& 0.31704574
\end{array}\right],~~M_K(\al_0)\approx\left[\begin{array}{cc} -0.02892554&-0.02387658\\[3pt] -0.02387658& -0.01979521
\end{array}\right],\\[5pt]
&\left[\begin{array}{cc}
0 & I+\al_0GG(\al_0)\\[3pt]
I+\al_0GG(\al_0) & -(I+\al_0G(\al_0)^{-1}M_{K}(\al_0))
\end{array}
\right]\\[5pt]
&\qquad~\,\approx \left[\begin{array}{cccc}
0 & 0 & -3.71039667 & -4.20673828\\
0 & 0 & 4.20673828  & 4.76947575\\
-3.71039667 & -4.20673828 & 9.27847371 & 8.05002984\\
4.20673828  & 4.76947575 & -11.02309397 & -9.61267519
\end{array}
\right],
\end{align*}
and \eqref{mat-GGM} holds.   The rank of the matrix is computed by using the function {\rm matrix$\_$rank$()$} from {\rm numpy.linag} of Python. 

\end{rem}
\beg{rem}
For \eqref{nu-si0} with $J$ and $\{v_i\}_{i=1}^l$ are not trivial zeros, we can also numerically solve the following problem on $\R^{l+1}$ to obtain  approximations of $\si_0$ and $\nu_{\si_0}$:
\[r_i=\ff {\int_{\R^d} v_i(x)\exp\left\{-\ff 2 {\si_0^2} \left(V_0(x)+ \ff {\be} 2\sum_{i,j=1}^l J_{ij} r_jv_i(x)\right)\right\}\d x} {\int_\R \exp\left\{-\ff 1 {\si_0^2} \left(V_0(x)+ \ff {\be} 2\sum_{i,j=1}^l J_{ij}r_j v_i(x)\right)\right\}\d x},~i=1,\cdots,l. \]
\end{rem}

Finally, we revisit \eqref{Daw-fixp} and \eqref{Ham-exp}  by using Theorem \ref{thm-bif}.
\beg{exa}\label{Daw-exa}
Consider \eqref{Daw-fixp}. Fix $\be$. Let 
\beg{equation}\label{nu_si}
\nu_{\si}(\d x)=\ff { \exp\left\{- \ff 2 {\si^2}  \left(\ff {x^4} 4-\ff {x^2} 2+\ff {\be} 2  x^2\right) \right\}} {\int_\R \exp\left\{- \ff 2 {\si^2}  \left(\ff {x^4} 4-\ff {x^2} 2+\ff {\be} 2  x^2\right)   \right\} \d x}\d x.
\end{equation}
Then $\nu_{\si}$ is a stationary probability measure for \eqref{Daw-fixp}. If there is $\si_0>0$ such that 
\begin{equation}\label{2nusi}
1=\ff {2\be} {\si_0^2}\int_{\R}x^2\nu_{\si_0}(\d x),
\end{equation}
then $\si_0\in [\sq{\ff {2\be} 3}, \sq{2\be}]$, and $\si_0$ is a bifurcation point, i.e. for any $\de'>0$, there is $(\mu_{\si},\si)$ with $\si \in (\si_0-\de',\si_0+\de')$ so that  $(\mu_{\si},\si)$ satisfy \eqref{Daw-fixp} and $\mu_{\si}\neq \nu_{\si}$.
\end{exa}
\beg{proof}
Choose $0<\hat\si<\si_0<\check{\si}<+\infty$, and set $\hat\al =\ff {2\be} {\check{\si}^2}$, $\check{\al}=\ff {2\be} {\hat\si^2}$, $\al=\ff {2\be} {\si^2}$, $\th(\al)=\ff {\al} {\be}$, and
\[V_0(x)=\ff {x^4} 4-\ff {x^2} 2,~V_1(x)=\ff {x^2} 2,~V_2(x,y)=0,~K_1(y)=y^2,~K_2(x,y)=-xy.\]
Let  $\bar V(x)=\ff {\hat\al} {8\be} x^4$. Then the H\"older inequality yields that
\beg{align*}
-\th V_0(x)+\be_1\ff {x^2} 2&\leq -\ff {\hat\al} {4\be} x^4+\left( \th +\be_1\right)\ff {x^2} 2\\
&\leq -\ff {\hat\al} {4\be} x^4+\ff {\hat\al} {8\be} x^4+\ff {\be } {2 \hat\al}\left( \th +\be_1\right)^2\\
&=-\bar V(x)+\ff {\be } {2\hat\al}\left( \th +\be_1\right)^2,~\th\in (\ff {\hat\al} {\be},\ff {\check{\al}} {\be}),~\be_1>0.
\end{align*}
It is clear that $K_1\in L^2(\bar\mu)$. Thus $\cT$ is of the following form 
\beg{equation}\label{tldT-Daw}
\cT(x;\rh,\al)= \ff {\exp\left\{-\ff {\al}  {\be} \left(\ff {x^4} 4-\ff {x^2} 2\right) - \al \left(\ff {x^2} 2- x \int_\R y\rh(y)\bar\mu(\d y)\right)+\bar V(x)\right\}} {\int_\R \exp\left\{-\ff {\al}  {\be} \left(\ff {x^4} 4-\ff {x^2} 2\right) - \al \left(\ff {x^2} 2-x \int_\R y\rh(y)\bar\mu\d y)\right)+\bar V(x)\right\}\bar\mu(\d x)},
\end{equation}
and (A1) and (A2)  hold.  For $\cT$ without $K_2$, we have that
\beg{align*}
\rh_{\al}(x)&=\ff {\exp\left\{-\ff {\al}  {\be} \left(\ff {x^4} 4-\ff {x^2} 2\right) - \ff {\al} 2  x^2 +\bar V(x)\right\}} {\int_\R \exp\left\{-\ff {\al}  {\be} \left(\ff {x^4} 4-\ff {x^2} 2\right) -\ff {\al} 2  x^2 +\bar V (x)\right\}\bar\mu(\d x)},\\
\nu_\si(f)&=\bar\mu(\rh_\al f)=\mu_{\al}(f),~f\in\sB_b(\R^d),
\end{align*}
and 
\[I+\al \pi_{\al}{\bf V}_{2,\al}\pi_{\al}\equiv I,~\al\in (\hat\al,\check\al).\]
Noting that $K_2(x,y)$ is symmetric,  eigenvalues of $I+\al\pi_\al({\bf V}_{2,\al}+{\bf K}_{2,\al})\pi_\al$ are semisimple.  According to Theorem \ref{thm-bif}, we need to show that  $0$ is an eigenvalue of $I+\al\pi_{\al_0}({\bf V}_{2,\al_0}+{\bf K}_{2,\al_0})\pi_{\al_0}$ with odd algebraic multiplicity and $I_{\cH_0}+\tld M_0$ is invertible. Since $\rh_\al$ is an even function and $\bar\mu$ is a symmetric measure, it is clear that (A3) holds. For ${\bf K}_{2,\al_0}$, we have for any $f\in L^2(\mu_{\al_0})$
\[{\bf K}_{2,\al_0} \pi_{\al_0} f(x)=-x \int_{\R}(\pi_{\al_0} f)(y)y \mu_{\al_0}(\d y)=-x \int_{\R}f(y)y \mu_{\al_0} (\d y),\]
and $\pi_{\al_0} {\bf K}_{2,\al_0} \pi_{\al_0}f={\bf K}_{2,\al_0}f$. Thus,  $0$ is an eigenvalue of $I+\al_0\pi_{\al_0}({\bf V}_{2,\al_0}+{\bf K}_{2,\al_0})\pi_{\al_0}$ if and only if $\al_0$ satisfies
\[
x-\al_0 x\int_{\R} y^2\mu_{\al_0}(\d y)=0,~x\in\R,
\]
equivalently, 
\beg{align}\label{eq:al0}
1=\al_0\int_{\R}y^2\mu_{\al_0}(\d y)\equiv\ff {2\be} {\si_0^2}\int_{\R}x^2\nu_{\si_0}(\d x).
\end{align}
This also implies that  $\{\sq{\al_0} x \}$ is a orthonormal basis of $\cH_0$ and 
\[P(\al_0)f(x)=\sq{\al_0} x  \int_{\R} (\sq{\al_0} y )f(y)\mu_{\al_0}(\d y)= \al_0 x\int_{\R}yf(y)\mu_{\al_0}(\d y).\]
For $\tld M_0$,
\beg{align*}
\pp_\al\log\rh_{\al_0}(x)&= -\ff 1 {\be}\left(\ff {x^4} 4-\ff {x^2} 2\right)-\ff {x^2} 2+\int_{\R} \left(\ff 1 {\be}\left(\ff {x^4} 4-\ff {x^2} 2\right)+\ff {x^2} 2\right)\mu_{\al_0}(\d x)\\
&= -\ff {x^4} {4\be} +\ff {1-\be} {2\be}x^2+\ff {m_4} {4\be}-\ff {1-\be} {2\be} m_2,
\end{align*}
where 
\[m_4=\int_{\R}x^4\mu_{\al_0}(\d x),\qquad m_2=\int_{\R}x^2\mu_{\al_0}(\d x).\]
Thus
\beg{equation}\label{tldM0-0}
\beg{split}
\tld M_0 &=\al_0\int_{\R}\left(-\ff {x^4} {4\be} +\ff {1-\be} {2\be}x^2+\ff {m_4} {4\be}-\ff {1-\be} {2\be} m_2\right)  (\al_0 x^2 ) \mu_{\al_0}(\d x)\\
&= \al_0^2\left(-\ff {m_6} {4 \be}+\ff {1-\be} {2 \be} m_4+\ff {m_4 m_2} {4 \be}-\ff {1-\be} {2 \be}m_2^2\right).
\end{split}
\end{equation}
It is clear that $\mu_{\al_0}$ is the unique invariant probability measure of the following SDE
\[\d X_t=-(X_t^3-X_t)\d t-\be X_t\d t+\si_0\d W_t.\]
It follows from the It\^o formula that
\beg{align*}
X_t^2-X_0^2&=\int_0^t(-2X_s^4+2(1-\be)X_s^2 + \si_0^2)\d s+2\int_0^t X_s\si_0\d W_s\\
X_t^4-X_0^4&=\int_0^t(-4X_s^6+4(1-\be)X_s^4+6 \si_0^2 X_s^2)\d s+4\int_0^t X_s^3\si_0\d W_s.
\end{align*}
Choosing $X_0\overset{d}{=}\mu_{\al_0}$, we find that
\beg{align*}
0&=-2 m_4+2(1-\be)m_2+\si_0^2,\\
0&= -4m_6+4(1-\be)m_4+6\si_0^2 m_2.
\end{align*}
Putting these with \eqref{eq:al0}, which yields $m_2=\al_0^{-1}$, into \eqref{tldM0-0} and taking into account $\si_0^2=\ff {2\be} {\al_0}$, we arrive at
\beg{equation}\label{tldM0-1}
\beg{split}
\tld M_0&=\al_0^2\left( \ff {(1-\be)^2- \si_0^2} {4 \be}m_2 -\ff {1-\be } {4 \be}m_2^2+\ff {(1-\be)\si_0^2} {8 \be} \right)\\
&= \ff {\al_0(1-\be)-(1+\be)} {4 \be}.
\end{split}
\end{equation}
The Jensen inequality, $m_2=\al_0^{-1}$ and $\si_0^2=\ff {2\be} {\al_0}$ imply that 
\beg{align*}
\ff 1 {\al_0^2}=m_2^2\leq m_4=(1-\be)m_2+\ff {\si_0^2} 2= \ff {1-\be} {\al_0}+\ff {\be} {\al_0}=\ff 1 {\al_0}=m_2.
\end{align*}
This yields that $\al_0\geq 1$ and $m_6=((1-\be) +\ff 3 2\si_0^2)m_2$. Due to the Hankel inequality, see e.g. \cite[(3.33)]{Daw} or from the nonnegative definiteness of the moment matrix: 
\begin{align*}
\left[\begin{array}{cccc} 
m_0 & m_1 & m_2 & m_3\\
m_1 & m_2 & m_3 & m_4\\
m_2 & m_3 &m _4 & m_5\\
m_3 & m_4 &m _5 &m_6
\end{array}\right]
=\int_{\R}\left[\begin{array}{c} 
1 \\
x \\
x^2\\
 x^3
\end{array}\right]\left[\begin{array}{cccc} 
1 & x & x^2 & x^3
\end{array}\right] \mu_{\al_0}(\d x),
\end{align*}
we have that  
\beg{equation}\label{Hank}
m_2m_4m_6-m_4^3+m_2^2m_4^2-m_2^3m_6\geq 0.
\end{equation}
Combining this with $\si_0^2=2\be\al_0^{-1}$, $m_4=m_2=\al_0^{-1}<1$ and $m_6=((1-\be) +\ff 3 2\si_0^2)m_2$, we derive from \eqref{Hank} that 
\[0\leq (m_6-m_2)(1-m_2)=(\ff 3 2\si_0^2-\be)\ff {\al_0-1} {\al_0^2}=\ff {(\al_0-1)(3-\al_0)\be} {\al_0^3}.\]
We find that $\al_0=\ff {2\be} {\si_0^2}\in [1,3]$, and  
\[1+\tld M_0=  \ff {3\be-1+\al_0(1-\be) } {4 \be}=\ff {(3-\al_0)\be+(\al_0-1)} {4\be}>0.\] 

\end{proof} 

Base on Example \ref{Daw-exa}, we have the following example for Vlasov-Fokker-Planck equation. 
\beg{exa}\label{bif-Hamil}
Consider \eqref{Ham-exp}. Fix $\be>0$. Suppose that there is $\si_0>0$ satisfying  \eqref{2nusi}.  Then  $\si_0\in [\sq{\ff {2\be} 3}, \sq{2\be}]$ and  is a bifurcation point for \eqref{Ham-exp}.
\end{exa}
\beg{proof}
We first remark that if $ \nu_{1,\si}(\d x)$ is a fixed point of \eqref{Daw-fixp}, 
then 
\[\nu_{\si}(\d x,\d y):=\ff {e^{-\ff {y^2} {\si^2} }} {\sq{\pi\si^2}} \nu_{1,\si}(\d x)\d y\] 
is a solution of  \eqref{Ham-exp}. Thus the assertion of this example follows from Example \ref{Daw-exa}.

We can also repeat the proof of Example \ref{Daw-exa}. Let $e_0(x,y)=x$ and $\al_0=\ff {2\be} {\si_0^2}$. Then $\cH_0={\rm span}[\sq{\al_0} e_0]$, $P(\al_0)$ is the  orthogonal projection from $L^2(\nu_{\si_0})$ to $\cH_0$, and $\tld M_0$ equals to the $\tld M_0$ in Example \ref{Daw-exa}.

\end{proof}

\subsection{Proofs of lemmas and corollaries in Subsection 3.1}

\beg{proof}[Proof of Lemma \ref{con-rhal}]
Throughout the proof of this lemma, we assume that $K_1,K_2\equiv 0$. This lemma is proved according to the implicit function theory, see e.g. \cite[Theorem 15.1 and Theorem 15.3]{Deim}. We first investigate the regularity of $\cT$.   For any $\rh_1,\rh_2\in L^2(\bar\mu)$ and $\al_1,\al_2\in (\hat{\al},\check{\al})$, we derive from \eqref{V-F3-ad}  and \eqref{funi-1} that
\beg{equation}\label{ad-ttrh}
\beg{split}
&\left|e^{-\th(\al_1)V_0 -\al_1 V( \rh_1\bar\mu)+\bar V }-e^{-\th(\al_2)V_0 -\al_2 V( \rh_2\bar\mu)+\bar V }\right|\\
&\leq  \left(1\we C_{1,2}\right) \max_{i=1,2}\left\{e^{-\th(\al_i)V_0+\bar V+|\al_iV(\rh_i\bar\mu)| }\right\}\\
&\leq  (1\we C_{1,2})\exp\left\{\max_{i=1,2}\left\{C_0\left(\th(\al_i), \al_i,\al_i\|\rh_i\|_{L^2(\bar\mu)} \right)\right\}\right\},
\end{split}
\end{equation}
where $C_{1,2}$ is a positive function on $\R^d$ defined as follows
\beg{align*}
C_{1,2}(x)&=|\th(\al_1)-\th(\al_2)||V_0(x)|+ |\al_1-\al_2||V_1(x)|\\
&\quad +  \|\al_1\rh_1 -\al_2\rh_2\|_{L^2(\bar\mu)} \|V_2(x,\cdot)\|_{L^2(\bar\mu)}.
\end{align*} 
Combining  \eqref{ad-ttrh} with \eqref{ine-cT-K}, we have that
\beg{equation}\label{CT-CT}
\beg{split}
&|\cT(\rh_1,\al_1)-\cT(\rh_2,\al_2) |\\
&\leq  
 \left(C_{1,2}\we 1+\bar\mu(C_{1,2}\we 1)\|e^{-\th(\al_2)V_0}\|_{L^1}^{-2} e^{2C_0(\th(\al_2),\al_2,\al_2\|\rh_2\|_{L^2(\bar\mu)})}\right) \\
&\quad \times e^{ C_0(\th(\al_1),\al_1,\al_1\|\rh_1\|_{L^2(\bar\mu)})} \|e^{-\th(\al_1)V_0}\|_{L^1}^{-2}  \max_{i=1,2}\left\{e^{C_0(\th(\al_i),\al_i,\al_i\|\rh_i\|_{L^2(\bar\mu)})}\right\}.
\end{split}
\end{equation}
This, together with $\sup_{\th_1\in R_{\th}}e^{-\th_1 V_0}\in L^1$ and the dominated convergence theorem,  yields that   $\cT$ is locally bonded from $L^2(\bar\mu)\times (\hat{\al},\check{\al})$ to $L^\infty(\bar\mu)$ and is continuous from $L^2(\bar\mu)\times (\hat{\al},\check{\al})$ to $L^r(\bar\mu)$ for any $r\geq 1$. 

For any $w,\rh\in L^2(\bar\mu)$ and $N\in\N$, we derive from \eqref{ine-mu-cT-K}, \eqref{ine-cT-K} and the inequality 
$$x\leq e^{x},~x\geq 0,$$
that 
\beg{equation}\label{ine-TV1}
\beg{split}
&\sup_{s\in[-N,N]}\left|\cT(x;\rh+sw ,\al) \left({\bf V}_2w\right)(x)\right| \\
&\leq \left(\sup_{s\in[-N,N]}\left|\cT(x;\rh+sw ,\al)\right|\right)\|V_2(x,\cdot)\|_{L^2(\bar\mu)}\|w\|_{L^2(\bar\mu)}\\
&\leq \left(\sup_{s\in[-N,N]}\left|\cT(x;\rh+sw ,\al)\right|\right)e^{\|V_2(x,\cdot)\|_{L^2(\bar\mu)}\|w\|_{L^2(\bar\mu)}}\\
&\leq \|e^{-\th(\al)V_0}\|_{L^1}^{-2} e^{ 2C_0(\th(\al ),\al ,\al(\|\rh\|_{L^2(\bar\mu)}+(N+1)\|w\|_{L^2(\bar\mu)})}.
\end{split}
\end{equation}
This yields that 
\beg{equation}\label{ad-ine-TV1}
\beg{split}
&\sup_{s\in[-N,N]}\left|\cT(x;\rh+sw ,\al) \bar\mu\left(\cT(\rh+sw,\al){\bf V}_2w\right)\right|\\
 &\leq \|e^{-\th(\al)V_0}\|_{L^1}^{-4} e^{ 4C_0(\th(\al ),\al ,\al(\|\rh\|_{L^2(\bar\mu)}+(N+1)\|w\|_{L^2(\bar\mu)})}.
\end{split}
\end{equation}
Consequently, the dominated theorem theorem implies  that $\cT(\cdot,\al)$ is G\^ateaus differentiable and 
\beg{align*}
\pp_w\cT(\rh,\al)=-\al\cT(\rh,\al)\left({\bf V}_2 w-\bar\mu(\cT(\rh,\al){\bf V}_2w)\right).
\end{align*}
We also have by \eqref{CT-CT} that there is $C(\al_1,\al_2,\|\rh_1\|_{L^2(\bar\mu)},\|\rh_2\|_{L^2(\bar\mu)})>0$ which is locally bounded for $\al_1,\al_2,\|\rh_1\|_{L^2(\bar\mu)},\|\rh_2\|_{L^2(\bar\mu)}$ so that
\beg{equation}\label{ine-TV2}
\beg{split}
&\left| \left(\cT(x;\rh_1,\al_1)-\cT(x;\rh_2,\al_2) \right)({\bf V}_2 w)(x)\right|\\
&\qquad \leq C_{\al_1,\al_2,\|\rh_1\|_{L^2(\bar\mu)},\|\rh_2\|_{L^2(\bar\mu)}}\left(C_{1,2}\we 1+\bar\mu(C_{1,2}\we 1)\right)\|w\|_{L^2(\bar\mu)}. 
\end{split}
\end{equation}
This  implies that $\pp\cT$ is continuous from $L^2(\bar\mu)\times (\hat{\al},\check{\al})$ to $\sL(L^2(\bar\mu))$. Thus $\cT(\cdot,\al)$ is continuously Fr\'echet  differentiable on $L^2(\bar\mu)$ with the Fr\'echet derivative $\nn\cT$ continuous from $L^2(\bar\mu)\times (\hat{\al},\check{\al})$ to $\sL(L^2(\bar\mu))$. \\
Similarly, we can derive from $\th'\in C(\hat{\al},\check{\al})$, \eqref{V-F3-ad}, \eqref{ine-mu-cT-K}, \eqref{ine-cT-K} and \eqref{ad-ttrh} that  if $V_0\in  L^r(\bar\mu)$ for some $r\geq 1$, then $\cT(\rh,\cdot)$ is differentiable from $(\hat{\al},\check{\al})$ to $L^r(\bar\mu)$ and 
\beg{align*}
\pp_\al \cT(\rh,\al)=-\cT(\rh,\al)\left(\th'(\al)V_0+V(\rh\bar\mu)-\bar\mu(\cT(\rh,\al)(\th'(\al)V_0+V(\rh\bar\mu)))\right)
\end{align*}
which is also continuous from $L^2(\bar\mu)\times (\hat{\al},\check{\al})$ to $L^r(\bar\mu)$ for any $r\geq 1$.

Let $\tld\Ph(\rh,\al)=\rh-\cT(\rh,\al)$. From the regularity of $\cT$ and $V_0\in L^2(\bar\mu)$, we find that $\tld\Ph$ is continuously differentiable on $L^2(\bar\mu)\times  (\hat{\al},\check{\al})$ with
\beg{align*}
\nn  \tld\Ph(\rh,\al) & = I+\al\cT(\rh,\al){\bf V}_2  - \al \cT(\rh,\al)\bar\mu\left(\cT(\rh,\al){\bf V}_2\cdot \right),\\
\pp_\al \tld\Ph(\rh,\al)&=\cT(\rh,\al)\left(\th'(\al)V_0+V(\rh\bar\mu)-\bar\mu(\cT(\rh,\al)(\th'(\al)V_0+V(\rh\bar\mu)))\right).
\end{align*} 
In particular, 
\beg{align*}
\nn_w   \tld\Ph(\rh_{\al_0},\al_0) &=w+\al_0\rh_{\al_0}\left({\bf V}_2w - \mu_{\al_0}\left({\bf V}_2w\right)\right)\\
& = w +\al_0 \rh_{\al_0} \pi_{\al_0}{\bf V}_2w.
\end{align*}
Due to $\rh_{\al_0}=\cT(\rh_{\al_0},\al_0)$ and $V_2\in L^2(\bar\mu\times\bar\mu)$, we find that $\rh_{\al_0}\pi_{\al_0}{\bf V}_2$ is an integral operator on $L^2(\bar\mu)$ with a kernel in $L^2(\bar\mu\times\bar\mu)$. Then $\rh_{\al_0}\pi_{\al_0}{\bf V}_2$ is a Hilbert-Schmidt operator on $L^2(\bar\mu)$, and  $\nn\tld\Ph(\rh_{\al_0},\al_0)$ is a Fredholm operator on $L^2(\bar\mu)$.  Thus, on $L^2(\bar\mu)$, $\nn\tld\Ph(\rh_{\al_0},\al_0)$ is invertible if and only if ${\rm Ker}\left(\nn\tld\Ph(\rh_{\al_0},\al_0)\right)=\{0\}$.  For $w\in {\rm Ker}\left(\nn\tld\Ph(\rh_{\al_0},\al_0)\right)$, there is 
\beg{align}\label{eq-rhV}
w=-\al_0 \rh_{\al_0} \pi_{\al_0}{\bf V}_{2}w.
\end{align}
Taking into account $\pi_{\al_0}{\bf V_2} w\in L^2(\bar\mu)$, which is derived from $V_2\in L^2(\bar\mu\times\bar\mu)$ and $\rh_{\al_0}\in L^\infty$ (due to Remark \ref{rem:1}), there is $v\in L^2(\bar\mu)\subset L^2(\mu_{\al_0})$  so that $\mu_{\al_0}(v)=0$ and $w=\rh_{\al_0} v$. Note that $\mu_{\al_0}(v)=0$ yields $\pi_{\al_0}v=v$. Thus there is $w\in L^2(\bar\mu)$ satisfying \eqref{eq-rhV} if and only if there is $v\in L^2(\mu_{\al_0})$ satisfying 
$$v =- \al_0 \pi_{\al_0}{\bf V}_{2,\al_0}v=- \al_0 \pi_{\al_0}{\bf V}_{2,\al_0}\pi_{\al_0}v.$$
Hence, ${\rm Ker}\left(\nn\tld\Ph(\rh_{\al_0},\al_0)\right)=\{0\}$ on $L^2(\bar\mu)$ if and only if $I+\al_0 \pi_{\al_0}{\bf V}_{2,\al_0}\pi_{\al_0}$ is invertible on $L^2(\mu_{\al_0})$. Therefore, it follows from \cite[Theorem 15.1 and Theorem 15.3]{Deim} that there is a neighborhood of $\al_0$ such that $\al\mapsto \rh_\al $ is continuously differentiable in $L^2(\bar\mu)$. \\
Let $\de>0$ such that  $\rh_\al$ is continuously differentiable in $L^2(\bar\mu)$ for $\al\in J_{\al_0,\de}$. Then  $\|\rh_\al\|_{L^2(\bar\mu)}+\al \|\pp_\al\rh_\al\|_{L^2(\bar\mu)}$ is bounded of $\al$ on $J_{\al_0,\de}$. It follows from \eqref{ine-cT-K} that  for each $\al\in J_{\al_0,\de}$
\beg{equation}\label{rhal0}
|\rh_\al|=|\cT(\rh_\al,\al)|\leq \|e^{-\th(\al)V_0}\|_{L^1(\bar\mu)}^{-2} e^{ 2C_0(\th(\al ),\al ,\al \|\rh_{\al} \|_{L^2(\bar\mu)})}.
\end{equation}
Since $\sup_{\th_1\in R_{\th}}e^{-\th_1 V_0}\in L^1$, there is $R_{\th}\ni \th_1\mapsto \|e^{-\th_1V_0}\|_{L^1}$ is continuous. Thus 
\beg{equation}\label{supe-V0L1}
\sup_{\al\in J_{\al_0,\de}}\|e^{-\th(\al)V_0}\|_{L^1(\bar\mu)}^{-2}<+\infty.
\end{equation}
Combining this with \eqref{rhal0}, we obtain \eqref{sup-rh}. 

Since  $\rh_\al=\cT(\rh_\al,\al)$ and that $\|\pp_\al\rh_\al\|_{L^2(\bar\mu)}$ is bounded of $\al$ on $ J_{\al_0,\de}$, we have 
\beg{align*}
\pp_\al\log\rh_\al&=-\th'(\al)V_0 - V(\rh_\al\bar\mu)-\al {\bf V}_2\pp_{\al} \rh_{\al}-\pp_\al\log\bar\mu(e^{-\th(\al)V_0-\al V(\rh_\al\bar\mu)+\bar V}).
\end{align*}
The H\"older inequality yields that
\beg{equation}\label{V-alV}
\beg{split}
& |V(x,\rh_\al\bar\mu)|+\al | \left({\bf V}_2\pp_{\al} \rh_{\al}\right)(x)|\\
&\leq  |V_1(x)|+\|V_2(x,\cdot)\|_{L^2(\bar\mu)}(\|\rh_\al\|_{L^2(\bar\mu)}+\al\|\pp_\al\rh_\al\|_{L^2(\bar\mu)}).
\end{split}
\end{equation}
This, together with that $\|\rh_\al\|_{L^2(\bar\mu)}+\al \|\pp_\al\rh_\al\|_{L^2(\bar\mu)}$ is bounded of $\al$ on $ J_{\al_0,\de}$, $\th\in C^1(\hat{\al},\check{\al})$, $V_0\in L^2(\bar\mu)$ and \eqref{V-F3-ad}, implies by  the dominated convergence theorem that
\beg{align*}
&\pp_\al\log\bar\mu(e^{-\th(\al)V_0-\al V(\rh_\al\bar\mu)+\bar V})\\
&=\ff {-\bar\mu(e^{-\th(\al)V_0-\al V(\rh_\al\bar\mu)+\bar V}(\th'(\al)V_0+V(\rh_\al\bar\mu)+\al {\bf V}_2\pp_\al\rh_\al))} {\bar\mu(e^{-\th(\al)V_0-\al V(\rh_\al\bar\mu)+\bar V})}\\
&= - \mu_\al (\th'(\al)V_0+V(\rh_\al\bar\mu)+\al {\bf V}_2\pp_\al\rh_\al).
\end{align*}
Hence, 
\beg{align}\label{pplogrh}
\pp_\al\log\rh_\al&=-\th'(\al)V_0 - V(\rh_\al\bar\mu)-\al {\bf V}_2\pp_\al\rh_\al\nonumber\\
&\quad\, +\mu_\al (\th'(\al)V_0+V(\rh_\al\bar\mu)+\al {\bf V}_2\pp_\al\rh_\al),\\
|\pp_\al\log\rh_\al(x)| &\leq  |\th'(\al)V_0(x)|+|V_1(x)|+\mu_{\al}(|\th'(\al)||V_0|+|V_1|) \nonumber\\
&\quad + (\|\rh_\al\|_{L^2(\bar\mu)}+\al\|\pp_\al\rh_\al\|_{L^2(\bar\mu)})^2\left(\|V_2(x,\cdot)\|_{L^2(\bar\mu)}+\|V_2\|_{L^2(\bar\mu\times\bar\mu)}\right).\nonumber
\end{align}
This, together with \eqref{sup-rh}, \eqref{expVV} and that $\|\rh_\al\|_{L^2(\bar\mu)}+\al \|\pp_\al\rh_\al\|_{L^2(\bar\mu)}$ is bounded of $\al$ on $J_{\al_0,\de}$, implies  \eqref{sup-pp-rh}.  Since $\mu_{\al_0}(\pp_\al\log\rh_{\al_0})=0$, we have that 
\[{\bf V_2}\pp_\al\rh_{\al_0}={\bf V}_{2,\al_0}\pp_{\al}\log\rh_{\al_0}={\bf V}_{2,\al_0}\pi_{\al_0}\pp_{\al}\log\rh_{\al_0},\]
and we also derive from \eqref{pplogrh} that
\beg{align*}
\pp_\al\log\rh_{\al_0}&=-\th'(\al_0)(V_0-\mu_{\al_0}(V_0)) - (V(\mu_{\al_0})-\mu_{\al_0}(V(\mu_{\al_0})))\nonumber\\
&\quad\, -\al_0\left({\bf V}_{2,\al_0}\pp_\al\log\rh_{\al_0} - \mu_{\al_0}\left({\bf V}_{2,\al_0}\pp_\al\log\rh_{\al_0}\right)\right)\\
& = -\pi_{\al_0}\left(\th'(\al_0)V_0+V(\mu_{\al_0})\right)-\al_0\pi_{\al_0} {\bf V}_{2,\al_0}\pi_{\al_0}(\pp_\al\log\rh_{\al_0}). 
\end{align*}
This implies \eqref{pp-logrh0}. 
\\
If \eqref{expVV} holds, then $\cT$ is continuous  from  $L^2(\bar\mu)\times (\hat{\al},\check{\al})$ to $L^r(\bar\mu)$ for any $r\geq 1$. Taking into account  $\rh_\al=\cT(\rh_\al,\al)$, we find that $\al\mapsto\rh_\al$ is continuous from $J_{\al_0,\de}$ to $L^r(\bar\mu)$. By \eqref{pplogrh}, we have for any $\al_1,\al_2\in J_{\al_0,\de}$ that
\beg{align*}
&|\pp_\al\log\rh_{\al_1}(x)-\pp_\al\log\rh_{\al_2}(x)|\\
&\leq \|V_2(x,\cdot)\|_{L^2(\bar\mu)}\Big(\|\rh_{\al_1}-\rh_{\al_2}\|_{L^2(\bar\mu)}+(\al_0+\de)\|\pp_{\al}\rh_{\al_1}-\pp_\al\rh_{\al_2}\|_{L^2(\bar\mu)} \\
&\quad +|\al_1-\al_2| \sup_{\al\in J_{\al_0,\de}}\|\pp_\al\rh_{\al}\|_{L^2(\bar\mu)}\Big)+|\th'(\al_1)-\th'(\al_2)||V_0(x)|\\
&\quad +C \left(\|\rh_{\al_1}-\rh_{\al_2}\|_{L^2(\bar\mu)}+|\al_1-\al_2|+\|\pp_\al\rh_{\al_1}-\pp_\al\rh_{\al_2}\|_{L^2(\bar\mu)}+|\th'(\al_1)-\th'(\al_2)|\right).
\end{align*}
where $C$ is a positive constant depending on $\al_0$,$\de$,$\|V_0\|_{L^2(\bar\mu)}$,$\|V_1\|_{L^2(\bar\mu)}$, $\|V_2\|_{L^2(\bar\mu\times\bar\mu)}$, and $\sup_{\al\in J_{\al_0,\de}}(|\th'(\al)|+ \|\rh_\al\|_{L^2(\bar\mu)}+\|\pp_\al\rh_\al\|_{L^2(\bar\mu)})$. This, together with \eqref{expVV}, implies that $\al\mapsto \pp_\al\log\rh_\al$ is continuous from $J_{\al_0,\de}$ to $L^r(\bar\mu)$ for any $r\geq 1$. Consequently, $\pp_\al\rh_\al=\rh_\al\pp_\al\log\rh_\al$ is continuous of  $\al$ from  $J_{\al_0,\de}$ to $L^r(\bar\mu)$.

\end{proof}

\beg{proof}[Proof of Lemma \ref{lem-ps-T}]
We first prove that, for each $w\in L^2(\mu_{\al})$, $ \Ph(w,\al)\in L^2(\bar\mu)$, i.e. $\rh_\al^{-1}\cT((w+1)\rh_\al;\al)\in L^2(\bar\mu)$, and
$$\rh_\al^{-1}\cT((w+1)\rh_\al,\al)=\ff { \Ps(w,\al)} {\mu_\al(\Ps(w,\al))}.$$  
Since Remark \ref{rem:1} and  $\rh_\al$ is a fixed point of $\cT(\cdot,\al)$ without $K_1,K_2$, there is $\rh_{\al}\in L^\infty$. It is clear that  
\beg{align}\label{muK1}
\mu_\al(|w K_1|) \leq  \|w\|_{L^2(\mu_\al)}\|K_1\|_{L^2(\mu_\al)}\leq \|\rh_\al\|_{\infty}^{\ff 1 2}\|w\|_{L^2(\mu_\al)}\|K_1\|_{L^2(\bar\mu )}.
\end{align}
Due to \eqref{K2mu0}, we have that 
\beg{equation}\label{ex-rh-0}
\beg{split}
&\exp\left\{-(\th(\al)V_0-\bar V)(x)-\al\left(V_1(x)+\mu_\al\left( ( V_2+K_2)(x,\cdot) \right)\right)\right\}\\
& =\exp\left\{-(\th(\al)V_0-\bar V)(x)-\al  V_1(x)-\al \mu_\al( V_2 (x,\cdot))\right\}\\
& = \rh_\al(x) \bar\mu(e^{-\th(\al)V_0 -\al V(\rh_\al\bar\mu)+\bar V }).
\end{split}
\end{equation}
Combing this with the  following two iequalities 
\beg{align}\label{V2-rh-b}
\int_{\R^d} |V_2(x,y)w(y)|\mu_\al(\d y)&\leq \|\rh_\al\|_\infty^{\ff 1 2} \|V_2(x,\cdot)\|_{L^2(\bar\mu)}\|w\|_{L^2( \mu_\al)},\\
\int_{\R^d} |K_2(x,y) w(y) |\mu_\al(\d y) &\leq \|\rh_\al\|_\infty^{\ff 1 2}\|K_2(x,\cdot)\|_{L^2(\bar\mu)}\|w\|_{L^2( \mu_\al)},\label{K2-rh-b}
\end{align}
we find that 
\beg{equation}\label{ex-inrh}
\beg{split}
&\exp\left\{-(\th(\al)V_0-\bar V)(x)-\al (V_1(x)+\mu_\al\left( (w+1)  (V_2+K_2)(x,\cdot)  \right))\right\}\\
& = \rh_\al(x) \bar\mu(e^{-\th(\al)V_0 -\al V(\rh_\al\bar\mu)+\bar V })\Ps(x;w,\al)\\
&\leq   \rh_\al(x)\bar\mu(e^{-\th(\al)V_0 -\al V(\rh_\al\bar\mu)+\bar V }) e^{\mathscr{V}(x;\al,w)},
\end{split}
\end{equation}
where
\beg{align}\label{VV}
\mathscr{V}(x;\al,w):=\al \|\rh_\al\|_\infty^{\ff 1 2} \|w\|_{L^2( \mu_\al)} \left( \|V_2(x,\cdot)\|_{L^2(\bar\mu)}+\|K_2 (x,\cdot)\|_{L^2(\bar\mu)}\right).
\end{align}
Let 
\[Z(\rh\rh_\al,\al)=\int_{\R^d}\exp\{-\th(\al)V_0(x)-\al \int_{\R^d}V(x,y)\rh(y)\mu_{\al}(\d y)+\bar V(x)\}\bar\mu(\d x).\]
Then, \eqref{ex-inrh}, together with \eqref{expVV}, \eqref{V1-z-V}, \eqref{muK1} and \eqref{ine-mu-cT-K}, implies that there is $C>0$ which depends on ${\al,\|\rh_\al\|_{\infty},\|K_1\|_{L^2(\bar\mu)},\|w\|_{L^2(\mu_\al)}}$ such that
\beg{align*}
Z((w+1)\rh_\al,\al)&= e^{-\al\mu_\al(K_1(w+1))}\bar\mu\left(\rh_\al \Ps(w,\al)\right) \bar\mu(e^{-V_0 -\al V(\rh_\al\bar\mu)+\bar V })\\
&\leq C\|\rh_\al\|_\infty \bar\mu\left( e^{\mathscr{V}(\al,w)}\right)\bar\mu(e^{-\th(\al)V_0 -\al V(\rh_\al\bar\mu)+\bar V })\\
&<\infty,\\
\|\rh_\al^{-1}\cT((w+1)\rh_\al,\al)\|_{L^2(\bar\mu)}&\leq\ff {e^{-\al\mu_\al(K_1(w+1))} } {Z((w+1)\rh_\al,\al)}\|e^{\mathscr{V}(\al,w)}\|_{L^2(\bar\mu)}\bar\mu(e^{-\th(\al)V_0 -\al V(\rh_\al\bar\mu)+\bar V })\\
&<\infty.
\end{align*}
According to \eqref{ex-rh-0},  
\beg{equation}\label{cT-s}
\beg{split}
\rh_\al^{-1}(x)\cT(x;(w+1)\rh_\al,\al)&=\ff { \Ps(x;w,\al) e^{-\al \mu_\al(K_1(w+1))}} {\int_{\R^d}\rh_\al(x)\Ps(x;w,\al) e^{-\al\mu_\al(K_1(w+1))}\bar\mu(\d x)}\\
&=\ff { \Ps(x;w,\al)} {\mu_\al(\Ps(w,\al))}.
\end{split}
\end{equation} 
In particular, $\rh_\al^{-1}\cT(\rh_\al,\al)=1$ (or $\cT(\rh_\al,\al)=\rh_\al$) and $\Ph(0,\al)=0$.

Next, we prove the regularity of $\Ph(\cdot,\al)$. For each $w\in L^2( \mu_{\al})$, following from \eqref{V2-rh-b}, \eqref{K2-rh-b}, \eqref{VV} and the inequality
$$x^2 \leq 2 e^x,~x\geq 0,$$
we find that
\beg{equation}\label{ine-VVK}
\beg{split}
\left|\log\Ps(x;w,\al)\right|^2&\leq \mathscr{V}(x;\al,w)^2\\
&\leq 2\al^2\|\rh_\al\|_{\infty}\|w\|_{L^2( \mu_\al)}^2\exp\left(\ff {\mathscr{V}(x;\al,w)} {\al^2\|\rh_\al\|_{\infty}\|w\|_{L^2( \mu_\al)}^2} \right).
\end{split}
\end{equation}
Then,  for  $w_1,w\in L^2( \mu_\al)$  and  constant $M>0$, there is a constant $C>0$ depending on $M$, $\|w\|_{L^2( \mu_\al)},\|w_1\|_{L^2(\mu_\al)}$ and $\|\rh_\al\|_\infty$  such that
\beg{align*}
&\sup_{s\in [-M,M]}\left|\ff {\d} {\d s}\Ps(sw +w_1,\al)\right|^2 =\sup_{s\in [-M,M]}\left|\Ps(sw +w_1,\al)\log\Ps(w_1,\al)\right|^2\\
&\leq C\exp\left\{ C\left( \|V_2(x,\cdot)\|_{L^2(\bar\mu)}+\|K_2(x,\cdot)\|_{L^2(\bar\mu)})\right)\right\}. 
\end{align*}
This, together with \eqref{expVV}, \eqref{V1-z-V} and the dominated convergence theorem, implies that  the mapping $s\mapsto \Ps(sw+w_1,\al)$ is differentiable in $L^2(\bar\mu)$, which also implies the mapping $ s\mapsto \mu_\al(\Ps(sw+w_1,\al))$ is differentiable.  Since
$$\Ph(w_1,\al)=w_1+1-\rh_\al^{-1} \cT((w_1+1)\rh_\al,\al)=w_1+1- \ff {\Ps(w_1,\al)} {\mu_\al(\Ps(w_1,\al))},$$
we have proved that $\Ph(\cdot,\al)$ is  G\^ateaux differentiable from $L^2( \mu_\al)$ to $L^2(\bar\mu)$ and the G\^ateaus derivative is given by \eqref{nnPh0}.\\
For $ w,\tld w\in L^2(\mu_\al)$,  following from  \eqref{ine-VVK} and \eqref{funi-1},  we have that
\beg{align*}
\left|\Ps(\tld w,\al)-\Ps(w,\al)\right|&=\left|\Ps(\tld w-w,\al)-1\right|\Ps(w,\al)\\
&\leq \left|\log\Ps(\tld w-w,\al)\right|\exp\left(|\log\Ps(\tld w-w,\al)|\right)\Ps(w,\al)\\
&\leq \sq 2\al \|\rh_\al\|_{\infty}^{\ff 1 2}\|w-\tld w\|_{L^2(\mu_\al)}\exp\left( \mathscr{V}(\tld w-w,\al)+ \mathscr{V}(w,\al)\right).
\end{align*}
This implies that $\Ps(\cdot,\al)$ is continuous from $L^2(\mu_\al)$ to $L^2(\bar\mu)$ since \eqref{expVV}, \eqref{V1-z-V} and \eqref{VV}. Moreover, we also have 
\beg{align*}
&\left| \Ps(\tld w,\al)\log\Ps(w_2,\al)-\Ps(w,\al) \log\Ps(w_1,\al)\right|\\
&\leq\left|\left(\Ps(\tld w,\al)-\Ps(w,\al)\right)\log\Ps(w_2,\al)\right|+ \Ps(w,\al)\left|\log\Ps(w_2,\al)-\log\Ps(w_1,\al)\right|\\
&\leq\left|\left(\Ps(\tld w,\al)-\Ps(w,\al)\right)\log\Ps(w_2,\al)\right|+ \Ps(w,\al)\left|\log\Ps(w_2-w_1,\al)\right|\\
&\leq 2 \al^2 \|\rh_\al\|_\infty \|w-\tld w\|_{L^2(\mu_\al)}\|w_2\|_{L^2(\mu_\al)}\exp\left( \mathscr{V}(\tld w-w,\al)+ \mathscr{V}(w,\al)\right)\\
&\quad +\sq 2\al  \|\rh_\al\|_\infty^{\ff 1 2} \|w_1-w_2\|_{L^2(\mu_\al)}\exp\left(  \mathscr{V}(w,\al)\right),~w,w_1,w_2\in L^2(\mu_\al).
\end{align*}
Due to \eqref{expVV} and \eqref{V1-z-V}, we have proven that the following mapping  is continuous from $L^2(\mu_\al)\times L^2(\mu_\al)$ to $L^2(\bar\mu)$
\beg{align*}
 (w,w_1)\mapsto \Ps(w,\al) \log\Ps(w_1,\al).
\end{align*}
According to \eqref{nnPh0}, we have proven that  $\nn \Ph(\cdot,\al)\in C(L^2(\mu_\al)\times L^2(\mu_\al);L^2(\bar\mu))$, and this implies that $ \Ph(\cdot,\al)$ is Fr\'echet differentiable with  Fr\'echet derivative given by \eqref{nnPh0}.

It is clear that $\Ps(0,\al)=1$. Thus \eqref{nnPh0} yields that for every $w\in L^2(\bar\mu)$
$$\nn_{w} \Ph(0,\al) =w -\log\Ps(w ,\al)+\mu_\al(\log\Ps(w ,\al)),$$
which implies \eqref{nnPh(0)}.

Noticing $\rh_\al\in L^\infty$ since Remark \ref{rem:1}, we can derive directly from \eqref{nnPh(0)}, \eqref{expVV} and \eqref{V1-z-V} that $\nn \Ph(0,\al)$ is a  Fredholm operator on $L^2(\bar\mu)$ and $L^2(\mu_\al)$.

\end{proof}

\beg{proof}[Proof of Corollary \ref{cor:mu0mu}]

We choose $\de$ small such that 
$$2\th(\al)-\th(\al_0)\in R_\th,\al\in J_{\al_0,\de}.$$ 
We prove the regularity of $ \Ph(\cdot,\al)$ from $L^2(\mu_{\al_0})$ to $L^2(\bar\mu)$. For $\al\in J_{\al_0,\de}$,  by \eqref{V-F3-ad}, we prove as \eqref{ine-mu-cT-K} and \eqref{ine-cT-K} that
\beg{equation}\label{rh2rh}
\beg{split}
&\ff {\rh_\al^2} {\rh_{\al_0}}(x)\\
&\leq \ff {\bar\mu(e^{-\th(\al_0)V_0-\al_0 V(\rh_{\al_0}\bar\mu)+\bar V})} {\bar\mu(e^{-\th(\al)V_0-\al V(\rh_{\al}\bar\mu)+\bar V})^2}\exp\{-(2\th(\al)-\th(\al_0) )V_0(x)+\bar V(x)\}\\
&\quad \times \exp\{\|V_2(x,\cdot)\|_{L^2(\bar\mu)}\|2\al\rh_\al-\al_0\rh_{\al_0}\|_{L^2(\bar\mu)} -(2\al-\al_0)V_1(x)\}\\
&\leq   \|e^{-\th(\al)V_0}\|_{L^1}^{-4} \exp\left\{ C_0(2\th(\al)-\th(\al_0) ,2\al-\al_0,\|2\al\rh_\al-\al_0\rh_{\al_0}\|_{L^2(\bar\mu)})\right\}\\
&\quad \times \exp\left\{ 2C_0(\th(\al),\al,\al\|\rh_{\al}\|_{L^2(\bar\mu)})+C_0(\th(\al_0),\al_0,\al_0\|\rh_{\al_0}\|_{L^2(\bar\mu)})\right\}.
\end{split}
\end{equation}
Taking into account that $\|\rh_\al\|_{L^2(\bar\mu)}$ is bounded of $\al$ on $  J_{\al_0,\de}$ and  \eqref{supe-V0L1}, and setting
\beg{align}\label{Cal0de}
C_{\al_0,\de}=&\sup_{\al\in J_{\al_0,\de}}\Big\{\|e^{-\th(\al)V_0}\|_{L^1}^{-4} \exp\left\{ C_0(2\th(\al)-\th(\al_0) ,2\al-\al_0,\|2\al\rh_\al-\al_0\rh_{\al_0}\|_{L^2(\bar\mu)})\right\} \nonumber\\
&\quad  \times \exp\left\{ 2C_0(\th(\al),\al,\al\|\rh_{\al}\|_{L^2(\bar\mu)})+C_0(\th(\al_0),\al_0,\al_0\|\rh_{\al_0}\|_{L^2(\bar\mu)})\right\}\Big\},
\end{align}
then $C_{\al_0,\de}<+\infty$, and 
\beg{equation*} 
\beg{split}
\left|\mu_{\al}\left((V_2+K_2)(x,\cdot)w(\cdot)\right)\right|& \leq \|(V_2+K_2)(x,\cdot)\|_{L^2(\bar\mu)}\|w\rh_{\al}\|_{L^2(\bar\mu)}\\
& = \|(V_2+K_2)(x,\cdot)\|_{L^2(\bar\mu)}\left\|w\ff {\rh_{\al}} {\sq\rh_{\al_0}}\right\|_{L^2( \mu_{\al_0})}\\
& \leq \sq{C_{\al_0,\de}}\|(V_2+K_2)(x,\cdot)\|_{L^2(\bar\mu)}\|w\|_{L^2(\mu_{\al_0})},~\al\in J_{\al_0,\de}.
\end{split}
\end{equation*}
By using this inequality,   \eqref{ine-VVK} holds with $\mathscr{V}$ replaced by 
$$\tld{\mathscr{V}}(x;\al_0,w):=\sq{C_{\al_0,\de}}\|(V_2+K_2)(x,\cdot)\|_{L^2(\bar\mu)}\|w\|_{L^2(\mu_{\al_0})}.$$
Then, repeating the  proof of the second assertion of  Lemma \ref{lem-ps-T}, we can prove the assertion on the regularity of $ \Ph(\cdot,\al)$, and \eqref{nnPh0} and \eqref{nnPh(0)} hold. 

By using \eqref{rh2rh}, we find that
\beg{align*}
&\sup_{\al\in J_{\al_0,\de}}\int_{\R^d\times\R^d}(V_2+K_2)^2(x,y) \left(\ff {\rh_{\al}} {\rh_{\al_0}}\right)^2(y)\mu_{\al_0}(\d x)\mu_{\al_0}(\d y)\\
&\qquad\leq 2\left(\sup_{\al\in J_{\al_0,\de}}\left\|\ff {\rh_{\al}} {\sq{\rh_{\al_0}}}\right\|_\infty^2\right)\|\rh_{\al_0}\|_\infty \left(\|V_2\|_{L^2_xL^2_y}^2+\|K_2\|_{L^2_xL^2_y}^2\right)\\
&\qquad =2C_{\al_,\de}\|\rh_{\al_0}\|_\infty \left(\|V_2\|_{L^2_xL^2_y}^2+\|K_2\|_{L^2_xL^2_y}^2\right)<+\infty.
\end{align*}
Combining this with
$$\mu_{\al}\left((V_2+K_2)(x,\cdot)w(\cdot)\right)=\mu_{\al_0}\left((V_2+K_2)(x,\cdot)\left(\ff {\rh_{\al}} {\rh_{\al_0}}w\right)(\cdot)\right),$$
we can derive from \eqref{nnPh(0)} that $\nn\Ph(0,\al)$ is a Fredholm operator on $L^2(\mu_{\al_0})$.  

\end{proof}

\beg{proof}[Proof of Lemma \ref{con-nnPh}]

It follows from (A2), \eqref{sup-pp-rh}, \eqref{rh2rh} and the H\"older inequality that
\beg{equation}\label{ad-in-pp1}
\beg{split}
&\int_{\R^d\times\R^d}|w_2(x) V_2(x,y) w_1(y)|\sup_{\al\in J_{\al_0,\de}}\left(|\pp_\al \rh_\al|(y)\rh_\al(x) \right) \bar\mu (\d x)\bar\mu(\d y)\\
&=\int_{\R^d\times\R^d}|(w_2\sq{\rh_{\al_0}})(x) V_2(x,y) (w_1\sq{\rh_{\al_0}})(y)|\\
&\quad\times\sup_{\al\in J_{\al_0,\de}}\left(|\pp_\al\log \rh_\al|(y)\ff {\rh_{\al}} {\sq{\rh_{\al_0}}}(y)\ff {\rh_\al} {\sq{\rh_{\al_0}}}(x) \right) \bar\mu (\d x)\bar\mu(\d y)\\
&\leq  \left(\sup_{\al\in J_{\al_0,\de}}\left\|\ff {\rh_\al} {\sq{\rh_{\al_0}}}\right\|_\infty^2\right) \left\|\sup_{\al\in  J_{\al_0,\de}} |\pp_\al\log \rh_\al| \right\|_{L^{\ff {2\ga_1} {\ga_1-2}}(\bar\mu)}\\
&\quad \times\|V_2\|_{L^2_xL^{\ga_1}_y}\|w_1\sq{\rh_{\al_0}}\|_{L^2(\bar\mu)}\|w_2\sq{\rh_{\al_0}}\|_{L^2(\bar\mu)}\\
&\leq C_{\al_0,\de}\left\|\sup_{\al\in  J_{\al_0,\de}} |\pp_\al\log \rh_\al| \right\|_{L^{\ff {2\ga_1} {\ga_1-2}}(\bar\mu)}\|V_2\|_{L^2_xL^{\ga_1}_y}\|w_1 \|_{L^2(\mu_{\al_0})}\|w_2 \|_{L^2(\mu_{\al_0})}.
\end{split}
\end{equation}
Thus, the dominated convergence theorem implies that, in $L^2(\mu_{\al_0})$,
\beg{align*}
\pp_\al \int_{\R^d}V_2(x,y)w(y)\mu_\al(\d y)&=\int_{\R^d}V_2(x,y)w(y)\pp_\al\rh_\al(y)\bar\mu(\d y)\\
&=\int_{\R^d}V_2(x,y)\ff {\pp_\al\rh_\al} {\rh_{\al_0}}(y)w(y) \mu_{\al_0}(\d y),~w\in L^2( \mu_{\al_0}).
\end{align*}
Similarly, we have that 
\beg{equation}\label{ad-in-pp2}
\beg{split}
&\int_{\R^d\times\R^d}|w_2(x)V_2(x,y)w_1(y)|\sup_{\al\in J_{\al_0,\de}}\left(|\pp_\al \rh_\al|(x) \rh_\al(y) \right) \bar\mu (\d x)\bar\mu(\d y)\\
&\leq   C_{\al_0,\de} \left\|\sup_{\al\in  J_{\al_0,\de}} |\pp_\al \log\rh_\al| \right\|_{L^{\ff {2\ga_1} {\ga_1-2}}(\bar\mu)}\|V_2\|_{L^2_yL^{\ga_1}_x}\|w_1\|_{L^2( \mu_{\al_0})}\|w_2\|_{L^2( \mu_{\al_0})}.
\end{split}
\end{equation}
The dominated convergence theorem implies that for every $w\in L^2(\mu_{\al_0})$
\beg{align*}
&\pp_\al\int_{\R^d\times\R^d} V_2(x,y)w(y)\mu_\al(\d x)\mu_\al(\d y)\\
&=\int_{\R^d\times\R^d}V_2(x,y)\pp_\al \rh_\al(x) w(y) \rh_\al(y) \bar\mu (\d x)\bar\mu (\d y)\\
&\quad +\int_{\R^d\times\R^d}V_2(x,y)(\pp_\al \rh_\al)(y) w(y) \rh_\al(x)\bar\mu (\d x)\bar\mu (\d y)\\
&=\int_{\R^d} \pp_\al\rh_\al(x) \left(\int_{\R^d}V_2(x,y)\rh_\al(y)w(y) \bar\mu(\d y)\right) \bar\mu(\d x)\\
&\quad +\int_{\R^d}\rh_\al(x)\left(\int_{\R^d}V_2(x,y)\pp_\al \rh_\al(y) w(y)\bar\mu(\d y)\right)\bar\mu (\d x).
\end{align*}
Hence, 
\beg{align*}
\pp_\al \left(\pi_{\al}  {\bf V}_{2,\al}w\right)&= \Big({\bf V}_{2}\cM_{\pp_\al \rh_\al} -(\1\otimes\pp_\al\rh_\al)  {\bf V}_{2}\cM_{\rh_\al}\\
&\qquad -(\1\otimes\rh_\al )  {\bf V}_{2}\cM_{\pp_\al\rh_\al}\Big)w,~w\in L^2(\mu_{\al_0}).
\end{align*}

Next, we discuss the continuity of $\al$ for $\pp_\al \left(\pi_{\al}  {\bf V}_{2,\al}\right)$. For ${\bf V}_{2}\cM_{\pp_\al \rh_\al}$, the H\"older inequality implies that
\beg{align*}
&\int_{\R^d\times\R^d}V_2(x,y)^2\left(\ff {\pp_\al\rh_{\al_1} -\pp_\al\rh_{\al_2}} {\rh_{\al_0}}\right)^2(y)\mu_{\al_0}(\d x)\mu_{\al_0}(\d y)\\
&\quad  \leq \|\rh_{\al_0}\|_\infty\int_{\R^d }\|V_2(\cdot,y)\|_{L^2(\bar\mu)}^2\left(\ff {\pp_\al\rh_{\al_1} -\pp_\al\rh_{\al_2}} {\rh_{\al_0}}\right)^2(y)\mu_{\al_0}(\d y)\\
&\quad\leq \|\rh_{\al_0}\|_\infty\|V_2\|_{L^2_xL^{\ga_1}_y}^2\left\|\ff {\pp_\al\rh_{\al_1} -\pp_\al\rh_{\al_2}} {\sq{\rh_{\al_0}}}\right\|_{L^{\ff {2\ga_1} {\ga_1-2}}}^2\\
&\quad\leq \|\rh_{\al_0}\|_\infty\|V_2\|_{L^2_xL^{\ga_1}_y}\left\| \pp_\al \log\rh_{\al_1}-\pp_\al \log\rh_{\al_2}\right\|_{L^{\ff {2\ga_1} {\ga_1-2}}}^2\left(\sup_{\al\in  J_{\al_0,\de}}\left\|\ff {\rh_\al^2} { \rh_{\al_0}}\right\|_\infty \right)\\
&\qquad +\|\rh_{\al_0}\|_\infty\|V_2\|_{L^2_xL^{\ga_1}_y}\left\|\ff {\rh_{\al_1}- \rh_{\al_2}} {\sq{\rh_{\al_0}}}\right\|_{L^{\ff {4\ga_1} {\ga_1-2}}}^2\left(\sup_{\al\in J_{\al_0,\de}}\left\| \pp_\al \log\rh_{\al }\right\|_{L^{\ff {4\ga_1} {\ga_1-2}}}^2\right).
\end{align*}
Due to \eqref{rh2rh} and Lemma \ref{con-rhal}, we see that ${\bf V}_{2}\cM_{\pp_\al \rh_\al}\in C(J_{\al_0,\de};\cL_{HS}(L^2(\mu_{\al_0})))$. We can prove similarly that ${\bf V}_{2}\cM_{\rh_\al}\in C(J_{\al_0,\de};\cL_{HS}(L^2(\mu_{\al_0})))$.  It follows from Lemma \ref{con-rhal}, \eqref{rh2rh} and 
\beg{align*}
\|\1\otimes \pp_\al\rh_{\al_1}-\1\otimes \pp_\al\rh_{\al_2}\|_{\sL_{HS}(L^2(\mu_{\al_0}))}&=\left\|\ff { \pp_\al\rh_{\al_1}-\pp_\al\rh_{\al_2}} {\sq{\rh_{\al_0}}}\right\|_{L^2(\bar\mu)}\\
& \leq C_{\al_0,\de}\left\| \pp_\al \log\rh_{\al_1}-\pp_\al \log\rh_{\al_2}\right\|_{L^2(\bar\mu)}\\
&\quad\, +\left(\sup_{\al\in J_{\al_0,de}}\|\pp_\al\log\rh_{\al}\|_{L^4}\right)\left\|\ff {\rh_{\al_1}-\rh_{\al_2}} {\sq{\rh_{\al_0}}}\right\|_{L^4},\\
\|\1\otimes  \rh_{\al_1}-\1\otimes  \rh_{\al_2}\|_{\sL_{HS}(L^2(\mu_{\al_0}))}&=\left\|\ff {  \rh_{\al_1}- \rh_{\al_2}} {\sq{\rh_{\al_0}}}\right\|_{L^2(\bar\mu)},
\end{align*}
that $\1\otimes \pp_\al\rh_{\al}$ and $\1\otimes  \rh_{\al}$ are also continuous in $\al$ from $J_{\al_0,\de}$ to $\cL_{HS}(L^2(\mu_{\al_0}))$.
Hence, $\pp_\al (\pi_{\al}{\bf V}_{2,\al})\in C( J_{\al_0,\de};\cL_{HS}(L^2(\mu_{\al_0})))$.

We can similarly prove that
\beg{align*}
\pp_\al (\pi_{\al}{\bf K}_{2,\al})= {\bf K}_{2}\cM_{\pp_\al\rh_\al}-(\1\otimes\pp_\al\rh_\al) {\bf K}_2\cM_{\rh_\al}-(\1\otimes\rh_\al)  {\bf K}_2\cM_{\pp_\al\rh_\al},
\end{align*}
and $\pp_\al( \pi_{\al}{\bf K}_{2,\al})\in C( J_{\al_0,\de};\cL_{HS}(L^2(\mu_{\al_0})))$.
Noticing that $I-\1\otimes\rh_\al=\pi_{\al}$, we arrive at
\beg{equation*}
\beg{split}
\pp_\al\nn\Ph(0,\al) & = \pi_{\al}({\bf V}_{2,\al}+{\bf K}_{2,\al})+\al \pp_\al( \pi_{\al}{\bf V}_{2,\al}+ \pi_{\al} {\bf K}_{2,\al})\\
& =(I-\1\otimes \rh_\al)({\bf V}_{2}+{\bf K}_{2})(\cM_{\rh_\al}+\al\cM_{\pp_\al\rh_\al})\\
&\quad -\al (\1\otimes\pp_\al\rh_\al) ({\bf V}_{2}+{\bf K}_{2})\cM_{\rh_\al}\\
& =\pi_{\al}({\bf V}_{2}+{\bf K}_{2})\cM_{\rh_\al+\al\pp_\al\rh_\al}-\al(\1\otimes\pp_\al\rh_\al) ({\bf V}_{2}+{\bf K}_{2})\cM_{\rh_\al}.
\end{split}
\end{equation*}

\end{proof}

\subsection{Proof of Theorem \ref{thm-bif}}

By the assumption of Theorem \ref{thm-bif} and Corollary \ref{cor:mu0mu}, $\nn\Ph(0,\al_0)$ is a Fredholm operator and $0$ is an isolate eigenvalue of $\nn\Ph(0,\al_0)$. Let $\Ga$ be a closed simple curve enclosing $0$ with diameter less than $1$ but no other eigenvalue of $\nn \Ph(0,\al_0)$ on $L^2_{\mathbb{C}}(\mu_{\al_0})$. Let $Q(\al)$ be the eigenprojection on $L^2_{\mathbb{C}}(\mu_{\al_0})$ given by $\Ga$ and $\nn \Ph(0,\al)$:
\beg{equation}\label{Prj}
Q(\al)=-\ff 1 {2\pi {\bf i}}\int_{\Ga}( \nn \Ph(0,\al) -\et)^{-1}\d\et.
\end{equation}
Due to Lemma \ref{con-nnPh}, \cite[Theorem IV. 2.23, Theorem 3.16, Section IV. 5]{Kato}, we have that
$$\lim_{|\al-\al_0|\ra 0^+}\|Q(\al)-Q(\al_0)\|_{\sL(L^2_{\mathbb{C}}(\mu_{\al_0}))}=0,$$
and there is $\de_1>0$ such that for every $\al\in J_{\al_0,\de_1}$ 
$$\dim\left( Q(\al)L^2_{\mathbb{C}}(\mu_{\al_0})\right)=\dim\left( Q(\al_0)L^2_{\mathbb{C}}(\mu_{\al_0})\right).$$ 
For $\al\in J_{\al_0,\de_1}-\{\al_0\}$, we call the spectrum of $\nn \Ph(0,\al)$ that is enclosed in the curve $\Gamma$ the $0$-group of $\nn \Ph(0,\al)$. 
\beg{defn}\label{def-odd-cro}
Let $\la_1,\cdots,\la_k$ be all the negative eigenvalues in the $0$-group of $\nn \Ph(0,\al)$ with algebraic multiplicities $m_1,\cdots,m_k$, respectively. Denote 
$$\si_{<}(\al)=(-1)^{\sum_{i=1}^k m_i},$$  
and set $\sum_{i=1}^k m_i=0$ if $k=0$. If $\nn \Ph(0,\al)$ is an isomorphism on $L^2_{\mathbb{C}}(\mu_{\al_0})$ for $\al\in   J_{\al_0,\de_1}-\{\al_0\}$ and $\si_{<}(\al)$ changes at $\al=\al_0$, then we say $\nn \Ph(0,\al)$ has an odd crossing number at $\al=\al_0$.
\end{defn}  

Due to the Krasnosel’skii Bifurcation Theorem (\cite[Theorem II.3.2]{Kie}), if $\nn\Ph(0,\al)$ has an odd crossing number at $\al_0$, then $\al_0$ is a bifurcation point of $ \Ph=0$.  To give a criteria for $\nn\Ph(0,\al)$ has an odd crossing number at $\al_0$, we use determinants for Fredholm operators, see \cite{Simon}. We denote by ${\rm det}(I+A)$ the Fredholm determinant of a trace class operator $A$ on $L^2_{\mathbb{C}}(\mu_{\al_0})$ and ${\rm det_2}(I+A)$ the regularized determinant for $A$ in the Hilbert-Schmidt class on $L^2_{\mathbb{C}}(\mu_{\al_0})$. According to the proof of Corollary \ref{cor:mu0mu}, $\al\pi_{\al}({\bf V}_{2}+{\bf K}_2)\cM_{\rh_\al}$ is a Hilbert-Schmidt operator on $L^2(\mu_{\al_0})$. Then we have the following lemma.
\beg{lem}\label{Ph-det2}
Suppose assumptions of Lemma \ref{con-nnPh} hold. Then $\nn \Ph(0,\al)$  has an odd crossing number at $\al=\al_0$ if and only if ${\rm det_2}(\nn \Ph(0,\al_0))=0$ and ${\rm det_2}(\nn\Ph(0,\al))$ changes sign at $\al=\al_0$.
\end{lem}
\beg{proof}
It follows from \cite[DEFINITION, THEOREM 9.2]{Simon} that
\beg{equation}\label{det2pi}
\beg{split}
{\rm det_2}(\nn \Ph(0,\al))&={\rm det}\left(\left(I+\al\pi_{\al}({\bf V}_{2}+{\bf K}_2)\cM_{\rh_\al}\right)e^{-\al\pi_{\al}({\bf V}_{2}+{\bf K}_2)\cM_{\rh_\al}}\right)\\
&=\prod\limits_{i}\left((1+ \ka_i(\al))e^{-\ka_i(\al)}\right),
\end{split}
\end{equation}
where $\{\ka_i(\al)\}$ are all the eigenvalues of $\al\pi_{\al}({\bf V}_{2}+{\bf K}_2)\cM_{\rh_\al}$ and the convergence in \eqref{det2pi} is absolute. Since $\al\pi_{\al}({\bf V}_{2}+{\bf K}_2)\cM_{\rh_\al}$ is a real Hilbert-Schmidt operator on $L^2_{\mathbb{C}}(\mu_{\al_0})$, where the ``real" operator means the operator that maps the real function in $L^2_{\mathbb{C}}(\mu_{\al_0})$ to a real function in $L^2_{\mathbb{C}}(\mu_{\al_0})$.  Then, for $\ka_i(\al)$ which is an eigenvalue of $\al\pi_{\al}({\bf V}_{2}+{\bf K}_2)\cM_{\rh_\al}$, the conjugate $\bar\ka_i(\al)$ is also an eigenvalue of $\al\pi_{\al}({\bf V}_{2}+{\bf K}_2)\cM_{\rh_\al}$ with the same algebraic multiplicity. Denote by ${\bf Im}(\ka_i(\al))$ the imaginary part of $\ka_i(\al)$, and  ${\bf Re}(\ka_i(\al))$ the real part of $\ka_i(\al)$. Let $D_{\Gamma}$ be the domain enclosed by the curve $\Gamma$. Then
\beg{align*}
{\rm det_2}(\nn  \Ph(0,\al)) &= \prod\limits_{ \ka_i(\al)\in\R}\left((1+\ka_i(\al)) e^{-\ka_i(\al)}\right) \\
&\qquad \times\prod\limits_{{\bf Im}(\ka_i(\al))>0}\left( |1+\ka_i(\al)|^2 e^{-(\ka_i(\al)+\bar\ka_i(\al))}\right)\\
&= \prod\limits_{\ka_i(\al)\in\R \atop 1+\ka_i(\al)\in D_{\Gamma}}\times\prod\limits_{\ka_i(\al)\in\R\atop 1+\ka_i(\al)\not\in D_{\Gamma}}\left((1+\ka_i(\al)) e^{-\ka_i(\al)}\right)\\
&\qquad \times\prod\limits_{{\bf Im}(\ka_i(\al))> 0}\left(|1+\ka_i(\al)|^2e^{-2{\bf Re}(\ka_i(\al)) }\right).
\end{align*}
Note that $\{1+\ka_i(\al)\}$ are eigenvalues of $\nn \Ph(0,\al)$ since \eqref{nnPh(0)} and the spectral mapping theorem, and that $\nn \Ph(0,\cdot)\in C(J_{\al_0,\de};\sL(L^2_{\mathbb{C}}(\mu_{\al_0})))$ due to Lemma \ref{con-nnPh}. We derive from the upper semicontinuity of the spectrum (see \cite[Remark IV.3.3]{Kato}) that, at $\al=\al_0$, the $\si_{<}(\al)$ for the $0$-group of $\nn \Ph(0,\al)$ changes if and only if the sign of ${\rm det_2}(\nn \Ph(0,\al))$ changes.   According to \eqref{det2pi}, \cite[THEOREM 9.2 (e)]{Simon} and  \eqref{nnPh(0)}, we have that ${\rm det_2}(\nn \Ph(0,\al_0))=0$ if and only if  $0$ is an eigenvalue of  finite  algebraic multiplicity of $\nn \Ph(0,\al_0)$. 

\end{proof}

Let $\la\neq 1$ be an eigenvalue of  $\nn \Ph(0,\al)$ as a Fredholm operator on $L_{\mathbb{C}}^2( \mu_{\al_0})$. According to \cite[Theorem 21.2.6 and Theorem 25.2.2']{Lax}, there is an integer $k_0$ such that
\beg{align*}
{\rm dim}~{\rm Ker}\left((\la I-\nn \Ph(0,\al))^{k_0}\right)&=\max_{k\in\N} {\rm dim}~{\rm Ker}\left((\la I-\nn \Ph(0,\al))^{k}\right)<\infty,\\
{\rm Ker}\left((\la I-\nn \Ph(0,\al))^{k_0}\right) & = {\rm Ker}\left((\la I-\nn \Ph(0,\al))^{k}\right),~k>k_0.
\end{align*}
The dimension of ${\rm Ker}\left((\la I-\nn \Ph(0,\al))^{k_0}\right)$ is the algebraic multiplicity of $\la $, and  functions in ${\rm Ker}\left((\la I-\nn \Ph(0,\al))^{k_0}\right)$ are called generalized eigenfunctions.  
The following lemma indicates that all the eigenvalues except $1$ and the associated generalized eigenfunctions  of  $\nn \Ph(0,\al)$ as an operator on $L_{\mathbb{C}}^2( \mu_{\al_0})$ are the same as that of $\nn \Ph(0,\al)$  on $L_{\mathbb{C}}^2( \bar\mu )$. We denote by ${\rm Ran}_{\mathbb{C}}(\pi_{\al}({\bf V}_{2,\al}+{\bf K}_{2,\al}))$  the range of  the complexified operator of $\pi_{\al}\left({\bf V}_{2,\al}+{\bf K}_{2,\al}\right)$.

\beg{lem}\label{Ker-nnPh}
The assumptions of Lemma \ref{lem-ps-T} hold. Let $\la\neq 1$ be an eigenvalue of  $\nn \Ph(0,\al)$ on $L^2(\mu_{\al_0})$ and $k_0$ be defined as above. Then 
\beg{align}\label{ker-sub}
{\rm Ker}\left((\la I-\nn \Ph(0,\al))^{k_0}\right)\subset {\rm Ran}_{\mathbb{C}}(\pi_{\al}({\bf V}_{2,\al}+{\bf K}_{2,\al})) \subset L_{\mathbb{C}}^2(\bar\mu).
\end{align}
\end{lem}
\beg{proof}

For any $w\in {\rm Ker}\left((\la I-\nn \Ph(0,\al))^{k_0}\right)$ and  $0\leq k\leq k_0$, we denote $w^{[k]}=(\la I-\nn \Ph(0,\al))^{k}w$.   Then $w^{[k_0]}=0$, $w^{[k]}=(\la I-\nn \Ph(0,\al))w^{[k-1]}$ and $w^{[0]}=w$. We first derive from $w^{[k_0]}=0$ that
\beg{align*}
0 =(\la I-\nn \Ph(0,\al)) w^{[k_0-1]} = (\la-1) w^{[k_0-1]} +\al\pi_{\al}({\bf V}_{2,\al}+{\bf K}_{2,\al})(w^{[k_0-1]}).
\end{align*}
It follows from  $\la\neq 1$ that  
\beg{align*}
w^{[k_0-1]}=-\ff {\al} {\la-1}\pi_{\al}({\bf V}_{2,\al}+{\bf K}_{2,\al})(w^{[k_0-1]}),
\end{align*}
which  implies that $w^{[k_0-1]}\in {\rm Ran}_{\mathbb{C}}(\pi_{\al}({\bf V}_{2,\al}+{\bf K}_{2,\al})) \subset L_{\mathbb{C}}^2(\bar\mu)$. If, for  $1\leq k\leq k_0$, there is $w^{[k]}\in {\rm Ran}_{\mathbb{C}}(\pi_{\al}({\bf V}_{2,\al}+{\bf K}_{2,\al}))$, then  we can derive from $w^{[k]}=(\la I-\nn \Ph(0,\al))w^{[k-1]}$ that 
\beg{align*}
w^{[k-1]}= \ff {w^{[k]}} {\la-1} -\ff {\al} {\la-1}(\pi_{\al}({\bf V}_{2,\al}+{\bf K}_{2,\al}))(w^{[k-1]}), 
\end{align*}
which implies that $w^{[k-1]}\in {\rm Ran}_{\mathbb{C}}(\pi_{\al}({\bf V}_{2,\al}+{\bf K}_{2,\al}))$. By iteration, we have that 
$$w^{[k]}\in {\rm Ran}_{\mathbb{C}}(\pi_{\al}({\bf V}_{2,\al}+{\bf K}_{2,\al}))\subset L_{\mathbb{C}}^2(\bar\mu),~0\leq k\leq k_0.$$
Particularly, $w=w^{[0]}\in {\rm Ran}_{\mathbb{C}}(\pi_{\al}({\bf V}_{2,\al}+{\bf K}_{2,\al}))$. Hence,  \eqref{ker-sub} holds.

\end{proof}

\beg{rem}
From this lemma, we have that $\overline{{\rm Ran}(\pi_{\al}({\bf V}_{2,\al}+{\bf K}_{2,\al}))}$ is an invariant subspace of $\nn \Ph(0,\al)$. Then we can use 
$${\rm det_2}\left(I+\al \pi_{\al}({\bf V}_{2,\al}+{\bf K}_{2,\al})\Big|_{\overline{{\rm Ran}(\pi_{\al}({\bf V}_{2,\al}+{\bf K}_{2,\al}))}}\right)$$
to character whether $\nn \Ph(0,\al)$ has an odd crossing  number  at $\al=\al_0$, when ${\bf V}_2$ and ${\bf K}_2$ have finite rank.
\end{rem}

\beg{lem}\label{det22}
Suppose assumptions of Corollary \ref{cor:mu0mu} hold. Then
\beg{align*}
{\rm det_2}(\nn \Ph(0,\al))={\rm det_2}\left( I+\al\pi_{\al}({\bf V}_{2}+{\bf K}_{2})\cM_{\rh_\al}\pi_{\al}\right).
\end{align*}
\end{lem}
\beg{proof} 
According to Corollary \ref{cor:mu0mu},
\beg{align*}
\nn \Ph(0,\al)&=I+\al\pi_{\al}({\bf V}_{2,\al}+{\bf K}_{2,\al})\\
&=I+\al\pi_{\al}({\bf V}_{2,\al}+{\bf K}_{2,\al})\pi_{\al}+\al\pi_{\al}({\bf V}_{2,\al}+{\bf K}_{2,\al})(I-\pi_{\al})\\
&=\left( I+\al\pi_{\al}({\bf V}_{2,\al}+{\bf K}_{2,\al})(I-\pi_{\al})\right)\left(I+\al\pi_{\al}({\bf V}_{2,\al}+{\bf K}_{2,\al})\pi_{\al}\right)
\end{align*}
By \cite[(2.40)]{GLZ08}, we find that
\beg{equation}\label{det2ph-pi}
\beg{split}
&{\rm det}_2(\nn \Ph(0,\al))\\
&={\rm det}_2\left( I+\al\pi_{\al}({\bf V}_{2,\al}+{\bf K}_{2,\al})(I-\pi_{\al})\right){\rm det_2}\left(I+\al\pi_{\al}({\bf V}_{2,\al}+{\bf K}_{2,\al})\pi_{\al}\right).
\end{split}
\end{equation}
For every $f\in L^2(\mu_{\al_0})$, 
$$\pi_{\al}({\bf V}_{2,\al}+{\bf K}_{2,\al})(I-\pi_{\al})f=\mu_\al(f)\pi_{\al}({\bf V}_{2,\al}+{\bf K}_{2,\al})\1.$$
We find that $\al\pi_{\al}({\bf V}_{2,\al}+{\bf K}_{2,\al})(I-\pi_{\al})$ is a finite rank operator,  and $0$ is the only eigenvalue since $\mu_{\al}(\pi_{\al}({\bf V}_{2,\al}+{\bf K}_{2,\al})\1)=0$.  Thus,  the trace $\trac(\pi_{\al}({\bf V}_{2,\al}+{\bf K}_{2,\al})(I-\pi_{\al}))=0$ 
and
\beg{align*}
{\rm det_2}\left( I+\al\pi_{\al}({\bf V}_{2,\al}+{\bf K}_{2,\al})(I-\pi_{\al})\right)&=\det\left(\left( I+\al\pi_{\al}({\bf V}_{2,\al}+{\bf K}_{2,\al})(I-\pi_{\al})\right)\right)=1.
\end{align*}
Substituting this into \eqref{det2ph-pi} and taking to account that 
$$({\bf V}_{2}+{\bf K}_{2})\cM_{\rh_\al}={\bf V}_{2,\al}+{\bf K}_{2,\al},$$
the corollary is proved.

\end{proof}

\beg{rem}
Combining Lemma \ref{Ph-det2} with Lemma \ref{det22}, we have that ${\rm det_2}(\nn \Ph(0,\al_0))=0$ if and only if $1$ is an eigenvalue of $-\al_0\pi_{\al_0}({\bf V}_{2,\al_0}+{\bf K}_{2,\al_0})\pi_{\al_0}$.
\end{rem}
Combining this remark and the following lemma, the proof of main assertions of Theorem \ref{thm-bif} is finished. 
\beg{lem}
Suppose assumptions of Theorem \ref{thm-bif} hold. Then ${\rm det}_{2}(I+\al \pi_{\al}({\bf V}_{2,\al}+{\bf K}_{2,\al})\pi_{\al })$ changes sign at $\al=\al_0$.
\end{lem}
\beg{proof}
Due to \eqref{rh2rh}, $\pi_{\al}$ is a bounded operator on $L^2(\mu_{\al_0})$ , and
\beg{align*}
\pp_\al\pi_{\al}f&=-\bar\mu(f\pp_{\al}\rh_{\al})=-\mu_{\al}(f\pp_{\al}\log\rh_{\al}),\\
\pi_{\al} f&=\pi_{\al_0} f-(\al-\al_0)\mu_{\al_0}(f\pp_{\al}\log\rh_{\al_0})\\
&\quad -\int_{\al_0}^{\al}\left(\bar\mu(f\pp_\al\rh_s)-\bar\mu(f\pp_\al \rh_{\al_0})\right)d s.
\end{align*}
The H\"older inequality implies that
\beg{align*}
\left|\bar\mu(f\pp_\al\rh_s)-\bar\mu(f\pp_\al \rh_{\al_0})\right|&\leq \|f\sq{\rh_{\al_0}}\|_{L^2(\bar\mu)}\left\|\ff {\pp_{\al}\rh_s-\pp_{\al}\rh_{\al_0}} {\sq{\rh_{\al_0}}}\right\|_{L^2(\bar\mu)}\\
&=\|f \|_{L^2( \mu_{\al_0})}\left\|\ff {\rh_{s}} {\sq{\rh_{\al_0}}}\pp_{\al}\log\rh_s-\ff { \rh_{\al_0}} {\sq{\rh_{\al_0}}} \pp_{\al}\log\rh_{\al_0} \right\|_{L^2(\bar\mu)}.
\end{align*}
According to \eqref{rh2rh} and Lemma \ref{con-rhal}, we have that
\beg{align*}
\sup_{s\in J_{\al_0,\de}}\left\|\ff {\rh_{s}} {\sq{\rh_{\al_0}}}\pp_{\al}\log\rh_s\right\|_{L^r(\bar\mu)}<+\infty,~r>2,\\
\ff {\rh_{s}} {\sq{\rh_{\al_0}}}\pp_{\al}\log\rh_s\ra\ff { \rh_{\al_0}} {\sq{\rh_{\al_0}}} \pp_{\al}\log\rh_{\al_0},~\mbox{in}~\bar\mu.
\end{align*}
Thus, the dominated convergence theorem implies that
\beg{align*}
\lim_{\al\ra \al_0}\sup_{\|f\|_{L^2(\al_0)}\leq 1} \left(\ff 1 {\al-\al_0} \left|\int_{\al_0}^{\al}\left(\bar\mu(f\pp_\al\rh_s)-\bar\mu(f\pp_\al \rh_{\al_0})\right)d s\right|\right)=0.
\end{align*}
Hence, on $\sL(L^2(\mu_{\al_0}))$
\[\pi_\al=\pi_{\al_0}-(\al-\al_0)\1\otimes_{\al_0}\pp_{\al}\log\rh_{\al_0}+o(|\al-\al_0|),\]
where $\otimes_{\al_0}$ is the tensor product on $L^2(\mu_{\al_0})$. Combining this with Lemma  \ref{con-nnPh}, we find that
\beg{equation}\label{VK-expa1}
\beg{split}
&\al\pi_{\al}({\bf V}_{2,\al}+{\bf K}_{2,\al}) \pi_{\al}=(\nn \Ph(0,\al)-I)\pi_\al\\
&\qquad= (\nn \Ph(0,\al_0)-I+(\al-\al_0)(\pp_\al\nn \Ph(0,\al_0))+o(|\al-\al_0|))\\
&\qquad\quad \times (\pi_{\al_0}-(\al-\al_0)\1\otimes_{\al_0}\pp_{\al}\log\rh_{\al_0}+o(|\al-\al_0|))\\
&\qquad= \al_0\pi_{\al_0}({\bf V}_{2,\al_0}+{\bf K}_{2,\al_0})\pi_{\al_0}+(\al-\al_0)\pp_\al\nn \Ph(0,\al_0) \pi_{\al_0}\\
&\qquad\quad -(\al-\al_0) \al_0\pi_{\al_0}({\bf V}_{2,\al_0}+{\bf K}_{2,\al_0}) \left(\1\otimes_{\al_0}\pp_{\al}\log\rh_{\al_0}\right) +o(|\al-\al_0|)
\end{split}
\end{equation}
Denote
\beg{align*}
A_0&=-\al_0\pi_{\al_0}({\bf V}_{2,\al_0}+{\bf K}_{2,\al_0})\pi_{\al_0},\\
A_1&=\pp_\al\nn\Ph(0,\al_0) \pi_{\al_0}-\al_0\pi_{\al_0}({\bf V}_{2,\al_0}+{\bf K}_{2,\al_0}) \left(\1\otimes_{\al_0}\pp_{\al}\log\rh_{\al_0}\right).
\end{align*}
Let $\cH_1=(I-P(\al_0))L^2_{\mathbb{C}}(\mu_{\al_0})$, and let $k_0={\rm dim}(\cH_0) $ be the algebraic multiplicity of the eigenvalue $1$.  According to \cite[Theorem 2.7]{GLZ08}, we derive from \eqref{VK-expa1} that
\beg{align*}
&{\rm det}_{2}(I+\al \pi_{\al}({\bf V}_{2,\al}+{\bf K}_{2,\al})\pi_{\al })\\
&=\left[ {\rm det}_{2,\cH_1}(I_{\cH_1}- (I-P(\al_0))A_0(I-P(\al_0)))+o(1)\right]e^{k_0}(-1)^{k_0}\\
&\qquad \times {\rm det}_{2,\cH_0}\left(P(\al_0)(A_0-I)P(\al_0)-P(\al_0)A_1P(\al_0)(\al-\al_0)+o(\al-\al_0)\right)\\
&=\left[ {\rm det}_{2,\cH_1}(I_{\cH_1}- (I-P(\al_0))A_0(I-P(\al_0)))+o(1)\right]e^{k_0}\\
&\qquad \times {\rm det}_{2,\cH_0}\left(P(\al_0)\nn \Ph(0,\al_0)P(\al_0)+P(\al_0)A_1P(\al_0)(\al-\al_0)+o(\al-\al_0)\right),
\end{align*}
where ${\rm det}_{2,\cH_0}$ and ${\rm det}_{2,\cH_1}$ are regularized determinants on $\cH_0$ and $\cH_1$ respectively. Noticing that $I_{\cH_1}- (I-P(\al_0))A_0(I-P(\al_0))$ is invertible on $\cH_1$ and $e^{k_0}$ is a constant, one can see that ${\rm det}_{2}(I+\al \pi_{\al}({\bf V}_{2,\al}+{\bf K}_{2,\al})\pi_{\al })$ changes sign if and only if 
\[{\rm det}_{2,\cH_0}\left(P(\al_0)\nn \Ph(0,\al_0)P(\al_0)+P(\al_0)A_1P(\al_0)(\al-\al_0)+o(\al-\al_0)\right)\]
changes sign.

Since $\1$ is the eigenvector of $\al_0\pi_{\al_0}({\bf V}_{2,\al_0}+{\bf K}_{2,\al_0})\pi_{\al_0}$ associated to the eigenvalue $0$, we have that $P(\al_0)\1=0$. Then
\[P(\al_0)\left(\1\otimes_{\al_0}\pp_{\al}\log\rh_{\al_0}\right)({\bf V}_{2}+{\bf K}_{2})\cM_{\rh_{\al_0}}=0.\]
Thus, according to Lemma \ref{con-nnPh}, 
\beg{align*}
&P(\al_0)\pp_\al\nn \Ph(0,\al_0) \pi_{\al_0}P(\al_0)\\
&= P(\al_0)\pi_{\al_0}({\bf V}_{2,\al_0}+{\bf K}_{2,\al_0})(I+\al_0\cM_{\pp_{\al}\log\rh_{\al_0}})\pi_{\al_0}P(\al_0)\\
&\quad -\al_0P(\al_0)\left(\1\otimes_{\al_0}\pp_{\al}\log\rh_{\al_0}\right)({\bf V}_{2}+{\bf K}_{2})\cM_{\rh_{\al_0}}P(\al_0)\\
&=P(\al_0)\pi_{\al_0}({\bf V}_{2,\al_0}+{\bf K}_{2,\al_0})\pi_{\al_0}P(\al_0)\\
&\quad +\al_0 P(\al_0)\pi_{\al_0}({\bf V}_{2,\al_0}+{\bf K}_{2,\al_0})\pi_{\al_0} \cM_{\pp_{\al}\log\rh_{\al_0}}\pi_{\al_0}P(\al_0)\\
&\quad +\al_0 P(\al_0)\pi_{\al_0}({\bf V}_{2,\al_0}+{\bf K}_{2,\al_0})(I-\pi_{\al_0}) \cM_{\pp_{\al}\log\rh_{\al_0}}\pi_{\al_0}P(\al_0)\\
&=P(\al_0)\pi_{\al_0}({\bf V}_{2,\al_0}+{\bf K}_{2,\al_0})\pi_{\al_0}P(\al_0)\\
&\quad +\al_0 P(\al_0)\pi_{\al_0}({\bf V}_{2,\al_0}+{\bf K}_{2,\al_0})\pi_{\al_0}  \cM_{\pp_{\al}\log\rh_{\al_0}}\pi_{\al_0}P(\al_0)\\
&\quad +\al_0 P(\al_0)\pi_{\al_0}({\bf V}_{2,\al_0}+{\bf K}_{2,\al_0})\left(\1\otimes_{\al_0}  \pp_{\al}\log\rh_{\al_0}\right)\pi_{\al_0}P(\al_0)
\end{align*}
where in the last equality, we have used
\beg{align*}
(I-\pi_{\al_0}) \cM_{\pp_{\al}\log\rh_{\al_0}}f=\mu_{\al_0}\left(f\pp_\al \log\rh_{\al_0}\right)=\left(\1\otimes_{\al_0}\pp_\al \log\rh_{\al_0}\right) f,~f\in L^2(\mu_{\al_0}).
\end{align*}
Due to Lemma \ref{Ker-nnPh}, for every  $f\in P(\al_0)L^2(\mu_{\al_0})$, there is $\pi_{\al_0}f=f$. Thus $\pi_{\al_0}P(\al_0)=P(\al_0)$. Then we find that
\beg{align*}
P(\al_0)A_1P(\al_0)&=P(\al_0)\pp_\al\nn\Ph(0,\al_0) \pi_{\al_0}P(\al_0) \\
&\quad -\al_0P(\al_0)\pi_{\al_0}({\bf V}_{2,\al_0}+{\bf K}_{2,\al_0}) \left(\1\otimes_{\al_0}\pp_{\al}\log\rh_{\al_0}\right)P(\al_0)\\
&=P(\al_0)\pi_{\al_0}({\bf V}_{2,\al_0}+{\bf K}_{2,\al_0}) \pi_{\al_0}P(\al_0)\\
&\quad +\al_0 P(\al_0)\pi_{\al_0}({\bf V}_{2,\al_0}+{\bf K}_{2,\al_0}) \pi_{\al_0} \cM_{\pp_{\al}\log\rh_{\al_0}}\pi_{\al_0}P(\al_0)\\
&=P(\al_0)\pi_{\al_0}({\bf V}_{2,\al_0}+{\bf K}_{2,\al_0}) \pi_{\al_0}P(\al_0)\\
&\quad +\al_0 P(\al_0)\pi_{\al_0}({\bf V}_{2,\al_0}+{\bf K}_{2,\al_0}) \pi_{\al_0} P(\al_0) \cM_{\pp_{\al}\log\rh_{\al_0}} P(\al_0)\\
&=P(\al_0)\pi_{\al_0}({\bf V}_{2,\al_0}+{\bf K}_{2,\al_0}) \pi_{\al_0}P(\al_0)\left(I+\al_0\cM_{\pp_{\al}\log\rh_{\al_0}}\right) P(\al_0)\\
&=-\ff 1 {\al_0} P(\al_0)A_0P(\al_0)\left(I+\al_0\cM_{\pp_{\al}\log\rh_{\al_0}}\right) P(\al_0).
\end{align*}
Hence,
\beg{align*}
&P(\al_0)\nn \Ph(0,\al_0)P(\al_0)+(\al-\al_0)P(\al_0)A_1P(\al_0)\\
&\qquad =P(\al_0)(I-A_0)P(\al_0)-\ff {\al-\al_0} {\al_0} P(\al_0)A_0P(\al_0)\left(I+\al_0\cM_{\pp_{\al}\log\rh_{\al_0}}\right) P(\al_0).
\end{align*}
Combining this with \cite[DEFINITION, THEOREM 9.2]{Simon} and $(I_{\cH_0}+\tld M_0)$ is invertible on $\cH_0$, we arrive at  
\beg{align*}
&{\rm det}_{2,\cH_0}\left(P(\al_0)\nn \Ph(0,\al_0)P(\al_0)+P(\al_0)A_1P(\al_0)(\al-\al_0)\right)\\
&\qquad = e^{\trac{\left(\tld A_0+\ff {\al-\al_0} {\al_0}\tld A_0(I_{\cH_0}+\tld M_0)\right)}}{\rm det}_{ \cH_0}\left(I_{\cH_0}-\tld A_0-\ff {\al-\al_0} {\al_0}\tld A_0(I_{\cH_0}+\tld M_0)\right)\\
&\qquad = e^{\ff {\al} {\al_0}\trac{ \left(\tld A_0\right)}+\ff {\al-\al_0} {\al_0}\trac{\left(\tld A_0 \tld M_0 \right)}} {\rm det}_{ \cH_0}\left(\tld A_0(I_{\cH_0}+\tld M_0)\right)\\
&\qquad \quad\times {\rm det}_{\cH_0}\left( (I_{\cH_0}+\tld M_0)^{-1}\left(\tld A_0^{-1}-I_{\cH_0}\right)-\ff {\al-\al_0} {\al_0}\right).
\end{align*}
Since the algebraic multiplicity of the  eigenvalue $0$ of $(I_{\cH_0}+\tld M_0)^{-1}\left(\tld A_0^{-1}-I_{\cH_0}\right)$ is odd, we find that ${\rm det}_{\cH_0}\left( (I_{\cH_0}+\tld M_0)^{-1}\left(\tld A_0^{-1}-I_{\cH_0}\right)-\ff {\al-\al_0} {\al_0}\right)$ changes sign. This implies that
\[{\rm det}_{2,\cH_0}\left(P(\al_0)\nn \Ph(0,\al_0)P(\al_0)+P(\al_0)A_1P(\al_0)(\al-\al_0)+o(\al-\al_0)\right)\]
changes sign  at $\al=\al_0$. Therefore, ${\rm det}_{2}(I+\al \pi_{\al}({\bf V}_{2,\al}+{\bf K}_{2,\al})\pi_{\al })$ changes sign at $\al=\al_0$. 

In particular, if $0$ is a semi-simple eigenvalue of $I+\al_0\pi_{\al_0}({\bf V}_{2,\al_0}+{\bf K}_{2,\al_0})\pi_{\al_0}$, then $\tld A_0^{-1}=I_{\cH_0}$ and the algebraic multiplicity of the  eigenvalue $0$ of $(I_{\cH_0}+\tld M_0)^{-1}\left(\tld A_0^{-1}-I_{\cH_0}\right)$ equals to $k_0={\rm dim}(\cH_0)$, which also equals to the algebraic multiplicity of the  eigenvalue $0$ of $I+\al_0\pi_{\al_0}({\bf V}_{2,\al_0}+{\bf K}_{2,\al_0})\pi_{\al_0}$.

\end{proof}

Finally, we show that $\al_0$ is also a bifurcation point in the sense of Definition \ref{bif-proba}.
\beg{lem}\label{lem-bif-bif}
Suppose assumptions of Theorem \ref{thm-bif} hold. Then $\al_0$ is also a bifurcation point for \eqref{eq-bif} in the sense of Definition \ref{bif-proba}.
\end{lem}
\beg{proof}
According to \eqref{bif-L2}, we let $(\rh_n,\al_n)$ be a sequence such that $\rh_n\not\equiv 0$, $\Ph(\rh_n,\al_n)=0$ and 
\[\lim_{n\ra +\infty}\left(|\al_n-\al_0|+\|\rh_n\|_{L^2(\mu_{\al_0})}\right)=0.\]
Let $g$ be a bounded and continuous function on $\R^d$. Due to Lemma \ref{con-rhal}, there is 
\[\lim_{n\ra+\infty}\bar\mu(g\rh_{\al_n})=\bar\mu(g\rh_{\al_0}).\]
According to \eqref{PHPH}, we have that $\nu_{\al_n}:=(\rh_n+1)\rh_{\al_n}\bar\mu$ is a solution of \eqref{eq-bif} and $\nu_{\al_n}\neq \mu_{\al_n}$. Let $J_{\al_0,\de}$ be the interval in Corollary \ref{cor:mu0mu}. Then
\beg{align*}
\left|\bar\mu\left(g\rh_n\rh_{\al_n}\right)\right|&=\left|\int_{\R^d}g\rh_n\ff {\rh_{\al_n}} {\sq{\rh_{\al_0}}} \sq{\rh_{\al_0}}\d\bar\mu\right|\leq \|g\|_\infty \left\|\ff {\rh_{\al_n}} {\sq{\rh_{\al_0}}}\right\|_\infty\left|\int_{\R^d}|\rh_n| \sq{\rh_{\al_0}}\d\bar\mu\right|\\
&\leq \|g\|_\infty\sup_{\al\in J_{\al_0,\de}}\left\|\ff {\rh_{\al}} {\sq{\rh_{\al_0}}}\right\|_\infty\left(\int_{\R^d}|\rh_n|^2 \rh_{\al_0}\d\bar\mu\right)^{\ff 1 2}.
\end{align*}
Combing this with \eqref{rh2rh} and \eqref{Cal0de}, we have that
\[\lim_{n\ra +\infty}\left|\bar\mu\left(g\rh_n\rh_{\al_n}\right)\right|=0.\]
Hence, 
\[\varlimsup_{n\ra +\infty}\left|\nu_{\al_n}(g)-\mu_{\al_0}(g)\right|\leq\lim_{n\ra +\infty}\left|\bar\mu\left(g\rh_n\rh_{\al_n}\right)\right|+\lim_{n\ra +\infty}|\bar\mu(g\rh_{\al_n})-\bar\mu(g\rh_{\al_0})|=0.\]
Therefore, $(\mu_{\al_0},\al_0)$ is a bifurcation point for \eqref{eq-bif} in the sense of Definition \ref{bif-proba}.

\end{proof}

\subsection{Proof of Corollary \ref{cor-finite}}
Existence of solutions for \eqref{nu-si0} follows from Corollary \ref{exa-singular}. For $\si_0$, we can choose $0<\hat\si<\si_0<\check{\si}$ and $\bar V(x)=C_{\si_0}(1+|x|)^{\ga_1}$ for some $C_{\si_0}>0$  such that (A1) holds.  Let $\al=\ff {\be} {\si^2}$. Then
\beg{align*}
\cT(x;\rh,\al)= \ff {\exp\left\{-\ff {2\al} {\be}V_0(x)- \al\int_{\R^d} (V_2+K_2)(x,y)\rh(y)\bar\mu(\d y) +\bar V(x)\right\}} {\int_\R \exp\left\{- \ff {2\al} {\be} V_0(x)- \al\int_{\R^d} (V_2+K_2)(x,y)\rh(y)\d y+\bar V(x)\right\}\bar\mu(\d x)}.
\end{align*}
Since  
\beg{align*}
\mu_{\al_0}(f\pi_{\al_0} {\bf V}_{2,\al_0}\pi_{\al_0} f)&=\sum_{i,j=1}^{l}J_{ij}\mu_{\al_0}((v_i-\mu_{\al_0}(v_i))f)\mu_{\al_0}(v_j(f-\mu_{\al_0}(f))\\
&=\sum_{i,j=1}^{l}J_{ij}\mu_{\al_0}((v_i-\mu_{\al_0}(v_i))f)\mu_{\al_0}((v_j-\mu_{\al_0}(v_j))f)\\
&\geq 0,~f\in L^2(\mu_{\al_0}),
\end{align*}
we find that $I+\al_0\pi_{\al_0}{\bf V}_{2,\al_0}\pi_{\al_0}$ is invertible on $L^2(\mu_{\al_0})$. Hence, the first assertion of this corollary can follow from Lemma \ref{con-rhal} directly, and we focus on the bifurcation point in the following discussion.

We prove that $0$ is an eigenvalue of $I+\al_0\pi_{\al_0}({\bf V}_{2,\al_0}+{\bf K}_{2,\al_0})\pi_{\al_0}$ with odd algebraic multiplicity.   It is clear that ${\bf K}_{2,\al_0}$ and ${\bf V}_{2,\al_0}$ are self-adjoint operators on $L^2(\mu_{\al_0})$. For all $f\in L^2(\mu_{\al_0})$, since $K_2(x,\cdot)$ is anti-symmetric, $V_2(\cdot,y)$ and $\mu_{\al_0}$ are symmetric, we have that ${\bf V}_{2,\al_0}\pi_\al {\bf K}_{2,\al_0}=0$ and
\beg{align*}
{\bf K}_{2,\al_0}\pi_{\al_0}{\bf V}_{2,\al_0} f&=\int_{\R^d} K_2(x,z)\mu_{\al_0}(\d z)\int_{\R^d}(V_2(z,y)-\mu_{\al_0}(V_2(\cdot,y))f(y)\mu_{\al_0}(\d y)\\
&=\int_{\R^d}\left(\int_{\R^d}K_2(x,z)(V_2(z,y)-\mu_{\al_0}(V_2(\cdot,y))\mu_{\al_0}(\d z)\right)f(y)\mu_{\al_0}(\d y)\\
&=0.
\end{align*}
Let $R_V$ be closure of the range of $\pi_{\al_0} {\bf V}_{2,\al_0}\pi_{\al_0}$ and $R_K$ be closure of the range of $\pi_{\al_0} {\bf K}_{2,\al_0}\pi_{\al_0}$. Then $R_V\perp R_K$, $R_V\subset {\rm Ker}(\pi_{\al_0} {\bf K}_{2,\al_0}\pi_{\al_0})$ and $R_K\subset {\rm Ker}(\pi_{\al_0} {\bf V}_{2,\al_0}\pi_{\al_0})$. Then there is subspace $\cH$ such that $L^2(\mu_{\al_0})=R_V\oplus R_K\oplus \cH$ and 
\[\cH\subset  {\rm Ker}(\pi_\al {\bf V}_{2,\al}\pi_{\al})\cap   {\rm Ker}(\pi_\al{\bf K}_{2,\al}\pi_{\al}).\]
For $0\neq f\in L^2(\mu_{\al_0})$ with $f+\al_0\pi_{\al_0}({\bf V}_{2,{\al_0}}+{\bf K}_{2,{\al_0}})\pi_{\al_0}f=0$,  due to the decomposition of $L^2(\mu_{\al_0})$, there is $f=f_1+f_2$ for some $f_1\in R_V$ and $f_2\in R_K$.  Then
\[\left(f_1+\al_0\pi_{\al_0}{\bf V}_{2,\al_0}\pi_{\al_0} f_1\right)+\left(f_2+\al_0\pi_{\al_0}{\bf K}_{2,\al_0}\pi_{\al_0} f_2\right)=0.\] 
This yields that 
\[\left\{\begin{array}{l} 
f_1+\al_0\pi_{\al_0}{\bf V}_{2,\al_0}\pi_{\al_0} f_1=0,\\
f_2+\al_0\pi_{\al_0}{\bf K}_{2,\al_0}\pi_{\al_0} f_2 =0.
\end{array}\right.\]
Since $I+\al_0\pi_{\al_0}V_{2,\al_0}\pi_{\al_0}$ is invertible, $f_1=0$. Thus, $0$ is an eigenvalue of $I+\al_0\pi_{\al_0}({\bf V}_{2,\al_0}+{\bf K}_{2,\al_0})\pi_{\al_0}$  if and only if $0$ is an eigenvalue of $I+\al_0\pi_{\al_0} {\bf K}_{2,\al_0}\pi_{\al_0}=I+\al_0 {\bf K}_{2,\al_0} $, and
\beg{equation}\label{K-VK-K}
{\rm Ker}(I+\al_0\pi_{\al_0}\left({\bf V}_{2,\al_0}+{\bf K}_{2,\al_0}\right)\pi_{\al_0})={\rm Ker}(I+\al_0 {\bf K}_{2,\al_0}).
\end{equation}
Moreover,
\beg{align*}
I+\al_0\pi_{\al_0}({\bf V}_{2,\al_0}+{\bf K}_{2,\al_0})\pi_{\al_0}&=\left(I+\al_0\pi_{\al_0}{\bf V}_{2,\al_0}\pi_{\al_0}\right)\left(I+\al_0 {\bf K}_{2,\al_0} \right)\\
&=\left(I+\al_0 {\bf K}_{2,\al_0} \right)\left(I+\al_0\pi_{\al_0}{\bf V}_{2,\al_0}\pi_{\al_0}\right).
\end{align*} 
Combining this with $I+\al_0\pi_{\al_0}V_{2,\al_0}\pi_{\al_0}$  is invertible on $L^2(\mu_{\al_0})$ and \cite[Theorem 21.2.6 and Theorem 25.2.2']{Lax}, the algebraic multiplicity of $0$ as eigenvalue of $I+\al_0\pi_{\al_0}({\bf V}_{2,\al_0}+{\bf K}_{2,\al_0})\pi_{\al_0}$ is the same as $0$ as eigenvalue of $I+\al_0 {\bf K}_{2,\al_0} $.  \\
Let $\cH_{K}={\rm span}[k_1,\cdots,k_m]$. Since $R_K\subset\cH_K$, $\cH_{K}$ is an invariant space of $I+\al_0{\bf K}_{2,\al_0}$. Due to that $k_1,\cdots,k_m$ are linearly independent, $I+\al_0{\bf K}_{2,\al_0}$ can be represented under the basis $[k_1,\cdots,k_m]$ as the matrix $I+\al_0GG(\si_0)$.  Then $0$ is an eigenvalue of $I+\al_0{\bf K}_{2,\al_0}$ if and only if $0$ is an eigenvalue of the matrix $I+\al_0GG(\si_0)$, and they have the same algebraic multiplicities.  By (1) in the assumption, we have that $0$ is an eigenvalue of $I+\al_0\pi_{\al_0}({\bf V}_{2,\al_0}+{\bf K}_{2,\al_0})\pi_{\al_0}$ with odd algebraic multiplicity.

Let $P_K$ be the orthogonal projection from $L^2(\mu_{\al_0})$ on to $\cH_K$. To give a representation of $\tld M_0$, we give a representation of $P_{K}\cM_{\pp_{\al}\log\rh_{\al_0}}\Big|_{\cH_K}$.  For all $f\in L^2(\mu_{\al_0})$, there is 
 \[P_Kf= \sum_{i,j=1}^m(G(\si_0)^{-1})_{ij}\mu_{\al_0}(fk_j)k_i.\] 
Since $\{\1\}\cup\{v_1,\cdots,v_l\}$ is linearly independent, $\pi_{\al_0}(v_1),\cdots,\pi_{\al_0}(v_l)$ are also linearly independent. Then $J(\al_0)$ is invertible. Let $\cH_V(\al_0)={\rm span}[\pi_{\al_0}(v_1),\cdots,\pi_{\al_0}(v_l)]$ and $P_V$ be the orthogonal projection from $L^2(\mu_{\al_0})$ on to $\cH_V(\al_0)$. Then  
\begin{align*}
&(I-P_V)\pi_{\al_0}{\bf V}_{2,\al_0}=0,\\ 
P_Vf=&\sum_{i,j=1}^l (J(\al_0)^{-1})_{ij}\mu_{\al_0}(f \pi_{\al_0}v_j)\pi_{\al_0}(v_i),~\mu_{\al_0}(f)=0.
\end{align*}
According to \eqref{pp-logrh0}, $\th'(\al)=\ff 2 {\be}$, $V_1=0$ and the representation of $P_V$, we have that
\begin{align*}
\pp_\al\log\rh_{\al_0}&=- (I+\al_0\pi_{\al_0} {\bf V}_{2,\al_0}\pi_{\al_0})^{-1} \pi_{\al_0}\left(\ff 2 {\be}V_0+{\bf V}_{2,\al_0}\1\right)\\
&=-(I+\al_0\pi_{\al_0} {\bf V}_{2,\al_0}\pi_{\al_0})^{-1} P_V\left(\pi_{\al_0}\left(\ff 2 {\be}V_0+{\bf V}_{2,\al_0}\1\right)\right)\\
&\quad  -(I+\al_0\pi_{\al_0} {\bf V}_{2,\al_0}\pi_{\al_0})^{-1} (I-P_V)\left(\pi_{\al_0}\left(\ff 2 {\be}V_0+{\bf V}_{2,\al_0}\1\right)\right)\\
&=-(I+\al_0\pi_{\al_0} {\bf V}_{2,\al_0}\pi_{\al_0})^{-1}\left( \ff 2 {\be}P_V\pi_{\al_0}V_0+\pi_{\al_0}{\bf V}_{2,\al_0}\1\right)\\
&\quad  -\ff 2 {\be}(I+\al_0\pi_{\al_0} {\bf V}_{2,\al_0}\pi_{\al_0})^{-1} (I-P_V)\pi_{\al_0}V_0\\
&=-\sum_{i=1}^l\left[(I +\al_0JJ(\al_0))^{-1}\left(\ff 2 {\be}w+\tld w\right)\right]_{i}\pi_{\al_0}(v_i)\\
&\quad -\ff 2 {\be}(\pi_{\al_0}V_0-\sum_{i=1}^l w_i\pi_{\al_0}(v_i)),
\end{align*}
where $I$ is an identity operator, which maybe on different space from line to line.  Choose an orthonormal basis $[\tld k_1,\cdots,\tld k_m]$ of $\cH_{K}$ with inner product induced by $L^2(\mu_{\al_0})$.   Let $S\in \R^{m}\otimes\R^m$ such that $k_i=\sum_{j=1}^mS_{ij}\tld k_j$.  Then $G(\si_0)=SS^*$.  By the definition of $M_K(\al_0)$, under the basis $[\tld k_1,\cdots,\tld k_m]$, the operator $P_K\cM_{\pp_{\al}\log\rh_{\al_0}}\Big|_{\cH_K}$ can be  represented   as a matrix on $\cH_K$, saying $S^{-1}M_K(\al_0)(S^*)^{-1}$.

Finally, we verify that $(I_{\cH_0}+\tld M_0)$ is invertible on $\cH_0$, and the corollary is proved by using Theorem \ref {thm-bif}. Since 
\begin{align*}
K_2(x,y)&=\sum_{i,j=1}^mG_{ij}k_i(x)k_j(y)=\sum_{i,j=1}^mG_{ij}\left(\sum_{r=1}^mS_{ir}\tld k_r(x)\right)\left(\sum_{n=1}^mS_{jn}k_n(y)\right)\\
&=\sum_{i,j,r,n=1}^mS_{ir}G_{ij}S_{jn} \tld k_r(x) \tld k_n(y),
\end{align*}
${\bf K}_{2,\al_0}$ can be represented as  the matrix $S^*GS$ under the basis $[\tld k_1,\cdots,\tld k_m]$.  Then $\cH_0$ can be represented as ${\rm Ker}(I+\al_0S^*GS) $, which is a subspace of $\R^m$. It is clear that $P(\al_0)\subset P_K$.  We denote by $P_K(\al_0)$ the representation of $P(\al_0)$ restricted on $\cH_K$. Then 
\beg{equation}\label{PKK}
P_K(\al_0)\R^m= {\rm Ker}(I+\al_0S^*GS),
\end{equation}
and $P(\al_0)\cM_{\pp_{\al}\log\rh_{\al_0}}\Big|_{\cH_K}$ can be represented as $P_{K}(\al_0)S^{-1}M_K(\al_0)(S^*)^{-1}$.  By using these matrices, we prove that \\
$\bullet$ $(I_{\cH_0}+\tld M_0)$ is not invertible on $\cH_0$ if and only if the following system has a solution $\left(\beg{array}{c}
w_1\\w_2
\end{array}\right)\in \R^{2m}$ with $w_2\neq 0$
\begin{align}\label{eq-GG}
&\left(I+\al_0S^*GS\right)w_2=0,\\
&\left(I+\al_0S^{-1}M_K(\al_0)(S^*)^{-1}\right)w_2=\left(I+\al_0S^*G S\right)w_1,\label{eq-GM}
\end{align}
Indeed, $(I_{\cH_0}+\tld M_0)$ is not invertible on $\cH_0$ if and only if  there exists $w_2\in\R^m$ with $w_2\neq 0$ such that \eqref{eq-GG} holds and $P_{K}(\al_0)\left(I+\al_0S^{-1}M_K(\al_0)(S^*)^{-1}\right)w_2=0$. Taking into account $ {\rm Ker}(I+\al_0S^*G S)\perp {\rm Ran}(I+\al_0S^*G S)$ and \eqref{PKK}, we have that 
\[\left(I+\al_0S^{-1}M_K(\al_0)(S^*)^{-1}\right)w_2\in {\rm Ran}(I+\al_0S^*G S).\]
Thus there exists $w_1$ such that \eqref{eq-GM} holds. Conversely, if there exists $(w_1,w_2)$ with $w_2\neq 0$ such that \eqref{eq-GG} and \eqref{eq-GM} hold, then 
\[P_K(\al_0)(I+\al_0S^{-1}M_K(\al_0)(S^*)^{-1})w_2=P_K(\al_0)\left(I+\al_0S^*G S\right)w_1=0.\] 
Thus $(I_{\cH_0}+\tld M_0)$ is not invertible on $\cH_0$. \\
Rewrite \eqref{eq-GG} and \eqref{eq-GM} in the following form:
\begin{align*}
T\left(\begin{array}{c} w_1\\ w_2 \end{array}\right):=\left[\begin{array}{cc}
0 & I+\al_0S^*G S\\
(I+\al_0S^*G S) & -(I+\al_0S^{-1}M_K(\al_0)(S^*)^{-1})
\end{array}
\right]\left(\begin{array}{c} w_1\\ w_2 \end{array}\right) =0.
\end{align*} 
It is clear that  $(w_1,0)$ is a solution of this system if and only if $w_1\in {\rm Ker}(I+\al_0S^*GS)$.  Thus  there exists $(w_1,w_2)$ with $w_2\neq 0$ such that \eqref{eq-GG} and \eqref{eq-GM} holds if and only if 
\begin{equation}\label{dim-T-G}
{\rm dim}{\rm Ker}\left(T\right)>{\rm dim}{\rm Ker}(I+\al_0S^*GS).
\end{equation}
Since 
\begin{align*}
&{\rm dim}{\rm Ker}(I+\al_0S^*G S)=m-{\rm rank}(I+\al_0S^*G S),\\
&{\rm dim}{\rm Ker}\left(T\right)=2m-{\rm dim}{\rm Ran}(T^*)=2m-{\rm rank}(T),
\end{align*}
we have that \eqref{dim-T-G} holds if and only if 
\[m+{\rm rank}(I+\al_0S^*GS)>{\rm rank}(T).\]
Hence, $(I_{\cH_0}+\tld M_0)$ is invertible on $\cH_0$ if and only if 
\begin{equation}\label{mSGST}
m+{\rm rank}(I+\al_0S^*G S)\leq {\rm rank}(T).
\end{equation}
Taking into account that 
\begin{align*}
{\rm rank}(T)&\leq {\rm rank}(I+\al_0S^*G S)\\
&\quad + {\rm rank}([I+\al_0S^{-1}M_K(\al_0)(S^*)^{-1},-(I+\al_0S^*G S)])\\
&\leq {\rm rank}(I+\al_0S^*G S)+m, 
\end{align*}
we find that \eqref{mSGST} holds if and only if  $m+{\rm rank}(I+\al_0S^*G S)={\rm rank}(T)$. By using $SS^*=G(\si_0)$, we have that
\begin{align*}
I+\al_0S^{-1}M_K(\al_0)(S^*)^{-1}&=I+\al_0(S^*)(S^*)^{-1}S^{-1}M_K(\al_0)(S^*)^{-1}\\
&=S^*\left(I+\al_0G(\si_0)^{-1}M_K(\al_0)\right)(S^*)^{-1},\\
I+\al_0S^*GS&=I+\al_0S^*G S(S^*)(S^*)^{-1}\\
&=S^*\left(I+\al_0GG(\si_0)\right)(S^*)^{-1}.
\end{align*}
Then 
\begin{align*}
{\rm rank}(I+\al_0S^*G S)&={\rm rank}(I+\al_0GG(\si_0)),\\
{\rm rank}(T)&={\rm rank}\left(\left[\begin{array}{cc}
(S^*)^{-1} & 0\\
0 & (S^*)^{-1}
\end{array}
\right]T\left[\begin{array}{cc}
S^* &0\\
0 & S^*
\end{array}
\right]\right)\\
&={\rm rank}\left(\left[\begin{array}{cc}
0 & I+\al_0GG(\si_0)\\
I+\al_0GG(\si_0) & -(I+\al_0G(\si_0)^{-1}M_K(\al_0))
\end{array}
\right]\right).
\end{align*}
Therefore, we have that $(I_{\cH_0}+\tld M_0)$ is invertible on $\cH_0$ if and only if \eqref{eq-rak} holds.


\section{Appendix: proofs of auxiliary lemmas}

The following lemma is devoted to the regularity  of $\ps(\mu)$ in Section 2. 
\beg{lem}\label{ef}
If $f\in W^{1,p}_{loc}$ for some $p>d$, then $e^f\in W^{1,p}_{loc}$ and $\nn e^f=e^f\nn f$.
\end{lem}
\beg{proof}
Since for any $\ze\in C_0$, there is $N>0$ such that  $\supp\{\ze\}\subset B_N$. Then $e^{f}\ze =e^{f \ze_{2N}}\ze$ and $f \ze_{2N}\in W^{1,p}$. Hence, we first assume that $f\in W^{1,p}$.  In this case,  there is a sequence $\{f_m\}\subset C^{\infty}_0$ such that 
\beg{align*}
\lim_{m\ra+\infty}\|f_m-f\|_{W^{1,p}}=0.
\end{align*}
Since $p>d$, it follows from the Morrey embedding theorem (\cite[Theorem 9.12]{Bre}) that $W^{1,p}\subset L^\infty$ with continuous injection. Then 
\beg{align}\label{sup-phm}
\|f\|_{\infty}\vee \sup_{m\geq 1}\|f_m\|_{\infty}&\leq C\left(\|f\|_{W^{1,p}}\vee \sup_{m\geq 1}\|f_m\|_{W^{1,p}}\right)<\infty,\\
\lim_{m\ra+\infty}\|f_m-f\|_{\infty}&\leq C\lim_{m\ra+\infty}\|f-f_m\|_{W^{1,p}}=0.\label{lim-supW}
\end{align}
By using the following fundamental inequality 
\beg{align}\label{funi-1}
|e^x-e^y|\leq (|x-y|\we 1)e^{x\vee y},~x,y\in\R,
\end{align}
we have that
\beg{align*}
&\|e^{f_{m }}-e^{f}\|_{L^{p}} + \|e^{f_m}\nn f_{m }-e^{f}\nn f\|_{L^{p}}\\
&  \leq  e^{\|f_{m }\|_{ \infty}\vee \|f\|_\infty} \| f-f_{m} \|_{L^{p}}+\|(e^{f_{m}}-e^{f})\nn f_{m}\|_{L^{p}}+e^{ \|f\|_{ \infty}} \|\nn f_{m }-\nn f\|_{L^{p} }\\
&  \leq   e^{\|f_{m }\|_{ \infty}\vee \|f\|_\infty} \left(\| f-f_m\|_{L^{p}}+ \| f-f_m\|_{\infty}\|\nn f_{m}\|_{L^{p}}\right)+e^{ \|f\|_{ \infty}} \|\nn f_{m }-\nn f\|_{L^{p} }.
\end{align*}
This, together with \eqref{sup-phm} and \eqref{lim-supW}, implies that $e^{f_m}$ converges to $e^f$ in $W^{1,p}$ and $\nn e^f=e^f \nn f$. For $f\in W^{1,p}_{loc}$, $e^{f\ze_{2N}}\in W^{1,p}$. Then $e^f\ze=e^{f\ze_{2N}}\ze\in W^{1,p}$. Hence, $e^f\in W^{1,p}_{loc}$ and $\nn e^f=e^f \nn f$. 

\end{proof}

The following lemma is devoted to the invariant probability measure of $L_\mu$.  It is fundamental and we give the proof for readers' convenient. 
\beg{lem}
Assume {\bf (H)}. Then for each $\mu\in\sP_{W_0}$, $\hat\cT(x,\mu)\bar\mu(\d x)$ is an invariant probability measure of $L_{\mu}$.
\end{lem}
\beg{proof} 
For every $g\in C_0^\infty$,  due to  $V_0\in\cW^{1,p}_{q,\bar\mu}$ and \eqref{nnV2}, $\nn\log(\ps(\mu)e^{-\bar V})\in L^q(\bar\mu)\cap L^p_{loc}$. Then $\<\nn \log(\ps(\mu)e^{-\bar V}),\nn g\>\in L^q(\bar\mu)\cap L^p$. 
Hence,  for all $g\in C^\infty_0(\R^d)$, there is $L_{\mu}g\in L^1(\bar\mu)\cap L^1$.
It follows from the integration by part formula that ($\xi_n$ is defined in Lemma \ref{ef})
\beg{equation}\label{add-in00}
\beg{split}
&\left|\int_{\R^d}\ze_n(x)(L_{\mu}g)(x) \ps(x,\mu)\bar\mu(\d x)\right|\\
&\qquad=\left|\int_{\R^d}\ze_n(x)\div(\ps(x,\mu)e^{-\bar V}\nn g)(x) e^{\bar V(x)}\bar\mu( \d x)\right| \\
&\qquad=\left|\ff 1 n\int_{\R^d}\chi'(|x|/n)\left\<\ff x {|x|},  \nn g(x)\right\>\ps(x,\mu)\bar\mu(\d x)\right|\\
&\qquad\leq \ff 2 n\left|\int_{n\leq |x|\leq 2n} | \nn g(x)|\ps(x,\mu) \bar\mu(\d x)\right|\\
&\qquad\leq \ff {2\|\ps(\mu)\|_\infty} n\|\nn g\|_\infty\left|\int_{n\leq |x|\leq 2n} \bar\mu(  \d x)\right|,
\end{split}
\end{equation}
It follows from the dominated convergence theorem  that 
$$\int_{\R^d} (L_{\mu}g)(x) \ps(x,\mu)\bar\mu(\d x)=\lim_{n\ra +\infty}\int_{\R^d} \ze_n(x)(L_{\mu}g)(x) \ps(x,\mu)\bar\mu(\d x)=0.$$

\end{proof}

\noindent\textbf{Acknowledgements}

\medskip

The  author was supported by the National Natural Science Foundation of China (Grant No. 12371153), and Program for Innovation Research in Central University of Finance and Economics.



\begin{thebibliography}{99}

\bibitem{Bog}
	V.~I.~Bogachev, Measure Theory, Volume I, Springer-Verlag Berlin Heidelberg, 2007.
	
\bibitem{BKR}
	V.~I.~Bogachev, N~.V. Krylov and M. R\"ockner,  Elliptic and parabolic equations for measures,  {\it Russ. Math. Surv.} {\bf 64} (2009), 973--1078

	
\bibitem{Bre}
	H. Brezis, Functional Analysis, Sobolev Spaces and Partial Differential Equations, Springer New York Dordrecht Heidelberg London, 2011	
	
\bibitem{BLPR} R. Buckdahn, J. Li, S. Peng and C. Rainer, {\it Mean-field stochastic differential equations and
associated PDEs}, Ann. Probab. 2 (2017), 824--878.

\bibitem{CMV}
	J.A. Carrillo, R.J. McCann, C. Villani, Kinetic equilibration rates for granular media and related equations: entropy
dissipation and mass transportation estimates, {\it Rev. Mat. Iberoam.} {\bf 19} (2003), 971--1018.

\bibitem{CGPS}
	 J.A. Carrillo, R.S. Gvalani, G.A. Pavliotis, A. Schlichting, Long-Time Behaviour and Phase Transitions for the Mckean-Vlasov Equation on the Torus. {\it Arch. Rational Mech. Anal.} {\bf 235} (2020), 635--690. 
	
\bibitem{ChD}
	L.-P. Chaintron, A. Diez. Propagation of chaos: A review of models, methods and applications. I. Models and methods. {\it Kinetic and Related Models}, (2022) {\bf 15(6)}: 895--1015.  
	
\bibitem{Chen}
	Chen, M.-F., From Markov Chains to Non-Equilibrium Particle Systems, (2nd Ed.), World Scientific, 2004.
	

\bibitem{ChPa}
	Chayes, L., Panferov, V., The McKean-Vlasov equation in finite volume. {\it J. Stat. Phys.} {\bf 138(1–3)} (2010), 351--380.

\bibitem{Daw}
	D. A. Dawson, Critical dynamics and fluctuations for a mean-field model of cooperative behavior. {\it J. Stat. Phys.} {\bf 31(1)} (1983), 29--85.
		
\bibitem{Deim}
	K. Deimling, Nonlinear functional analysis, Springer-Verlag Berlin Heidelberg, 1985. 	
	
\bibitem{DGPS}
	M.G. Delgadino, R.S. Gvalani, G.A. Pavliotis, S. A. Smith, Phase Transitions, Logarithmic Sobolev Inequalities, and Uniform-in-Time Propagation of Chaos for Weakly Interacting Diffusions. {\it Commun. Math. Phys.} {\bf 401} (2023) 275--323.
	
	
\bibitem{DuTu}
	M. H. Duong, J. Tugaut, Stationary solutions of the Vlasov-Fokker-Planck equation: Existence, characterization and phase-transition, {\it Applied Mathematics Letters}  {\bf 52} (2016), 38--45
	


	
\bibitem{FenZ}
	Feng, S. and Zheng, X. G., Solutions of a class of non-linear Master equations, {\it Stoch. Proc. Appl.} {\bf 43} (1992), 65--84. 




\bibitem{GLZ08}
	F. Gesztesy, Y. Latushkin, K. Zumbrun, Derivatives of (modified) Fredholm determinants and stability of standing and traveling waves, {\it J. Math. Pures Appl.} {\bf 90} (2008), 160--200. 


\bibitem{HT10a}
	S. Herrmann, J. Tugaut, Non-uniqueness of stationary measures for self-stabilizing processes, {\it Stoch. Proc. Appl.} {\bf 120} (2010), 1215--1246.


	
	

\bibitem{Kato}
	T. Kato, Perturbation Theory For Linear Operators, Corrected Printing of the Second Edition, Berlin-New York: Springer-Verlag, 1980
	
\bibitem{Kie}
	H. Kielh\"ofer, Bifurcation Theory: An Introduction with Applications to Partial Differential Equations, Second Edition, New-York, Springer, 2014. 	


\bibitem{Kra}
	M. A. Krasnosel'skii, Topological methods in the theory of nonlinear integral equations, Pergamon, New York, 1964.	

\bibitem{Lax}
	P. D. Lax, Functional Analysis, John Wiley $\&$Sons, 2022.
	
\bibitem{Mcc}
	J.~L. McCauley, Stochastic Calculus and Differential Equations for Physics and Finance,  Cambridge University Press, New York, 2013.
 
\bibitem{McK}
	H.~P.~McKean, Jr.,  Propagation of chaos for a class of non-linear parabolic equations. In: Stochastic Differential Equations (Lecture Series in Differential Equations, Session 7, Catholic Univ.), pp. 41--57. Air Force Office Sci. Res., Arlington, VA, 1967
	
\bibitem {RZ} 
	M. R\"ockner, X. Zhang, Well-posedness of distribution dependent SDEs with singular drifts, {\it Bernoulli} {\bf 27(2)} (2021), 1131--1158.

\bibitem{Simon}
	B. Simon, Trace Ideals and Their Applications, second ed., Mathematical Surveys and Monographs, vol. 120, Amer. Math. Soc., Providence, RI, 2005.

\bibitem{Szn}
	A.-S., Sznitman, Topics in propagation of chaos. In: \'Ecole d'\'Et\'e de Probabilit\'es de Saint-Flour XIX-1989, volume 1464 of Lecture Notes in Math., pp. 165--251. Springer, Berlin, 1991


\bibitem{Tam}
	Y. Tamura, On asymptotic behaviors of the solution of a nonlinear diffusion equation. {\it J. Fac. Sci. Univ. Tokyo Sect. IA Math.} {\bf 31(1)} (1984), 195--221. 
	
	
\bibitem{Tug10}	
	J. Tugaut, Convergence to the equilibria for self-stabilizing processes in double well landscape, {\it Ann. Probab.} {\bf 41} (2010), 1427--1460.  

\bibitem{Tug13}
	J. Tugaut, Self-stabilizing processes in multi-wells landscape in  $\R^d$-convergence, {\it Stoch. Proc. Appl.} {\bf 123} (2013), 1780--1801. 

\bibitem{Tug14a}
	J. Tugaut, Phase transitions of McKean-Vlasov processes in double-wells landscape, {\it Stochastics}, {\bf 86 (2)} (2014), pp. 257--284

\bibitem{Tug14b}	
	J. Tugaut, Self-stabilizing processes in multi-wells landscape in $\R^d$-Invariant probabilities, {\it J. Theoret. Probab.} {\bf 27 (1)} (2014), pp. 57-79
		
	
\bibitem{Wan18}
	F.-Y. Wang, Distribution dependent SDEs for Landau type equations, {\it Stoch. Proc. Appl.} {\bf 128} (2018), 595--621.
	
\bibitem{Wan23}
F.-Y. Wang, Exponential ergodicity for non-dissipative McKean-Vlasov SDEs, {\it Bernoulli} {\bf 29(2)} (2023), 1035--1062.

    
\bibitem{ZSQ}
	S.-Q. Zhang, Existence and non-uniqueness of stationary distributions for distribution dependent SDEs, {\it Electron. J. Probab.} {\bf 28(93)}  (2023),1--34.
 
  
\end{thebibliography}
\end{document}